\documentclass[11pt]{article}

\usepackage{hyperref}
\usepackage{bm}
\usepackage{booktabs}
\usepackage[margin=15pt,font=small,labelfont={bf,sf},justification=justified]{caption}
\usepackage[superscript,biblabel]{cite}
\usepackage{color}
\usepackage{float}
\usepackage{placeins}
\usepackage{graphicx}
\usepackage{lscape}
\usepackage[bbgreekl]{mathbbol}
\usepackage{amsmath}
\usepackage{mathtools}
\usepackage{multibib}
\usepackage{multirow}
\usepackage{tabularx}
\usepackage{titleref}
\usepackage{url}
\usepackage{siunitx}
\usepackage{subcaption}
\usepackage[export]{adjustbox}
\usepackage{epsfig}

\newif\ifbembo
\newif\ifcharter
\newif\iferewhon
\newif\iflibertine
\newif\iflibertinealt
\newif\ifpalantino
\newif\iftimesnewroman

\bembofalse
\chartertrue
\erewhonfalse
\libertinefalse
\libertinealtfalse
\palantinofalse
\timesnewromanfalse
\ifbembo
  \usepackage[p,osf]{ETbb} % osf in text, tabular lining figures in math
  \usepackage[scaled=.95,type1]{cabin} % sans serif in style of Gill Sans
  \usepackage[varqu,varl]{zi4}% inconsolata typewriter
  \usepackage[T1]{fontenc}
  \usepackage[libertine,vvarbb]{newtxmath}
  \usepackage[cal=boondoxo,bb=boondox]{mathalfa}
\fi

\ifcharter
  \usepackage[scaled=.98,p]{XCharter}
  \usepackage[scaled=1.04,varqu,varl]{inconsolata}
  \usepackage[type1]{cabin}
  \usepackage[uprightscript,xcharter,vvarbb,scaled=1.05]{newtxmath}
  \usepackage[cal=euler,scr=boondoxo,bb=boondox]{mathalpha}  % prefer bb=bboldxLight when available
\fi

\iferewhon
  \usepackage[p,osf,scaled=.98,space]{erewhon} % scaling by .98 not really necessary
  \usepackage[varqu,varl]{inconsolata} % typewriter
  \usepackage[type1,scaled=.95]{cabin} % sans serif like Gill Sans
  \usepackage[utopia,vvarbb]{newtxmath}
\fi

\iflibertine
  \usepackage[sb]{libertinus}
  \usepackage[T1]{fontenc}
  \usepackage{textcomp}
  \usepackage[varqu,varl]{zi4}% inconsolata for mono, not LibertineMono
  \usepackage{libertinust1math} % slanted integrals, by default
  \usepackage[scr=boondoxo,bb=boondox]{mathalpha} %Omit bb=boondox for default libertinus bb
\fi

\iflibertinealt
  \usepackage[sb]{libertinus}
  \usepackage[T1]{fontenc}
  \usepackage{textcomp}
  \usepackage{libertinust1math} % slanted integrals, by default
  \usepackage[scr=boondoxo,bb=boondox]{mathalpha} %Omit bb=boondox for default libertinus bb
\fi

\ifpalantino
  \usepackage[largesc]{newpxtext}
  \usepackage[vvarbb]{newpxmath}
\fi

\iftimesnewroman
  \usepackage{newtxtext} %T1 is default encoding
  \usepackage[scaled=0.95]{inconsolata} % typewriter
  \usepackage[vvarbb]{newtxmath} % vvarbb gives STIX Bbb
\fi

\usepackage{microtype}
\usepackage{titling}
\usepackage{authblk}
\pretitle{\begin{center}\bfseries\Large}
\posttitle{\end{center}}

\predate{\begin{center}\normalsize}
\postdate{\end{center}}

\usepackage{titling}
\usepackage{authblk} % note: titling needs to be loaded before authblk
\pretitle{\begin{center}\bfseries\sffamily\LARGE}
\posttitle{\end{center}}
\usepackage{abstract}
\allowdisplaybreaks
\usepackage{sectsty} % note: the main alternative package, titlesec, does not work with titleref
\allsectionsfont{\sffamily}

\usepackage[left=0.5in,right=0.5in,top=0.5in,bottom=0.5in,includefoot,heightrounded]{geometry}

\makeatletter
\patchcmd{\LS@rot}{90}{-90}{}{}
\patchcmd{\endlandscape}{90}{-90}{}{}
\makeatother

\title{Composite B-Spline Regularized Delta Functions for the Immersed Boundary Method: Divergence-Free Interpolation and Gradient-Preserving Force Spreading}
\author[1]{Cole Gruninger}
\author[1--7,*]{Boyce E. Griffith}
\affil[1]{Department of Mathematics, University North Carolina, Chapel Hill, NC, USA}
\affil[2]{Department of Applied Physical Sciences, University of North Carolina, Chapel Hill, NC, USA}
\affil[4]{Department of Biomedical Engineering, University of North Carolina, Chapel Hill, NC, USA}
\affil[5]{Carolina Center for Interdisciplinary Applied Mathematics, University of North Carolina, Chapel Hill, NC, USA}
\affil[6]{Computational Medicine Program, University of North Carolina School of Medicine, Chapel Hill, NC, USA}
\affil[7]{McAllister Heart Institute, University of North Carolina School of Medicine, Chapel Hill, NC, USA\vspace{0.5\baselineskip}}
\affil[*]{To whom correspondence should be addressed; email: \texttt{boyceg@email.unc.edu}}

\newcommand{\eulVel}{\textbf{u}}

\newcommand{\lagVel}{\textbf{U}}
\newcommand{\eulForce}{\textbf{f}}
\newcommand{\lagForce}{\textbf{F}}

\newcommand{\xx}{\mathbf{x}}
\newcommand{\XX}{\mathbf{X}}
\newcommand{\XXt}{\mathbf{X}_{\text{tracer}}}
\newcommand{\parens}[1]{\mathopen{}\left(#1\right)\mathclose{}}

\newcommand{\IBfour}{$\text{IB}_4$}
\newcommand{\IBsix}{$\text{IB}_6$}
\newcommand{\BSone}{$\text{BS}_1$}
\newcommand{\BStwo}{$\text{BS}_2$}
\newcommand{\BSthree}{$\text{BS}_3$}
\newcommand{\BSfour}{$\text{BS}_4$}
\newcommand{\BSfive}{$\text{BS}_5$}
\newcommand{\BSsix}{$\text{BS}_6$}

\newtheorem{thm}{Theorem}
\newtheorem{cor}{Corollary}
\newenvironment{proof}{\textbf{Proof.}}{}

\begin{document}
\maketitle
\begin{abstract}
  This paper presents an approach to enhance volume conservation in the immersed boundary (IB) method by employing regularized delta functions derived from composite B-splines. These delta functions are constructed using tensor product kernels, similar to the conventional IB method. However, the kernels are B-splines whose polynomial degree varies according to the normal and tangential directions of each velocity component. The conventional IB method, while effective for fluid-structure interaction applications, has long been challenged by poor volume conservation, particularly evident in simulations of pressurized, closed membranes. We demonstrate that composite B-spline regularized delta functions significantly enhance volume conservation through two complementary properties: they provide continuously divergence-free velocity interpolants and maintain the gradient character of forces corresponding to mean pressure jumps across interfaces. By correctly representing these forces as discrete gradients, they eliminate a key source of spurious flows that typically plague immersed boundary computations. Our approach maintains the local nature of the classical IB method, avoiding the computational overhead associated with the non-local Divergence-Free Immersed Boundary (DFIB) method's construction of an explicit velocity potential which requires additional Poisson solves for interpolation and force spreading operations. Through a series of numerical experiments, we show that sufficiently regular composite B-spline kernels can maintain initial volumes to within machine precision. We provide a detailed analysis of the relationship between kernel regularity and the accuracy of force spreading and velocity interpolation operations. Our findings indicate that composite B-splines of at least $C^1$ regularity produce results comparable to the DFIB method in dynamic simulations, with errors in volume conservation primarily dominated by truncation error of the employed time-stepping scheme. This work offers a computationally efficient alternative for improving volume conservation in IB methods, particularly beneficial for large-scale, three-dimensional simulations. The proposed approach requires only substituting the identity of the regularized delta function in an existing IB code, making it an accessible improvement for a wide range of applications in computational fluid dynamics and fluid-structure interaction.
\end{abstract}
\section{Introduction}\label{sec:introduction}
The immersed boundary (IB) method is a mathematical formulation and numerical discretization procedure for modeling systems that involve fluid-structure-interaction (FSI)\cite{peskin2002}. The IB framework combines a Lagrangian description of the immersed structure with an Eulerian description of the fluid. In the continuous formulation of the IB method, the coupling between the two representations is mediated by convolutions with singular delta functions kernels supported along the immersed structure. The discrete formulation of the IB method replaces these singular delta functions by regularized delta functions, and integrals are discretized using a quadrature scheme. This discretization results in two key operations: the \textit{spreading} of Lagrangian force densities from the moving structure to the background fluid grid, and the \textit{interpolation} of fluid velocities from the fluid grid back to the Lagrangian mesh, ensuring the structure moves in concert with the local fluid flow.\par

The IB method has proven highly effective for modeling FSI problems that challenge traditional body-fitted mesh approaches. Since its inception, the method has been applied to a wide range of applications: modeling cardiac mechanics \cite{griffith_luo_2009,hasan2017,griffith_heart2012,chen2016,crowl2011,kaiser2019}; simulating platelet adhesion and aggregation \cite{skorczewski2014,fogelson2008}; studying insect flight \cite{jones2015,Santhanakrishnan2018}; and investigating undulatory swimming \cite{alben2013,bhalla2014,tytell2014,bale2015,hoover2017,nangia2017}. Despite these successes, the IB method has been plagued by poor volume conservation, a limitation that is particularly evident in simulations of pressurized, closed membranes. This issue manifests as a gradual loss of volume over time, a problem that Peskin noted early on when simulating the heart's contraction\cite{peskin1993}.\par 
An elementary consequence of the Reynolds transport theorem is that any closed surface moving with an incompressible fluid must maintain constant volume. Peskin and Printz recognized that the main issue with Peskin's original formulation of the IB method was that, although the fluid velocity field may be \textit{discretely} divergence free, the interpolated velocity is generally not \textit{continuously} divergence free. To address this issue, Peskin and Printz introduced a modified finite-difference approximation to the Eulerian divergence operator that ensures that the interpolated velocity field is continuously divergence free at least in an average sense\cite{peskin1993}. This modified divergence stencil, designed for collocated fluid discretizations, dramatically reduced the volume conservation error of the original method. However, despite the improvements made in volume conservation, the method has not been widely used in the community. One reason could be is the complexity of the finite difference stencil associated with the divergence operator, which must be derived specifically for each regularized delta function employed. Furthermore, Griffith demonstrated that the volume conservation improvements of the Peskin and Printz method were quantitatively similar to those achieved by a standard staggered-grid spatial discretization of the fluid variables \cite{griffith2012}. \par 
Since Peskin and Printz's work, there have been many other efforts to improve the volume conservation properties of the IB method. Cortez and Minion introduced the blob projection method, which solves for a velocity field in which the spread force is projected onto the space of divergence-free vector fields \cite{cortez2000}. They utilized regularized ``blob'' functions to spread forces from the immersed structure, allowing the true projection to be computed analytically. However, their tests demonstrated that the improvement in volume conservation was comparable to Peskin and Printz's method. Lee and Leveque utilized ideas from the immersed-interface method \cite{li1994,li2001,li1997} to improve the volume conservation of the IB method by incorporating the correct pressure jump in their projection method's Poisson solver \cite{long2003}. This approach generated volume conservation errors that converged to zero at a second-order accurate in space. However, their modified IIM implementation is more complicated than the standard IB method. It requires decomposing the Lagrangian force density into its tangential and normal components and computing correction terms for each pressure Poisson solve.
\par 
More recently, Bao et al. introduced an IB method that yields an interpolated velocity field that is continuously divergence-free\cite{bao2017}. We refer to Bao et al.'s method as the non-local Divergence Free Immersed Boundary (DFIB) method. The DFIB method is non-local because it achieves a continuously divergence-free interpolated velocity field by first solving a Poisson problem for a discrete velocity potential. This discrete velocity potential is then interpolated to create a continuous velocity potential, which ultimately yields a continuously divergence-free velocity provided that the regularized delta function used for interpolation is at least $C^2$. 

The DFIB method defines the force spreading operator as the discrete adjoint of the interpolation operator, thereby ensuring energy conservation in Lagrangian-Eulerian interactions \cite{bao2017,peskin2002}. For force spreading, the method determines the Eulerian force by assuming it is discretely divergence-free and constructing it by solving a vector Poisson equation. This approach dramatically improves volume conservation of closed, pressurized membranes, even conserving the initial area to near machine precision\cite{bao2017}.

However, the DFIB method's computational cost increases substantially due to the additional Poisson solves required for both velocity interpolation and force spreading. In two spatial dimensions, one extra Poisson solve is needed for each interpolation operation, and two additional Poisson solves are required for each force spreading operation. In three dimensions, three extra Poisson solves are necessary for both velocity interpolation and force spreading operations. Moreover, the DFIB method has been restricted to uniform periodic Cartesian grids, further limiting its applicability. Additionally, because the resulting Eulerian force density is constructed to be discretely divergence-free, compressive forces from the immersed structure are not spread to the background fluid grid, resulting in a physical pressure field that does not reflect the structure's presence and complicating extraction of physical pressure related to fluid-structure interaction.

Similar in spirit to the DFIB method is the Curl-Flow method developed by Chang et al.~\cite{chang2022}, though it is designed primarily for particle advection rather than coupled fluid-structure interaction models. Curl-Flow computes an analytically divergence-free velocity field using an approach similar to the DFIB method, but employs a parallel sweeping strategy combined with a single scalar Poisson solve to compute the vector potential. Unlike the DFIB method, Curl-Flow can be adapted to a variety of physical boundary conditions as well as non-uniform grids, including static boundaries such as cut cells.

We also mention the Fourier Spectral Immersed Boundary (FSIB) Method developed recently by Chen and Peskin~\cite{chen2024}. The FSIB method is grid-free, leveraging only the Fourier representation of the velocity field with discretization achieved by truncating the number of Fourier modes employed. At every timestep, each Fourier mode associated with the velocity field satisfies $\hat{\eulVel}(\hat{\mathbf{k}})\cdot\hat{\mathbf{k}} = 0$, ensuring that the interpolated velocity field is continuously divergence-free. The FSIB method can alternatively be derived by using a sinc function as the regularized delta function in the classic immersed boundary method on a periodic grid, with the Non-Uniform Fast Fourier Transform (NUFFT) employed to improve efficiency of force spreading and velocity interpolation operations. However, because the method leverages Fourier representations, it is limited to periodic boundary conditions, and an aliasing procedure must be chosen to accurately compute the nonlinear convective term and avoid aliasing errors. Despite these limitations, the FSIB method is more accurate than the standard IB method and exhibits faster convergence rates, particularly for problems such as an elastic ellipsoid immersed in periodic 3D fluid, where solutions have been observed to converge with second-order accuracy in the $L^2$ grid norm. While heuristic arguments suggest the classic IB method should compute a velocity converging at a $3/2$ rate rather than second order in the $L^2$ grid norm~\cite{mori2008}, the reason for this higher-order convergence rate observed with the FSIB method has yet to be explained theoretically.

This work introduces a different approach towards mitigating volume conservation errors associated with the IB method by adopting regularized delta functions constructed using composite B-splines. These regularized delta functions not only provide divergence-free velocity interpolation but also maintain gradient structure when spreading forces that represent pressure jumps across the interface. The consistent transfer of continuous gradients to discrete ones eliminates a key source of spurious flows that typically plague immersed boundary computations. Composite B-splines are in a sense smooth generalizations of the Raviart-Thomas elements\cite{raviart1977} and they have been utilized by the Isogeometric Analysis community to implement divergence-free conforming discretizations of equations modeling incompressible flow and FSI \cite{evans2013,tong2022,casquero2018,casquero2021}. Handscomb appears to be one of the earliest to use composite B-splines to produce divergence-free interpolants of discretely divergence-free velocity fields \cite{handscomb1984}. More recently, Schroeder et al. extended these ideas to the context of Eulerian variables discretized on a MAC grid \cite{schroeder2022}. They developed general divergence-free interpolation schemes based on generalized properties of composite B-splines, producing interpolants that are not only continuously divergence-free but also capable of reproducing discrete velocity fields defined at the edge centers of the MAC grid. Building on this work, Schroeder and colleagues have further introduced continuously curl-free interpolants and divergence-free interpolants adapted to various finite difference stencils of the divergence and curl operators on the MAC grid \cite{chowdhury2024}. 
Inspired by the work of Handscomb and Schroeder et al., we employ composite B-spline regularized delta functions and demonstrate greatly enhanced volume conservation properties of the immersed boundary (IB) method. We find that this enhancement stems from two primary factors.
First, supporting the findings of Handscomb and Schroeder et al., the Lagrangian marker points used to discretize the immersed boundary are advected with a continuously divergence-free velocity field. Second, and what we believe we have identified for the first time, these composite B-spline regularized delta functions regularize distributional force densities that correspond to distributional gradients into discrete gradients on the MAC grid. Consequently, when a force spread on the grid corresponds to a distributional gradient, the composite B-spline regularized delta functions ensure that this force is spread as a discrete gradient, allowing it to be handled entirely by the pressure field rather than incorporating spurious viscous components. This proper treatment enables the IB method utilizing composite B-spline regularized delta functions to attain hydrostatic equilibrium of pressurized membranes, something the IB method employed with isotropic regularized delta functions struggles to achieve.
The use of composite B-spline regularized delta functions is completely local and competitive with, and in some cases more accurate than, the DFIB method, without requiring Poisson solves for force spreading and velocity interpolation operations. Importantly, our approach only substitutes the identity of the regularized delta function used to mitigate volume errors; this is the only modification required. The composite B-spline regularized delta functions are compactly supported and thus efficient to use like Peskin's classic IB method. Furthermore, for quasi-static problems, we demonstrate that composite B-spline regularized delta functions can produce pointwise accurate Lagrangian force densities that converge solely under Lagrangian grid refinement, similar to the DFIB method.

\section{Continuous Equations of Motion}\label{sec:equations}
The immersed boundary method models fluid-structure interaction between a thin (co-dimension one) elastic structure and a surrounding fluid. The structure's motion is described in Lagrangian coordinates, while the fluid is described in Eulerian coordinates on a Cartesian grid. Let $\XX(s,t)$ denote the Cartesian position at time $t$ of a material point identified by the curvilinear coordinate $s$. We denote by $\Gamma$ the curve traced out by these material points, which we assume to be closed and continuously differentiable. We will assume throughout the rest of the paper that the curvilinear coordinate $s$ belongs to the interval $[0,2\pi)$. In the fluid domain $\Omega$, we define the Eulerian velocity field $\eulVel(\xx,t)$ and pressure field $p(\xx,t)$. The coupling between structure and fluid is mediated through interaction equations involving the Dirac delta function. For a viscous, incompressible fluid with constant density and a massless immersed structure, the equations of motion take the form:
\begin{align}
\label{eq:mom}
\rho\frac{D\eulVel}{Dt} &= \mu\Delta\eulVel(\xx,t) - \nabla p(\xx,t) + \eulForce(\xx,t), & \xx\in\Omega \\
\label{eq:div_free}
\nabla\cdot\eulVel\parens{\xx,t} &= 0, & \xx\in\Omega \\ 
\label{eq:vel_interp}
\frac{\partial \XX(s,t)}{\partial t} = \eulVel(\XX(s,t),t) &= \int_{\Omega} \eulVel(\xx,t)\,\delta(\xx - \XX(s,t))\thinspace\text{d}\xx, & s\in I\\
\label{eq:spread_cont}
\eulForce(\xx,t) &= \int_{0}^{2\pi}\lagForce(s,t)\,\delta(\xx-\XX(s,t))\thinspace\text{d}s. & \xx\in\Omega
\end{align}
Equations \eqref{eq:mom} and \eqref{eq:div_free} represent the Navier-Stokes equations for an incompressible, viscous fluid characterized by its constant density $\rho$ and viscosity $\mu$. The left-hand side of equation \eqref{eq:mom} contains the material derivative, $\frac{D\eulVel}{Dt}$, which describes the total rate of change of the velocity field. It is defined as $\frac{D\eulVel}{Dt} = \frac{\partial \eulVel}{\partial t} + \left(\eulVel\cdot\nabla\right)\eulVel$. On the right-hand side of equation \eqref{eq:mom}, the forcing term involves an integral transform with a Dirac delta function kernel. This term generates an Eulerian force density equivalent to the Lagrangian force density $\lagForce(s,t)$ defined on the immersed structure. The Lagrangian force densities are determined from the configuration of the interface. Often, the Lagrangian force density is taken to be the negative Fr\'{e}chet derivative of an energy functional $E(\XX(s,t))$, i.e. $\lagForce(s,t)= -\frac{\partial E}{\partial\XX}(s,t)$. The Dirac delta function appears again in equation \eqref{eq:vel_interp}, in which it acts as an integral kernel to ensure that the Lagrangian configuration moves according to the local fluid velocity. We remark that although the equations of motion associated with the IB method have been outlined here in two spatial dimensions with thin, co-dimension one immersed structures, extensions to volumetric bodies and to three spatial dimensions are straightforward.

\section{Discrete Equations of Motion}\label{sec:numerical implementations}

\subsection{Eulerian Spatial Discretization}\label{sec:spatial_disc}
\begin{figure}
    \centering
    \includegraphics[width=0.47\textwidth]{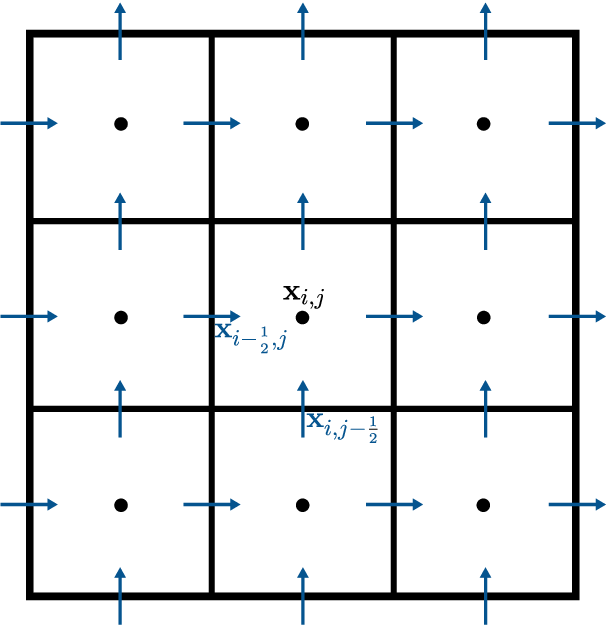}
    \caption{Illustration of the edge-centered and cell center locations defined about the Cartesian cell $\xx_{i,j}$.}
    \label{fig:Mac_grid}
\end{figure}

For simplicity of notation and discussion, we consider the fluid domain $\Omega$ as a periodic square with side length $L$ discretized on a uniform Cartesian grid of size $N\times N$, with grid increments $\Delta x = \Delta y = h = \frac{L}{N}$. The Eulerian variables are discretized using the MAC staggered-grid discretization introduced by Harlow and Welch \cite{harlow1965}.
We emphasize that these discretization choices are made purely for notational convenience and do not represent fundamental limitations of our approach. The IB method can readily accommodate physical boundary conditions rather than periodic ones\cite{griffith2012immersed,kallemov2016,griffith2009simulating}. More importantly, composite B-spline regularized delta functions can be used in any IB context where isotropic regularized delta functions are applicable, without deterioration in performance. This includes anisotropic Cartesian grids and locally adapted Cartesian grids, provided that the support of the regularized delta function does not extend past a coarse-fine interface. 
The MAC discretization approximates the pressure $p_{i,j} = p\left(\xx_{i,j}\right)$ at the centers of Cartesian grid cells $\xx_{i,j} = \left(\left(i + \frac{1}{2}\right)h, \left(j + \frac{1}{2}\right)h\right)$. The discrete velocity $\eulVel_{i,j} = \left(u_{i,j},v_{i,j}\right)$ is defined on the centers of the Cartesian grid cell edges with the $x$-component of the velocity $u_{i,j} = u(\xx_{i-\frac{1}{2},j})$ located at $\xx_{i-\frac{1}{2},j} = \left(ih, (j + \frac{1}{2})h\right)$ and the $y$-component of the velocity $v_{i,j} = v(\xx_{i,j-\frac{1}{2}})$ located at $\xx_{i,j-\frac{1}{2}} = \left(\left(i + \frac{1}{2}\right)h, jh\right)$. For an illustration of the MAC discretization, see Figure \ref{fig:Mac_grid}. Following Bao et al.\cite{bao2017}, we denote the cell-center degrees of freedom using $\mathbb{C}$, and the edge-centered degrees of freedom using $\mathbb{E}$. Additionally, we define the nodal degrees of freedom $\mathbb{N}$ which are located at $\xx_{i-\frac{1}{2},j-\frac{1}{2}} = \left(ih,jh\right)$. Although no discrete Eulerian variables are approximated on $\mathbb{N}$, the discrete curl of the velocity field and the discrete scalar potential needed for the DFIB method are both approximated at the grid nodes. 
We let $\nabla_h$ denote the discrete gradient operator, $\nabla_h\cdot\mbox{}$ the discrete divergence operator, $\nabla_h\times\mbox{}$ the discrete curl operator, and $\Delta_h$ the discrete Laplace operator.\begin{samepage} The application of each discrete differential operator is given by
\begin{align}
    \nabla_h p_{i,j} &= \frac{1}{h}\begin{bmatrix}
        p_{i,j} - p_{i-1,j} \\[1em]
        p_{i,j} - p_{i,j-1}
    \end{bmatrix}, \\
    \nabla_h\cdot\eulVel_{i,j} &= \frac{1}{h}\left(u_{i+1,j} - u_{i,j} + v_{i,j+1} - v_{i,j}\right), \\
    \nabla_h\times\eulVel_{i,j} &= \frac{1}{h}\left(v_{i,j}-v_{i-1,j} + u_{i,j-1}-u_{i,j}\right),\\
    \Delta_h\eulVel_{i,j} &= \nabla_h\cdot\nabla_h\eulVel_{i,j} = \frac{1}{h^2}\left(\eulVel_{i+1,j} + \eulVel_{i,j+1} + \eulVel_{i-1,j} + \eulVel_{i,j-1} - 4 \eulVel_{i,j} \right).
\end{align}
\end{samepage}
The mappings between degrees of freedom for each discrete differential operator are:
\begin{align}
\nabla_h &: \mathbb{C}\to \mathbb{F}, \quad \text{(discrete gradient)} \\
\nabla_h\cdot &: \mathbb{F}\to\mathbb{C}, \quad \text{(discrete divergence)} \\
\nabla_h\times &: \mathbb{F}\to\mathbb{N}. \quad \text{(discrete curl)}
\end{align}
The discrete Laplacian $\Delta_h$ maps a variable to grid locations corresponding to its own degrees of freedom. While we previously described it acting on an edge-centered velocity, it can also be applied to quantities defined at other locations on the Cartesian grid. When applied to a cell-centered quantity, it produces a cell-centered result. Similarly, when used on a quantity defined on the grid nodes, it yields a result also defined on the nodes. In addition the spatial finite difference operators above, we also define the discrete rotation of a scalar field
\begin{align}
    \nabla_h^{\perp} : \mathbb{N}\to\mathbb{F}, \quad \text{(discrete rotation)}
\end{align}
in which $\nabla_h^{\perp}$ is a ninety degree rotation of the discrete gradient operator:
\begin{align}
    \nabla_h^{\perp}a_{i,j} = \frac{1}{h}\begin{bmatrix}
        a_{i,j+1}-a_{i,j} \\
        -\left(a_{i+1,j} - a_{i,j}\right)
    \end{bmatrix}.
\end{align}
\par
The nonlinear convective term $\left(\eulVel\cdot\nabla\right)\eulVel$ is discretized in advective form using simple second-order accurate central finite differences\cite{e2001}. We denote the discrete convective term as $\mathbf{N}(\mathbf{u})_{i,j}$. Since the velocity components are defined on $\mathbb{F}$ and are staggered in space, we need to interpolate the vertical velocity component $v$ to the x-edge centers where the horizontal component $u$ is defined, and interpolate $u$ to the y-edge centers where $v$ is defined. We perform these interpolations using simple averages. \begin{samepage}The discretization of the convective term is given by 
\begin{align}
    \mathbf{N}(\eulVel)_{i,j} = \frac{1}{2h}\begin{bmatrix}
    &u_{i,j}\left(u_{i+1,j} - u_{i-1,j}\right) \\ 
    &+ \frac{1}{4}\left(v_{i-1,j} + v_{i,j} + v_{i-1,j+1} + v_{i,j+1}\right)\left(u_{i,j+1} - u_{i,j-1}\right) \\[1em]
    &\frac{1}{4}\left(u_{i,j-1} + u_{i,j} + u_{i+1,j} + u_{i+1,j-1}\right)\left(v_{i+1,j}-v_{i-1,j}\right) \\ 
    &+ v_{i,j}\left(v_{i,j+1}-v_{i,j-1}\right)
   \end{bmatrix}.
\end{align}
\end{samepage}

\subsection{Mimetic Properties of the Staggered MAC grid}
A fundamental property of the staggered MAC grid, combined with the finite difference approximations of divergence, curl, and gradient operators defined above, is the existence of a discrete Helmholtz decomposition. While we present the proof for periodic boundary conditions and a uniform Cartesian grid, analogous decompositions for anisotropic Cartesian grids and for other boundary conditions so long as the potential functions' boundary conditions handled appropriately\cite{dobrowolski2001,stoyan2004}. We begin by stating key identities proven by Bao et al.\cite{bao2017} for two spatial dimensions (noting that analogous results hold in three dimensions):
\begin{align}
    \nabla_h\times\nabla_h \phi_{i,j} &= 0, \\
    \nabla_h\cdot\left(\nabla_h\times\mathbf{u}_{i,j}\right) &= 0, \\
    \sum_{i,j} \nabla_h\phi_{i,j}\cdot\eulVel_{i,j}\,h^2 &= -\sum_{i,j}\phi_{i,j}\nabla_h\cdot\eulVel_{i,j}h^2, \label{eq:sum_by_parts_div_grad} \\
    \sum_{i,j}\nabla_h^{\perp}a_{i,j}\eulVel_{i,j}h^2 &= \sum_{i,j}a_{i,j}\nabla_h\times\eulVel_{i,j}h^2.
\end{align}

A key result from Bao et al.\cite{bao2017}, originally proven in three dimensions but readily adapted to two dimensions, is:

\begin{thm}\label{thm:theorem1}
    For any zero-mean, discretely divergence-free velocity field $\eulVel_{i,j}$ defined on the faces $\mathbb{F}$ of the MAC grid ($\nabla_h\cdot \eulVel_{i,j} = 0$), there exists a scalar potential $a_{i,j}$ defined on the grid nodes $\mathbb{N}$ such that 
    \begin{align}
        \nabla_h^{\perp}a_{i,j} = \eulVel_{i,j}.
    \end{align}
    Moreover, this potential is unique up to a constant factor.
\end{thm}

The proof relies on the fact that discretely harmonic functions on any portion of the MAC grid (nodes, faces, or cell-centers) must be constant. This theorem leads directly to the discrete Helmholtz decomposition:

\begin{cor}
    Any zero-mean vector field $\eulVel_{i,j}$ defined on the faces $\mathbb{F}$ of the MAC grid admits a discrete Helmholtz decomposition
    \begin{align}
        \eulVel_{i,j} = \nabla_h\phi_{i,j} + \nabla_h^{\perp}a_{i,j},
    \end{align}
    where $\phi_{i,j}$ is defined on $\mathbb{C}$ and $a_{i,j}$ on $\mathbb{N}$. The scalar potentials $a_{i,j}$ and $\phi_{i,j}$ are unique up to constants. 
\end{cor}

\begin{proof}
    We construct $\phi_{i,j}$ as the solution to the discrete Poisson equation
    \[
    \Delta_h\phi_{i,j} = \nabla_h\cdot\eulVel_{i,j}, 
    \]
    which exists because $\nabla_h\cdot\eulVel_{i,j}$ has zero mean (evident from identity \eqref{eq:sum_by_parts_div_grad}). The difference $\eulVel_{i,j} - \nabla_h\phi_{i,j}$ is discretely divergence-free since $\nabla_h\cdot\nabla_h = \Delta_h$ and $\Delta_h\phi_{i,j} = \nabla_h\cdot \eulVel_{i,j}$. Therefore, by Theorem \ref{thm:theorem1}, there exists a scalar potential $a_{i,j}$ on $\mathbb{N}$ satisfying
    \[
    \nabla_h^{\perp}a_{i,j} = \eulVel_{i,j} - \nabla_h\phi_{i,j},
    \]
    which completes the decomposition.
\end{proof}

In three dimensions, this result generalizes with the scalar potential $a_{i,j}$ replaced by a vector potential $\mathbf{a}_{i,j,k}$, and $\nabla_h^{\perp}$ replaced by the discrete curl operator. The uniqueness of $\mathbf{a}_{i,j,k}$ (up to a constant) requires the additional gauge condition $\nabla_h\cdot\mathbf{a}_{i,j,k} = 0$.

\subsection{Lagrangian Spatial Discretization}
To discretize the immersed boundary, we use $M$ Lagrangian marker points positioned at $\XX(s_k,t)$, where $k = 0, 1, ..., M-1$. We uniformly sample the parameter values $s_k$ along the interval $[0, 2\pi)$, with a constant increment $\Delta s = 2\pi/M$ between consecutive values.\par
We choose the number of Lagrangian marker points $M$ and the increment size $\Delta s$ based on a desired value of the ratio $\frac{\Delta X}{h}$. This ratio represents the physical distance between Lagrangian markers at the start of a simulation relative to the increment of the background Cartesian grid. We refer to this ratio as the \textit{mesh factor} and denote it by $M_{\text{fac}} = \frac{\Delta X}{h}$.\par
The Lagrangian force densities are similarly discretized using these same parameter values, so that each marker point has an associated Lagrangian force density $\lagForce_k(t) = \lagForce(s_k,t)$. The discretization of the Lagrangian force densities will be described in each numerical test below. 

The mesh factor is typically chosen based on the dominant flow characteristics~\cite{lee2022,gruninger2024}. For flows with relatively large pressure loads on the boundary, $M_{\text{fac}} \le 1$ is used to ensure the fluid velocity interacts with a sufficiently dense set of Lagrangian marker points, preventing "leaks" from occurring. Conversely, for shear-dominant flows such as flow past a sphere or cylinder, or channel flow modeled using Lagrangian markers, $M_{\text{fac}} > 1$ has been shown to provide more accurate solutions.

When studying convergence properties of the IB method, the mesh factor should be held constant to ensure consistent convergence behavior. We note that the use of composite B-spline regularized delta functions does not alter the recommended values of the mesh factor parameter.

\subsection{Regularized Delta Functions}
In the continuous setting, the immersed boundary method implements fluid-structure interaction through convolutions with singular Dirac delta function kernels as described in equations \eqref{eq:vel_interp} and \eqref{eq:spread_cont}. The presence of the singular Dirac delta function poses a significant numerical challenge. Thus, the first step in constructing the discretized the IB method is to replace this singular function with a regularized version, $\delta_h(\xx)$, in which the regularization parameter is chosen to be identical to the Eulerian meshwidth $h$. These regularized delta functions are typically expressed in tensor product form:
\begin{align}
\label{eq:tens_delta}
\delta_h(\xx) = \frac{1}{h^2}\phi\left(\frac{x}{h}\right)\psi\left(\frac{y}{h}\right),
\end{align}
in which $\phi(r)$ and $\psi(r)$ are one-dimensional kernel functions. Although the formula allows for different functions $\phi$ and $\psi$ in the $x$ and $y$ directions respectively, in the literature, these are almost always chosen to be the same function, yielding an approximately isotropic kernel function. Indeed, to our knowledge, there are no published descriptions of the immersed boundary method in which the functions $\phi$ and $\psi$ differ. 
Following Schroeder et al. \cite{schroeder2022}, who showed that composite B-splines produce continuously divergence-free interpolants of discretely divergence-free velocity fields on the MAC grid, our locally divergence-free IB method uses B-spline kernels for $\phi$ and $\psi$. These kernels differ in polynomial degree by one. In our tests, we also use isotropic kernels where $\phi$ and $\psi$ are identical, allowing us to compare our method with Peskin's original approach. \par 

Each of the one-dimensional kernels we consider in this paper are derived from two different families of kernel functions. The first family we consider is the cardinal B-spline family of kernels, whose invention is attributed to Schoenberg \cite{schoenberg1946}, but was mostly popularized by de Boor \cite{deBoor1978}. B-splines have the attractive property that they can be constructed recursively. This construction starts from the piecewise constant B-spline, $\text{BS}_1 (r)$:
\begin{equation}
\text{BS}_1 (r) = \begin{cases}
1 \quad -\frac{1}{2}\leq r\leq \frac{1}{2}, \\
0 \quad \text{otherwise}.
\end{cases}
\end{equation}
The rest of the B-spline family is then generated using recursive convolution:
\begin{equation}
\label{eq:recur_rel}
\text{BS}_{n+1}(r) = \int_{-\infty}^{\infty}\text{BS}_1(r-q)\text{BS}_n(q)\,\text{d}q.
\end{equation}
From this recursive identity, many properties of the B-spline family may be concluded. For example, taking the derivative of the $\left(n+1\right)^{\text{th}}$ B-spline results in a central difference of the $n^{\text{th}}$ B-spline 
\begin{align}
    \label{eq:deriv_prop}
    \frac{\text{d}}{\text{d}r}\text{BS}_{n+1}(r) = \text{BS}_n\left(r + \frac{1}{2}\right) - \text{BS}_n\left(r - \frac{1}{2}\right).
\end{align}
We remark that equation \eqref{eq:recur_rel} allows us to infer that the $n^{\text{th}}$ B-spline kernel is made up of $n$ nonzero polynomials of degree $n - 1$. The $n^{\text{th}}$ B-spline is $n-2$ times continuously differentiable, is an even function, and is compactly supported with support contained in the interval $\frac{-n}{2} < r < \frac{n}{2}$. One may also use equation $\eqref{eq:recur_rel}$ to show that as $n$ approaches infinity, the sequence of B-splines converges to a rescaled Gaussian \cite{unser1992}. \par
\begin{figure}
    \includegraphics[width=\textwidth]{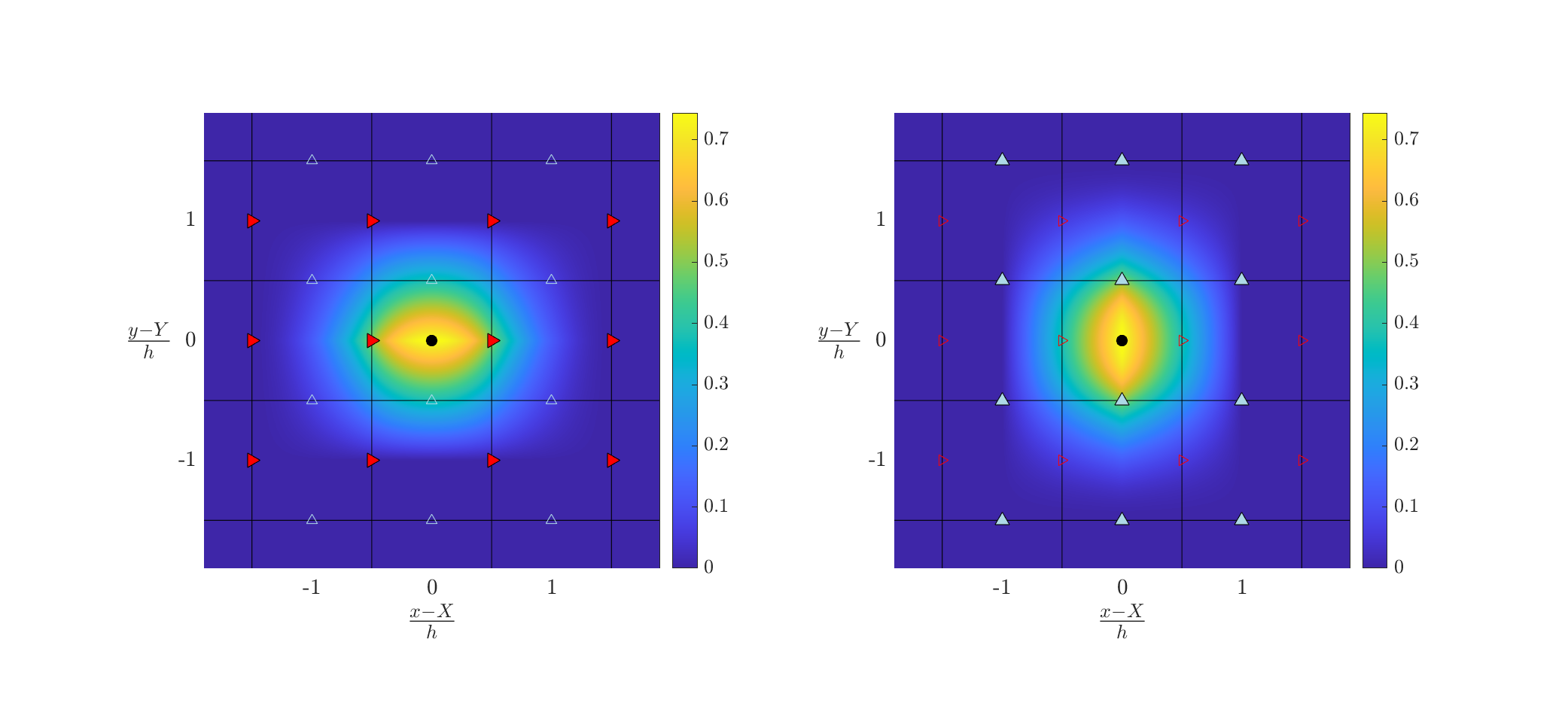}
    \caption{Illustration of the composite B-spline kernel $\mathrm{BS}_3\mathrm{BS}_2$ used for velocity interpolation and force spreading operations. Left: the composite $\mathrm{BS}_3\mathrm{BS}_2$ kernel used to interpolate $x$-components of velocity and spread $x$-components of Lagrangian force density. Right: the composite $\mathrm{BS}_3\mathrm{BS}_2$ kernel used to interpolate $y$-components of velocity and spread $y$-components of Lagrangian force density.}
    \label{fig:comp_bspline_ill}
\end{figure}

 
The regularized delta functions constructed using B-spline kernel functions are composite in nature. To interpolate the $x$-component of the velocity or spread the $x$-component of the force, we set $\phi(x) = \text{BS}_{k+1}(x)$ and $\psi(y) = \mathrm{BS}_{k}(y)$ in equation \eqref{eq:tens_delta}. Similarly, to interpolate the $y$-component of the velocity or spread the $y$-component of the force, we set $\phi(x) =\text{BS}_{k}(x)$ and $\psi(y) = \mathrm{BS}_{k+1}(y)$. For an illustration, see Figure \ref{fig:comp_bspline_ill}. This approach ensures that the interpolated velocities resulting from the discretely divergence-free velocities defined on the MAC grid are continuously divergence-free \cite{schroeder2022}. We detail the discrete velocity interpolation and force spreading operations in the following section. \par 
The second family of one-dimensional kernel functions we consider is the IB family of regularized delta functions. Peskin conceived this family, designing them to be computationally efficient, accurate, and physically relevant \cite{peskin2002}. In this work, we employ two kernels from this family: the four-point IB kernel, denoted as \IBfour, and the relatively new six-point IB kernel, denoted as \IBsix. The \IBfour~kernel is continuously differentiable and serves as our benchmark for comparing the standard IB method's performance against the composite B-spline regularized delta functions. The \IBsix~kernel, which we use to implement the divergence-free interpolation and force-spreading scheme developed by Bao et al. \cite{bao2017}, satisfies all the properties of the \IBfour~kernel but with improved grid translational invariance. Moreover, the \IBsix~kernel is three times continuously differentiable \cite{bao2016}. Regularized delta functions from the IB family are taken to be isotropic. Specifically, we use the same IB kernel function (either \IBfour~or \IBsix) for both one-dimensional kernels $\phi$ and $\psi$ in equation \eqref{eq:tens_delta}. This kernel is applied consistently for both $x$ and $y$ components of the interpolated velocities and spread forces. 

\subsection{Velocity Interpolation and Force Spreading Operations}
This section details the discretization of velocity interpolation and force spreading operations, as described by equations \eqref{eq:vel_interp} and \eqref{eq:spread_cont} in the continuous equations of motion. We present two discretization approaches: the non-local DFIB method introduced by Bao et al. \cite{bao2017} and the standard IB method, introduced by Peskin \cite{peskin2002}. Both methods replace the singular regularized delta function appearing in the continuous equations with a regularized version, as described in the previous section. The DFIB method, however, introduces an additional step: it first solves for a discrete potential function satisfying theorem \ref{thm:theorem1}, which is then used to ensure that the interpolated velocity field is continuously divergence-free. The force spreading operations for each method are constructed to be adjoint to the corresponding velocity interpolation operation so that energy and momentum are discretely conserved \cite{peskin2002}.

\subsubsection{Non-local DFIB Velocity Interpolation and Force Spreading}
We now describe the non-local DFIB velocity interpolation and force spreading operations. Although Bao et al. originally detailed this method in three spatial dimensions \cite{bao2017}, we present it here in two spatial dimensions to align with our test scenarios. Following their approach, we first outline the velocity interpolation operation and then construct the force spreading operation as its discrete adjoint.

Let $\eulVel_{i,j}$ be a discretely divergence-free MAC vector field and let $\eulVel_0$ be its discrete mean flow:
\begin{align}
    \eulVel_0 = \frac{1}{\left|\Omega\right|}\sum_{i,j=0}^{N-1}\eulVel_{i,j}\,h^2,
\end{align}
in which $\left|\Omega\right|$ is the area of $\Omega$. By Theorem \ref{thm:theorem1}, there exists a scalar potential $a_{i,j}$ defined at the nodes $\mathbb{N}$ such that
\begin{align}
    \nabla_h^{\perp} a_{i,j} = \eulVel_{i,j} - \eulVel_0.
\end{align}
This potential is unique up to an additive constant and can be computed by solving the discrete Poisson equation
\begin{align}
    \Delta_h a_{i,j} &= -\nabla_h\times\eulVel_{i,j}.
\end{align}

To construct a continuously divergence-free interpolant, we first interpolate the scalar potential $a_{i,j}$ using a regularized delta function:
\begin{align}
    A(\XX) = \sum_{i,j = 0}^{N-1} a_{i,j}\,\delta_h\left(\xx_{i-\frac{1}{2},j-\frac{1}{2}} - \XX\right)\,h^2.
\end{align}
The interpolated velocity $\lagVel(\XX)$ is then computed by taking the continuous perpendicular gradient and adding back in the mean flow:
\begin{align}
    \label{eq:DFIB_vel_interp}
    \lagVel(\XX) = \mathcal{J}_{\text{DFIB}}\left[\XX\right]\eulVel = \eulVel_0 + \sum_{i,j=0}^{N-1}a_{i,j}\,\nabla^{\perp}\delta_h\left(\xx_{i-\frac{1}{2},j-\frac{1}{2}} - \XX\right)\,h^2.
\end{align} 

When the regularized delta function is at least twice continuously differentiable, Clairut's theorem regarding mixed partial derivatives ensures the interpolated velocity field is continuously divergence-free. In practice, the scalar potential is never explicitly interpolated, instead the continuum perpendicular gradient of the regularized delta function is computed directly at the grid nodes $\mathbb{N}$ and equation \eqref{eq:DFIB_vel_interp} is used directly. 

The DFIB method's force spreading operation is constructed as the discrete adjoint of the velocity interpolation operation. This ensures that energy is conserved between Eulerian and Lagrangian interactions, similar to the standard IB method. Bao et al. \cite{bao2017} provided the derivation of the DFIB force spreading operation in three dimensions. Below, we adapt Bao et al.'s derivation to the two-dimensional case for completeness. \par 
Following Bao et al.'s approach, we obtain the two-dimensional DFIB force spreading operation starting from the principle that it should be the discrete adjoint of the velocity interpolation operation. That is,
\begin{align}
    \label{eq:adjoint_eq}
    \sum_{i,j=0}^{N-1}\eulVel_{i,j}\cdot\eulForce_{i,j}\,h^2 = \sum_{k=0}^{M-1}\lagVel_k\cdot\lagForce_k\,\Delta s,
\end{align}
in which $\lagVel_k = \lagVel(\XX_k)$. By replacing both the Eulerian $\eulVel_{i,j}$ and Lagrangian $\lagVel_k$ velocities with their representations based on the discrete scalar potential $a_{i,j}$, the left-hand side of equation \eqref{eq:adjoint_eq} becomes 
\begin{align}
    \sum_{i,j=0}^{N-1}\eulVel_{i,j}\cdot\eulForce_{i,j}\,h^2 &= \eulVel_0\cdot\eulForce_0\left|\Omega\right| + \sum_{i,j=0}^{N-1}\nabla_h^{\perp}a_{i,j}\cdot \eulForce_{i,j}\,h^2, \\
    &= \eulVel_0\cdot\eulForce_0\left|\Omega\right| + \sum_{i,j=0}^{N-1}a_{i,j}\left(\nabla_h\times\eulForce_{i,j}\right)\,h^2.
\end{align}
In the equation above, we note that we have applied summation by parts and the assumption of periodic boundary conditions to transfer the discrete perpendicular gradient operator acting on $a_{i,j}$ to a discrete curl operator acting on $\eulForce_{i,j}$. We introduce the discrete scalar potential $a_{i,j}$ on the right hand side of equation \eqref{eq:adjoint_eq} by replacing $\lagVel_k$ by the divergence-free interpolated velocity
\begin{align}
    \sum_{k=0}^{M-1}\lagVel_k\cdot\lagForce_k\,\Delta s &= \eulVel_0\cdot\sum_{k=0}^{M-1}\lagForce_k\,\Delta s + \sum_{k=0}^{M-1}\sum_{i,j=0}^{N-1}a_{i,j}\nabla_h^{\perp}\delta_h(\xx_{i-\frac{1}{2},j-\frac{1}{2}} - \XX_k)\cdot\left(\lagForce_k\Delta s\right)\,h^2, \\
    &= \eulVel_0\cdot\sum_{k=0}^{M-1}\lagForce_k\,\Delta s + \sum_{i,j=0}^{N-1}a_{i,j}\sum_{k=0}^{M-1}\left(\nabla_h\delta_h\left(\xx_{i-\frac{1}{2},j-\frac{1}{2}}-\XX_k\right)\times \lagForce_k\Delta s\right)\,h^2.
\end{align}
Therefore, equation \eqref{eq:adjoint_eq} holds if we set 
\begin{align}
    \label{eq:eul_force_mean}
    \eulForce_0 = \frac{1}{\left|\Omega\right|}\sum_{k=0}^{M-1}\lagForce_k\,\Delta s,
\end{align}
and 
\begin{align}
    \label{eq:curl_force_rel}
    \nabla_h\times \eulForce_{i,j} = \sum_{k=0}^{M-1}\left(\nabla\delta_h\left(\xx_{i-\frac{1}{2},j-\frac{1}{2}}-\XX_k\right)\times \lagForce_k\Delta s\right). 
\end{align}
To solve for the Eulerian force density $\eulForce_{i,j}$, we make the additional assumption, 
\begin{align}
    \nabla_h\cdot\eulForce_{i,j} = 0.
\end{align}
This discrete divergence-free condition placed on $\eulForce_{i,j}$ allows us to obtain the following vector Poisson equation for the Eulerian force density $\eulForce_{i,j}$ to be spread to the background grid
\begin{align}
    \label{eq:pois_f}
    -\Delta_h\eulForce_{i,j} = \nabla_h^{\perp}\sum_{k=0}^{M-1}\left(\nabla\delta_h\left(\xx_{i-\frac{1}{2},j-\frac{1}{2}}-\XX_k\right)\times \lagForce_k\Delta s\right).
\end{align}
Upon solving this equation, the unique solution for the Eulerian force density $\eulForce_{i,j}$ is determined by adding $\eulForce_0$ to the computed solution via equation \eqref{eq:eul_force_mean}. We denote the DFIB force spreading operation using the notation 
\begin{align}
    \eulForce_{i,j} = \mathcal{S}_{\text{DFIB}}\left[\XX\right]\lagForce.
\end{align}
The fact that the DFIB method spreads an Eulerian force density that is always discretely divergence-free is not overly restrictive, since only the divergence-free component of the spread force affects the flow field. However, because the spreading operation filters out the discrete non-solenoidal component of the spread force, the computed pressure differs from the physical pressure, which should account for the non-solenoidal forces arising from the presence of the immersed boundary. 

\subsubsection{Local IB Velocity Interpolation and Force Spreading}\label{sec:IB_Lag_spread_interp}
Given a discrete Eulerian velocity field $\eulVel$, we compute the horizontal $U_k$ and vertical $V_k$ components of the velocity associated with the $k^{\text{th}}$ Lagrangian marker point via
\begin{align}
    U_k &= \sum_{i,j = 0}^{N-1} u_{i,j}\,\delta_h\left(\mathbf{x}_{i-\frac{1}{2},j} - \XX(s_k,t)\right)\,h^2,\\
    V_k &= \sum_{i,j = 0}^{N-1} v_{i,j}\,\delta_h\left(\mathbf{x}_{i,j-\frac{1}{2}} - \XX(s_k,t)\right)\,h^2,
\end{align} 
in which $\mathbf{U}_k = \left(U_k,V_k\right)$. While the Eulerian velocity field satisfies a discrete incompressibility constraint ($\nabla_h\cdot\eulVel_{i,j} = 0$), the interpolated velocity field $\mathbf{U}(\XX)$ is generally not continuously divergence-free (i.e., $\nabla_{\XX}\cdot\mathbf{U}(\XX) \neq 0$) when using standard isotropic regularized delta functions. This failure to maintain a divergence-free velocity field can lead to spurious volume changes during Lagrangian marker advection. For composite B-spline regularized delta functions, the interpolation operation preserves the divergence-free property, converting discretely divergence-free velocity fields into continuously divergence-free interpolants \cite{schroeder2022,handscomb1984}. 

For completeness we provide a proof of this property. The proof assumes that either the computational domain is periodic or that all points $\XX$ where the velocity is being interpolated are located at least one kernel width away from the boundary of the staggered Cartesian grid. Under these assumptions, using composite B-spline regularized delta functions to interpolate the velocity yields a continuously divergence-free result. Because of their tensor product nature, the interpolated velocity remains continuously divergence-free even when the Cartesian grid is anisotropic or for locally adapted grids, provided that the support of the regularized delta function does not overlap a coarse-fine interface.

Under these assumptions, using composite B-spline regularized delta functions to interpolate the velocity yields
\begin{align}
U\left(\XX\right) &= \sum_{i,j}u_{i,j}\text{BS}_{n+1}\left(\frac{x_{i-\frac{1}{2}}-X}{h}\right)\text{BS}_{n}\left(\frac{y_{j}-Y}{h}\right), \\
V\left(\XX\right) &= \sum_{i,j}v_{i,j}\text{BS}_{n}\left(\frac{x_i - X}{h}\right)\text{BS}_{n+1}\left(\frac{y_{j-\frac{1}{2}}-Y}{h}\right).
\end{align}
Taking the continuous divergence of $\lagVel(\XX)$ with respect to $\XX$, we obtain
\begin{align}
    \frac{\partial U}{\partial X}(\XX) + \frac{\partial V}{\partial Y}(\XX) = -\frac{1}{h}\sum_{i,j}u_{i,j}\text{BS}^{\prime}_{n+1}\left(\frac{x_{i-\frac{1}{2}}-X}{h}\right)\text{BS}_{n}\left(\frac{y_{j}-Y}{h}\right) +v_{i,j}\text{BS}_{n}\left(\frac{x_{i}-X}{h}\right)\text{BS}^{\prime}_{n+1}\left(\frac{y_{j-\frac{1}{2}}-Y}{h}\right) .
\end{align}
Using the derivative identity \eqref{eq:deriv_prop} associated with the B-spline sequence, the continuous divergence may be expanded to 
\begin{align}
    \frac{\partial U}{\partial X}(\XX) + \frac{\partial V}{\partial Y}(\XX) &= -\frac{1}{h}\sum_{i,j}u_{i,j}\left(\text{BS}_n\parens{\frac{x_{i-\frac{1}{2}}-X}{h} + \frac{1}{2}} - \text{BS}_n\parens{\frac{x_{i-\frac{1}{2}}-X}{h} - \frac{1}{2}}\right)\text{BS}_n\parens{\frac{y_{j}-Y}{h}}\notag \\
    &- \frac{1}{h}\sum_{i,j} v_{i,j}\text{BS}_{n}\parens{\frac{x_{i} - X}{h}}\left(\text{BS}_n\parens{\frac{y_{j-\frac{1}{2}}-Y}{h} + \frac{1}{2}} - \text{BS}_n\parens{\frac{y_{j-\frac{1}{2}}-Y}{h} - \frac{1}{2}}\right), \\
    &= -\frac{1}{h}\sum_{i,j}u_{i,j}\left(\text{BS}_n\parens{\frac{x_{i}-X}{h}} - \text{BS}_n\parens{\frac{x_{i-1}-X}{h}}\right)\text{BS}_n\parens{\frac{y_{j}-Y}{h}}\notag \\
    &- \frac{1}{h}\sum_{i,j} v_{i,j}\text{BS}_{n}\parens{\frac{x_{i} - X}{h}}\left(\text{BS}_n\parens{\frac{y_{j}-Y}{h}} - \text{BS}_n\parens{\frac{y_{j-1}-Y}{h}}\right).
\end{align}
\begin{samepage}Applying summation by parts to the sums $\sum_{i,j}u_{i,j}\text{BS}_n\parens{\frac{x_{i-1}-X}{h}}$ and $\sum_{i,j}v_{i,j}\text{BS}_n\parens{\frac{y_{j-1}-Y}{h}}$ transforms the above into 
\begin{align}
    \frac{1}{h}\sum_{i,j}\left(u_{i+1,j}-u_{i,j} + v_{i,j+1}-v_{i,j}\right)\text{BS}_n\parens{\frac{x_i - X}{h}}\text{BS}_n\parens{\frac{y_j - X}{h}}, \\
    = \sum_{i,j}\nabla_h\cdot\mathbf{u}_{i,j}\text{BS}_n\parens{\frac{x_i - X}{h}}\text{BS}_n\parens{\frac{y_j - X}{h}},
\end{align}
which is zero according to our assumption that $\eulVel_{i,j}$ is discretely divergence free.\end{samepage} Thus, the interpolated velocity field satisfies $\nabla_{\XX}\cdot\mathbf{U}(\XX) = 0$ pointwise and will help mitigate volume conversation errors which result from the velocity interpolation operation. We emphasize that this proof relies on only two key mathematical ingredients: the derivative identity \eqref{eq:deriv_prop} and summation by parts. Consequently, any composite combination of kernel functions satisfying the same derivative identity will yield an interpolation operator that converts discretely divergence-free velocity fields into continuously divergence-free ones. For instance, one could construct such an interpolant using a composite combination of the \IBfour~kernel function and its convolution with the step function $\text{BS}_1$. 

Additionally, we note that although the interpolated velocity field is pointwise continuously divergence-free, the resulting interpolated velocity is still computed with first-order interpolation errors. This limitation arises because the velocity generally has a jump discontinuity in its normal derivative across the immersed boundary~\cite{gruninger2024,griffith2007,li1994}. If the velocity were smooth, using composite B-splines would yield second-order accuracy since they satisfy the first moment conditions~\cite{peskin2002}. However, since the regularity of the true velocity field is inherently limited by the presence of the immersed boundary, employing higher-order interpolation schemes does not improve performance. For this reason, we avoid using the higher-order continuously divergence-free interpolation schemes suggested by Chowdhury, Shinar, and Schroeder~\cite{chowdhury2024}.

Analogous to the description of the DFIB method, the force spreading operator is chosen to be the discrete adjoint of the velocity interpolation operator to ensure conservation of energy in Lagrangian-Eulerian interactions\cite{peskin2002}. \begin{samepage}For a discrete Lagrangian force density $\lagForce = \left(F^x,F^y\right)$, the spread force $\eulForce = \left(f^x, f^y\right)$ is computed according to
\begin{align}
    f^x_{i,j} = \sum_{k=0}^{M-1} F^x(s_k,t)\,\delta_h\left(\xx_{i-\frac{1}{2},j} - \XX(s_k,t)\right)\,\Delta s, \\
    f^y_{i,j} = \sum_{k=0}^{M-1} F^y(s_k,t)\,\delta_h\left(\xx_{i,j-\frac{1}{2}} - \XX(s_k,t)\right)\,\Delta s.
\end{align}
\end{samepage}
A property of composite B-spline regularized delta functions which to our knowledge is previously unrecognized, is their ability to preserve gradient structure when spreading Lagrangian force densities to Eulerian force densities on the grid. Specifically, when a Lagrangian force density represents a mean pressure jump across an interface, the force spreading operation transfers it to a discrete gradient of a scalar $\phi_{i,j}$ defined at the cell centers of the MAC grid:
\begin{align}
    \nabla_h\phi_{i,j} = \mathbf{f}_{i,j} = \int_{0}^{2\pi}\lagForce(s,t)\,\delta_h\left(\xx_{i,j} - \XX(s,t)\right)\,\text{d}s = \int_{\Gamma}[\![p]\!]_0\hat{\mathbf{n}}\,\delta_h\left(\xx_{i,j} - \XX\right)\,\text{d}S,
\end{align}
in which, for notational simplicity, we write $\mathbf{f}_{i,j}$ and $\delta_h(\xx_{i,j} - \XX)$ with the understanding that vector components are evaluated at their respective staggered grid locations on the MAC grid. Here $\Gamma$ is the set that describes the physical locations of the immersed boundary and $[\![p]\!]_0$ is the mean value of the jump in the pressure along the immersed boundary, defined as 
\begin{align}
    [\![p]\!]_0 = \int_{0}^{2\pi}[\![p]\!](s)\,\text{d}s.
\end{align}
In the continuous setting with singular delta functions, these forces represent distributional gradients of piecewise constant functions with discontinuities along the immersed boundary. Such forces naturally arise in fluid-structure interaction equilibria, particularly when elastic and pressure forces balance along a fluid-filled membrane at steady state. Although these forces should theoretically induce no fluid motion, standard isotropic regularized delta functions project a significant portion onto discretely divergence-free vector fields, generating spurious flows. While these spurious velocities converge to zero under grid refinement, they remain significant at practical resolutions due to their first-order pointwise convergence rate\cite{mori2008}. The combination of these spurious flows with a velocity interpolation operator that does not preserve divergence-free conditions leads to cumulative volume conservation errors. This phenomenon is particularly evident in classical IB simulations of elastic membranes, where the enclosed volume decreases linearly in time at a rate proportional to the pressure jump across the membrane interface\cite{griffith2012,long2003,peskin1993}.

In contrast to isotropic kernels, composite B-spline regularized delta functions preserve gradient structure by mapping these forces to discrete gradient fields that properly balance the pressure Lagrange multiplier. To demonstrate this property, consider a Lagrangian force density $\lagForce(s,t)$ satisfying
\begin{align}
    \label{eq:force_mean_pressure_jump}
    \int_{0}^{2\pi}\lagForce(s,t)\,\text{d}s = \int_{\Gamma}[\![p]\!]_0\hat{\mathbf{n}}\,\text{d}S,
\end{align}
in which $[\![p]\!]_0$ represents a mean pressure jump across the immersed boundary $\Gamma$ with an outward pointing unit normal vector $\hat{\mathbf{n}}$ and $\text{d}S$ is the infinitesimal arc-length increment. We note that such a force takes the form
\begin{align}
    \lagForce(s,t) = [\![p]\!]_0\hat{\mathbf{n}}(s,t)\left|\frac{\partial\XX}{\partial s}(s,t)\right|.
\end{align}
To show that composite B-spline regularized delta functions map these forces to discrete gradient fields, we can either demonstrate that $\left(I - \nabla_h\Delta_h^{-1}\nabla_h\cdot\right)\eulForce = 0$, or equivalently, using the discrete Helmholtz decomposition available on the staggered, MAC grid, show that $\nabla_h\times\eulForce = 0$. We choose the latter approach for its computational simplicity. Taking the discrete curl of the corresponding Eulerian force density $\mathbf{f}$, we get 
\begin{equation}
    \begin{multlined}
        \nabla_h\times\int_{\Gamma}[\![p]\!]_0\hat{\mathbf{n}}\,\delta_h\left(\xx_{i,j} - \XX\right)\,\text{d}S = \frac{1}{h^3}\int_{\Gamma}[\![p]\!]_0\hat{\mathbf{n}}\cdot\begin{bmatrix}
            -\text{BS}_{n+1}\left(\frac{x_{i-1/2} - X}{h}\right)\left(\text{BS}_n\left(\frac{y_j - Y}{h}\right) - \text{BS}_n\left(\frac{y_{j-1} - Y}{h}\right)\right) \\
            \left(\text{BS}_n\left(\frac{x_i - X}{h}\right) - \text{BS}_n\left(\frac{x_{i-1} - X}{h}\right)\right)\text{BS}_{n+1}\left(\frac{y_{j-1/2} - Y}{h}\right)
        \end{bmatrix}
        \,\text{d}S \\ 
        = \frac{1}{h^3}\int_{\Gamma}[\![p]\!]_0\hat{\mathbf{n}}\cdot\begin{bmatrix}
            -\text{BS}_{n+1}\left(\frac{x_{i-1/2} - X}{h}\right)\left(\text{BS}_n\left(\frac{y_{j-1/2} - Y}{h} + \frac{1}{2}\right) - \text{BS}_n\left(\frac{y_{j-1/2} - Y}{h} - \frac{1}{2}\right)\right) \\
            \left(\text{BS}_n\left(\frac{x_{i-1/2} - X}{h} + \frac{1}{2}\right) - \text{BS}_n\left(\frac{x_{i-1/2} - X}{h} - \frac{1}{2}\right)\right)\text{BS}_{n+1}\left(\frac{y_{j-1/2} - Y}{h}\right)
        \end{bmatrix}
        \,\text{d}S
    \end{multlined}
\end{equation}
Applying the derivative property associated with the B-spline family \eqref{eq:deriv_prop} and the chain rule to the above right hand side of the equation, we get 
\begin{equation}
    -\frac{1}{h^2}\int_{\Gamma}[\![p]\!]_0\hat{\mathbf{n}}\cdot\begin{bmatrix}
        -\text{BS}_{n+1}\left(\frac{x_{i-1/2} - X}{h}\right)\frac{\text{d}}{\text{d}Y}\text{BS}_{n+1}\left(\frac{y_{j-1/2} - Y}{h}\right)\\
        \frac{\text{d}}{\text{d}X}\text{BS}_{n+1}\left(\frac{x_{i-\frac{1}{2}} - X}{h}\right)\text{BS}_{n+1}\left(\frac{y_{j-1/2} - Y}{h}\right)\,\text{d}{S}
    \end{bmatrix}
\end{equation}
Applying the divergence theorem, we obtain 
\begin{equation}
    \begin{aligned}
    &-\frac{1}{h^2}\int_{\Gamma}[\![p]\!]_0\hat{\mathbf{n}}\cdot\begin{bmatrix}
        -\text{BS}_{n+1}\left(\frac{x_{i-1/2} - X}{h}\right)\frac{\text{d}}{\text{d}Y}\text{BS}_{n+1}\left(\frac{y_{j-1/2} - Y}{h}\right)\\
        \frac{\text{d}}{\text{d}X}\text{BS}_{n+1}\left(\frac{x_{i-\frac{1}{2}} - X}{h}\right)\text{BS}_{n+1}\left(\frac{y_{j-1/2} - Y}{h}\right)
    \end{bmatrix}\,\text{d}{S} \\
        &= -\frac{1}{h^3}\iint[\![p]\!]_0\left(\frac{\text{d}}{\text{d}X}\text{BS}_{n+1}\left(\frac{x_{i-1/2} - X}{h}\right)\frac{\text{d}}{\text{d}Y}\text{BS}_{n+1}\left(\frac{y_{j-1/2} - Y}{h}\right)- \frac{\text{d}}{\text{d}X}\text{BS}_{n+1}\left(\frac{x_{i-1/2} - X}{h}\right)\frac{\text{d}}{\text{d}Y}\text{BS}_{n+1}\left(\frac{y_{j-1/2} - Y}{h}\right)\right)\,\text{d}X\text{d}Y \\
        & \hspace{17em}= 0,
    \end{aligned}
\end{equation}
where the double integral is taken over the region enclosed by $\Gamma$. Due to the tensor product nature of composite B-spline regularized delta functions, this property readily extends to anisotropic Cartesian grids. The preservation of the gradient nature of the spread force also holds in the presence of physical boundary conditions, provided that the support of the composite B-spline regularized delta function does not extend past the computational boundary. Similarly, this result holds on locally adaptive Cartesian grids, provided that the delta function support does not overlap any coarse-fine interface.

This analysis demonstrates that composite B-splines preserve the gradient structure of forces corresponding to mean pressure jumps when the force spreading operator is computed exactly. In practice, both velocity interpolation and force spreading operators are discretized using the periodic trapezoidal rule, whose convergence rate depends on the Fourier series decay of the integrand \cite{boyd2001}. Consequently, smoother kernel functions provide more accurate force spreading operations and better preservation of gradient structure for forces satisfying \eqref{eq:force_mean_pressure_jump}. Further analysis of the spurious flows generated by isotropic regularized delta functions, including scaling estimates of induced spurious vorticity and velocity magnitude along with numerical validation, is provided in section \ref{sec:spur_flow_iso} of the appendix.

\subsection{Temporal Discretization}
The equations of motion \eqref{eq:mom}--\eqref{eq:spread_cont} are discretized in time using a semi-implicit scheme, previously presented for various IB method applications\cite{bao2017}. We describe the timestepping scheme using the notation for the standard IB method, noting that the same scheme applies to the DFIB method.
To advance from timestep $t^n$ to $t^{n+1} = (n+1)\Delta t$, we first approximate the Lagrangian configuration at the midpoint $t^{n+\frac{1}{2}}$
\begin{align}
\XX^{n+\frac{1}{2}} = \XX^{n} + \frac{\Delta t}{2}\mathcal{J}[\XX^n]\eulVel^n.
\end{align}
We then use $\XX^{n+\frac{1}{2}}$ to approximate the Lagrangian force density $\lagForce^{n+\frac{1}{2}}$ at $t^{n+\frac{1}{2}}$. The updated Lagrangian force density is then passed to the force spreading operator to obtain an approximation for $\eulForce^{n+\frac{1}{2}}$ at the midpoint
\begin{align}
\eulForce^{n+\frac{1}{2}} = \mathcal{S}[\XX^{n+\frac{1}{2}}]\lagForce^{n+\frac{1}{2}}.
\end{align}
Next, we solve the following discrete approximation to the momentum equation and incompressibility constraint for $\eulVel^{n+1}$ and $p^{n+\frac{1}{2}}$
\begin{samepage}
\begin{align}
\label{eq:discrete_ns}
\rho\left(\frac{\eulVel^{n+1}-\eulVel^n}{\Delta t} + \mathbf{N}(\eulVel)^{n+\frac{1}{2}}\right) &= -\nabla_h p^{n+\frac{1}{2}} + \frac{\mu}{2}\Delta_h\left(\eulVel^{n+1} + \eulVel^n\right) + \eulForce^{n+\frac{1}{2}},\\
\nabla_h\cdot\eulVel^{n+1} &= 0,
\end{align}
\end{samepage}
where $\mathbf{N}(\eulVel)^{n+\frac{1}{2}}$ is a second-order Adams-Bashforth (AB2) approximation to the advection term at the timestep midpoint
\begin{align}
\mathbf{N}(\eulVel)^{n+\frac{1}{2}} = \frac{3}{2}\mathbf{N}(\eulVel)^n - \frac{1}{2}\mathbf{N}(\eulVel)^{n-1}.
\end{align}
Finally, we compute an approximation to the Lagrangian configuration at $t^{n+1}$ using the midpoint approximation
\begin{align}
\XX^{n+1} = \XX^{n} + \frac{\Delta t}{2}\mathcal{J}[\XX^{n+\frac{1}{2}}]\left(\eulVel^{n+1} + \eulVel^{n}\right).
\end{align}
As the nonlinear advective term uses the AB2 scheme, we employ an explicit second-order Runge-Kutta (RK2) scheme for the first timestep, as previously described\cite{peskin2002}.
While both the DFIB and standard IB methods use the same time-stepping scheme, the DFIB method requires three additional scalar Poisson solves per timestep. We note that in three spatial dimensions, the DFIB method becomes substantially more computationally expensive. It requires nine additional scalar Poisson solves and nine additional scalar interpolation/spreading steps. This makes it approximately twice as costly per timestep in two dimensions and three times more expensive in three dimensions compared to the ordinary IB method \cite{bao2017}.

To solve equation \eqref{eq:discrete_ns} at each timestep, we utilize the projection method. We first solve for the pressure by taking the discrete divergence of equation \eqref{eq:discrete_ns} and solving the resulting scalar Poisson equation for the pressure using the Fast Fourier Transform (FFT) available on the periodic grid. After computing the pressure $p^{n+\frac{1}{2}}$, we compute its discrete gradient and solve the resulting vector equations for the velocity variables using the FFT. For non-periodic grids, we solve the fully coupled saddle-point system posed by \eqref{eq:discrete_ns} and use the projection method as a preconditioner as described by Griffith~\cite{griffith2009accurate}.

\subsection{Area Conservation Measurements}\label{sec:area_measurements}
A consequence of the incompressibility of the fluid is that the initial area (or volume) of a closed curve (or surface) remains constant if advected by the fluid. We use this concept to evaluate how well an immersed boundary simulation, combined with a specific regularized delta function, maintains the continuous incompressibility constraint of the fluid. To compute the areas of closed immersed boundaries, we introduce tracer points $\XXt$, positioned according to the initial configuration of the immersed boundary but advected solely by interpolated fluid velocity. This means the tracer points do not generate any Lagrangian forces that are spread back onto the fluid grid. We measure and track the area enclosed by the tracer points by invoking Green's theorem in the form 
\begin{equation}
    A(t;\XXt) = \frac{1}{2}\iint_{D_{\text{tracer}}}\nabla\cdot\begin{bmatrix} x \\ y\end{bmatrix}\,\text{d}x\,\text{d}y = \frac{1}{2}\int_{0}^{2\pi}\left(X_{\text{tracer}}(s,t)\frac{\partial Y_{\text{tracer}}}{\partial s}(s,t) - Y_{\text{tracer}}(s,t)\frac{\partial X_{\text{tracer}}}{\partial s}(s,t)\right)\,\text{d}s
\end{equation}
in which $D_{\text{tracer}}$ is the domain enclosed by the closed curve discretized by the tracers $\XXt$. To compute this quantity, we use cubic splines to interpolate the coordinates of the tracers. We then differentiate the spline interpolant of the tracers and compute the integral given by Green's theorem by replacing $X_\text{tracer}$, $\frac{\partial X_{\text{tracer}}}{\partial s}$, $Y_{\text{tracer}}$, and $\frac{\partial Y_{\text{tracer}}}{\partial s}$ with their spline representations and integrate these splines exactly. The number of tracer points $N_{\text{tracer}}$ we use is an integer multiple of the number of Lagrangian marker points $N_{\text{IB}}$ employed to discretize the physical immersed boundary. The integer multiple is chosen large enough so that the initial area enclosed by the immersed boundary is computed to about machine precision. In each of our tests, we report the relative change in the area at each timestep
\begin{equation}
    \Delta A(t;\XXt) = \frac{\left|A(t;\XXt) - A_{\text{initial}}\right|}{A_{\text{initial}}},
\end{equation} 
in which $A_{\text{initial}}$ is the exact initial area enclosed by immersed boundary.

\section{Numerical Tests and Discussion}\label{sec:numerical tests}
This section presents empirical tests demonstrating the effectiveness of composite B-splines in preserving areas of immersed boundaries defined by closed curves. To benchmark the composite B-spline regularized delta functions, we compare results with the IB method using the isotropic \IBfour~regularized delta function and the (non-local) DFIB method using the isotropic \IBsix~regularized delta function. All physical quantities in this paper are reported using the centimeter-gram-second (CGS) unit system. The fluid is characterized by a constant density of $\rho = 1.0$ and a dynamic viscosity value of $\mu = 0.1$, unless explicitly stated otherwise for a particular simulation.\par

\subsection{Advecting an Immersed Boundary Using Taylor Vortices}\label{sec:advect_test_sec}
We first test the ability of the interpolation scheme to reconstruct incompressible material trajectories for a case that does not involve fluid-structure interaction by comparing the standard IB method and the DFIB method when applied to advecting Lagrangian tracers $\XX_{\text{tracer}}$ according to a Taylor-Green vortex flow on the periodic unit square. The associated velocity and pressure fields are given by:
\begin{align}
    u(x,y,t) &= 1 + 2e^{-8\pi^2\frac{\mu}{\rho}t} \sin(2\pi(y-t)) \cos(2\pi(x-t)), \\
    v(x,y,t) &= 1 - 2e^{-8\pi^2\frac{\mu}{\rho}t} \cos(2\pi(y-t)) \sin(2\pi(x-t)), \\
    p(x,y,t) &= -e^{-16\pi^2\frac{\mu}{\rho}t} (\cos(4\pi(x-t)) + \cos(4\pi(y-t))).
\end{align}
The tracers are initialized so that they discretize the circle centered at $\left(\frac{1}{2},\frac{1}{2}\right)$ and radius of $r = \frac{1}{4}$. We emphasize that the tracers are purely advected using the interpolation strategies of both the IB and DFIB methods. No Lagrangian force densities are spread onto the Cartesian grid.\par  
To evaluate area conservation properties, we simulate advection of the tracers on a uniform grid with increment $h = \frac{1}{32}$. We use eight time step sizes: $\Delta t = \frac{h}{8}, \frac{h}{16}, \frac{h}{32}, \frac{h}{64}, \frac{h}{128}, \frac{h}{256},$ $\frac{h}{512}$, and $\frac{h}{1024}$. For each method and choice of regularized delta function, we compute the temporal mean of the relative area error $\overline{\Delta A}\parens{\XX_{\text{tracer}}}$ over the interval $t = [0,1]$. In this context, area computation errors may arise from two different sources: (1) the interpolated velocity field not being continuously divergence-free, or (2) the time-stepping scheme introducing errors in tracer positions. Both the DFIB and IB methods update the positions of the tracers using the explicit midpoint rule, which has a global truncation error that is second-order accurate. Consequently, we expect the errors introduced by the time-stepping scheme to be of order $\mathcal{O}(\Delta t^2)$ until the error associated with the incompressibility of the interpolated velocity field becomes dominant. Additionally, there are errors associated with the cubic spline interpolation of the tracer points. However, we've observed that these interpolation errors are much smaller in magnitude than the errors committed by our choice of time stepping scheme. \par
\begin{figure}[ht]
    \centering
    \includegraphics[width=0.7\textwidth]{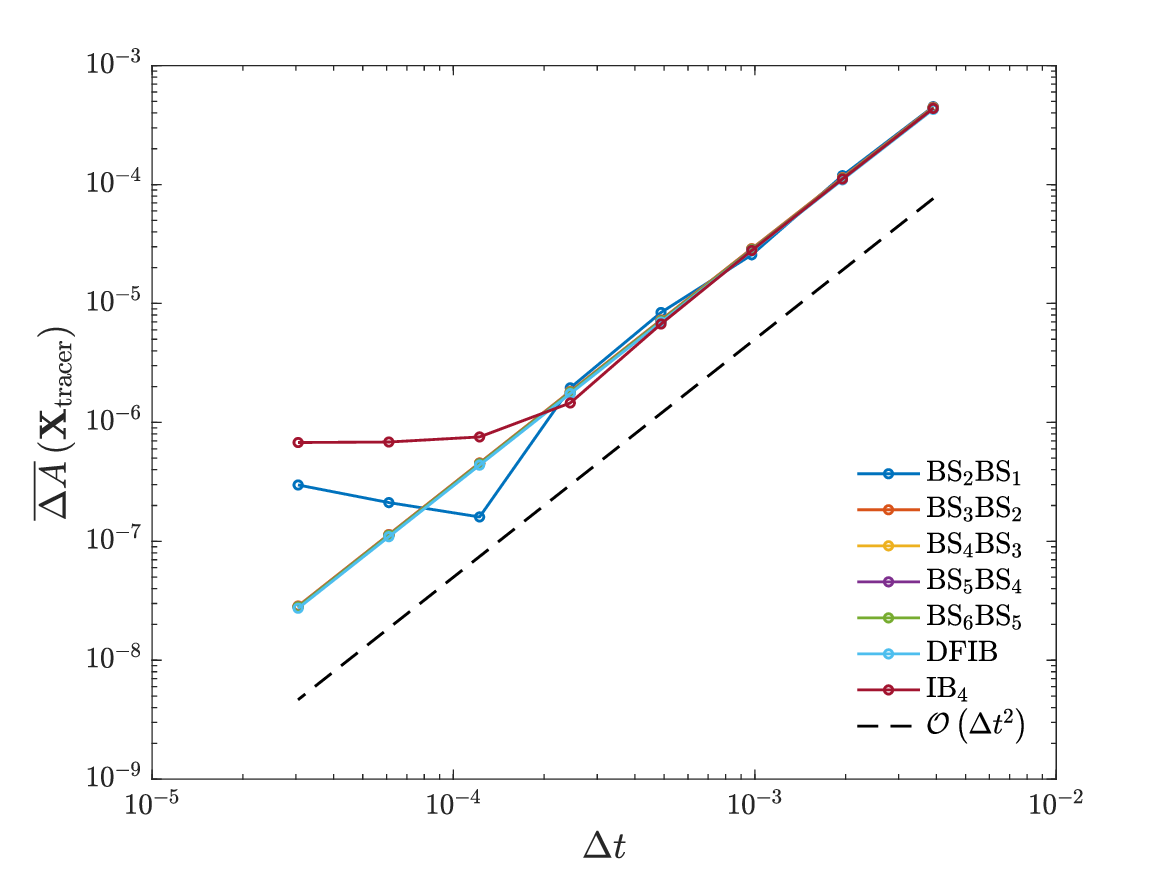}
    \caption{Mean relative area errors $\overline{\Delta A}\left(\XX_{\text{tracer}}\right)$ for regularized delta functions in the IB method and for the DFIB method.}
    \label{fig:Advect_test}
\end{figure}
The results are illustrated in Figure \ref{fig:Advect_test}. The relative area errors of the DFIB and IB methods implemented with each of the composite B-splines, except the discontinuous \BStwo\BSone~kernel, are completely superimposed and exhibit second-order convergence with respect to the time step size $\Delta t$. For larger time step sizes, the IB method with the \BStwo\BSone~and~\IBfour regularized delta functions also exhibits second-order convergence. However, for the \IBfour~kernel, upon decreasing the time step size below $\Delta t = \frac{h}{128}$, the relative area error begins to level off as the error associated with the incompressibility of the interpolated velocity dominates. Similarly, the relative error levels off for the \BStwo\BSone~kernel below a time step size of $\Delta t = \frac{h}{256}$.\par 
Although the \BStwo\BSone~kernel provides a continuously divergence-free interpolant, we observe that the error levels off as $\Delta t$ decreases. This unexpected behavior can be attributed to the discontinuous nature of the resulting interpolated velocity field. These discontinuities in the interpolated velocity field lead to discontinuities in the resulting time derivatives, which in turn reduce the asymptotic local truncation error of the explicit midpoint rule to only first-order in time. Consequently, the asymptotic global truncation error remains $\mathcal{O}(1)$ as $\Delta t$ approaches zero. Supplemental analysis regarding this phenomenon is provided in section \ref{sec:error_analysis_ts} in the appendix.

\subsection{A Pressurized Circular Membrane at Equilibrium}\label{sec:pres_membrane}
Next we consider a FSI problem involving a quasi-static pressurized membrane initialized in its circular equilibrium configuration with center $\left(\frac{1}{2},\frac{1}{2}\right)$ and radius $r = \frac{1}{4}$ in the periodic unit square $\Omega = [0,1]^2$, with zero initial background flow. The Lagrangian force density is described by
\begin{equation}
\label{eq:lag_force_press_membrane}
\lagForce(s,t) = \kappa\frac{\partial^2\XX(s,t)}{\partial s^2},
\end{equation}
and is discretized using centered, second-order accurate finite differences
\begin{align}
\label{eq:discrete_lag_force}
\lagForce_k = \frac{\kappa}{\Delta s^2}\left(\XX_{k+1} + \XX_{k-1} - 2\XX_k\right).
\end{align}
This discrete formulation models the Lagrangian markers as being linked by linear springs with zero rest lengths and uniform stiffness $\frac{\kappa}{\Delta s}$. At equilibrium, this force density should generate no fluid motion and produce a pressure field with a jump of magnitude $\kappa$ across the interface. This test case therefore allows us to verify our earlier analysis of how composite B-spline regularized delta functions preserve gradient structure when spreading such forces to the grid.\par
\begin{figure}[htbp]
    \centering 
    \includegraphics*[width=0.75\textwidth]{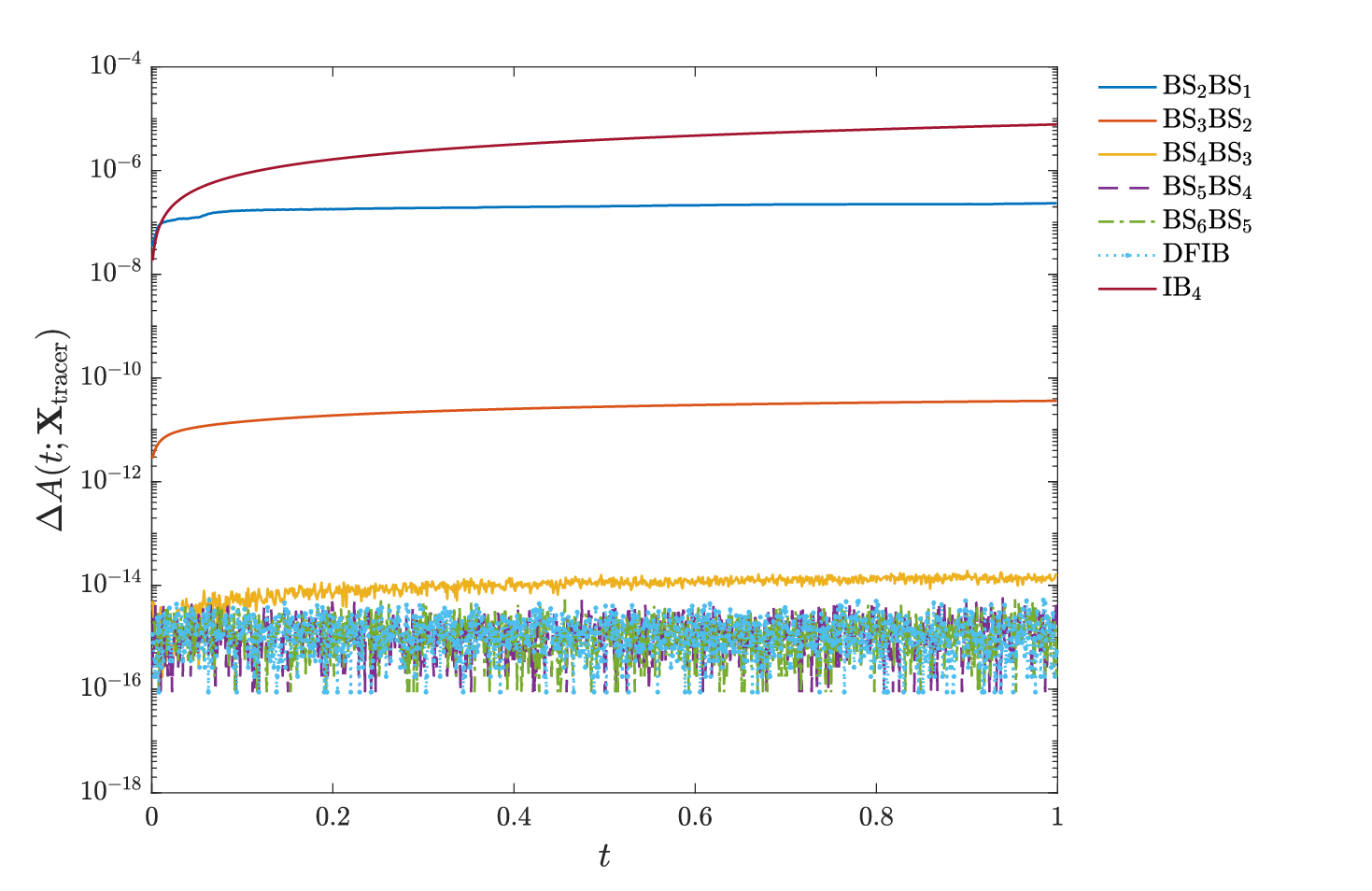}
    \caption{Semi-log plot of time-dependent relative area errors.}
    \label{fig:Membrane_Area_tracer}
\end{figure}

 
We set $\kappa = 1$, $\Delta t = h/8$ with $h = 1/128$, and unless otherwise noted, $M_{\text{fac}} = 1/2$. Because the immersed membrane is initialized at equilibrium, we consider any deviation from the initial area enclosed by the membrane to be attributed to errors associated with spatial discretization of the IB method employed. Griffith \cite{griffith2012} simulated a pressurized membrane at equilibrium using the IB method in which the Eulerian variables were discretized according to the MAC scheme \cite{harlow1965} and the velocity interpolation and force spreading operators were implemented using Peskin's four-point delta function \cite{peskin2002}. Griffith found that although the MAC discretization of the IB method exhibited better volume conservation properties compared to a collocated discretization, the IB method discretized on the MAC grid still exhibited persistent volume loss. Bao et al.\cite{bao2017} found that so long as the the immersed boundary was resolved enough, the DFIB method was able to maintain the initial volume of the immersed boundary to within machine precision. In Fig.~\ref{fig:Membrane_Area_tracer} we compare the DFIB method implemented with the $C^3$ six-point IB kernel to the standard IB method implemented with the four-point IB and composite B-spline regularized delta functions. 

Consistent with the findings of Bao et al., we observe that the DFIB method preserves the initial area of the circle to within machine precision throughout the simulation. Likewise, the two composite B-spline pairs with the highest regularity, $\text{BS}_5\text{BS}_4$ and $\text{BS}_6\text{BS}_5$, also maintain the initial area of the circle to within machine precision. Although composite B-spline pairs of lower regularity do not achieve the same level of accuracy as their higher regularity counterparts, they still significantly outperform the standard four-point IB kernel, with the least regular composite B-spline pairing, $\text{BS}_2\text{BS}_1$, performing about an order of magnitude better than the standard four-point kernel. Furthermore, we note that the relative area errors for the \BStwo\BSone~and \IBfour~regularized delta functions level off about where their mean errors level off in the pure advection test illustrated in Fig \ref{fig:Advect_test}. \par
\begin{figure}[htbp]
    \centering
    \includegraphics[width=0.475\textwidth]{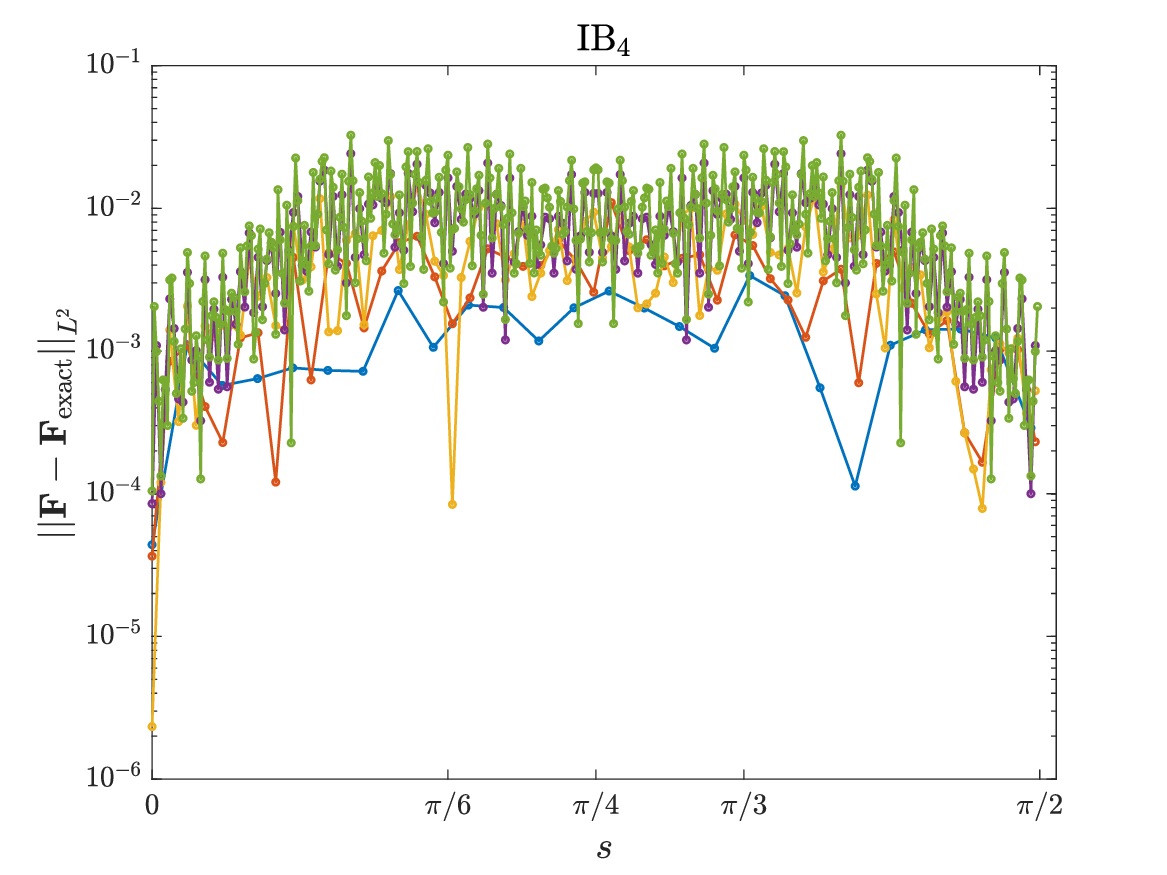}
    \includegraphics[width=0.475\textwidth]{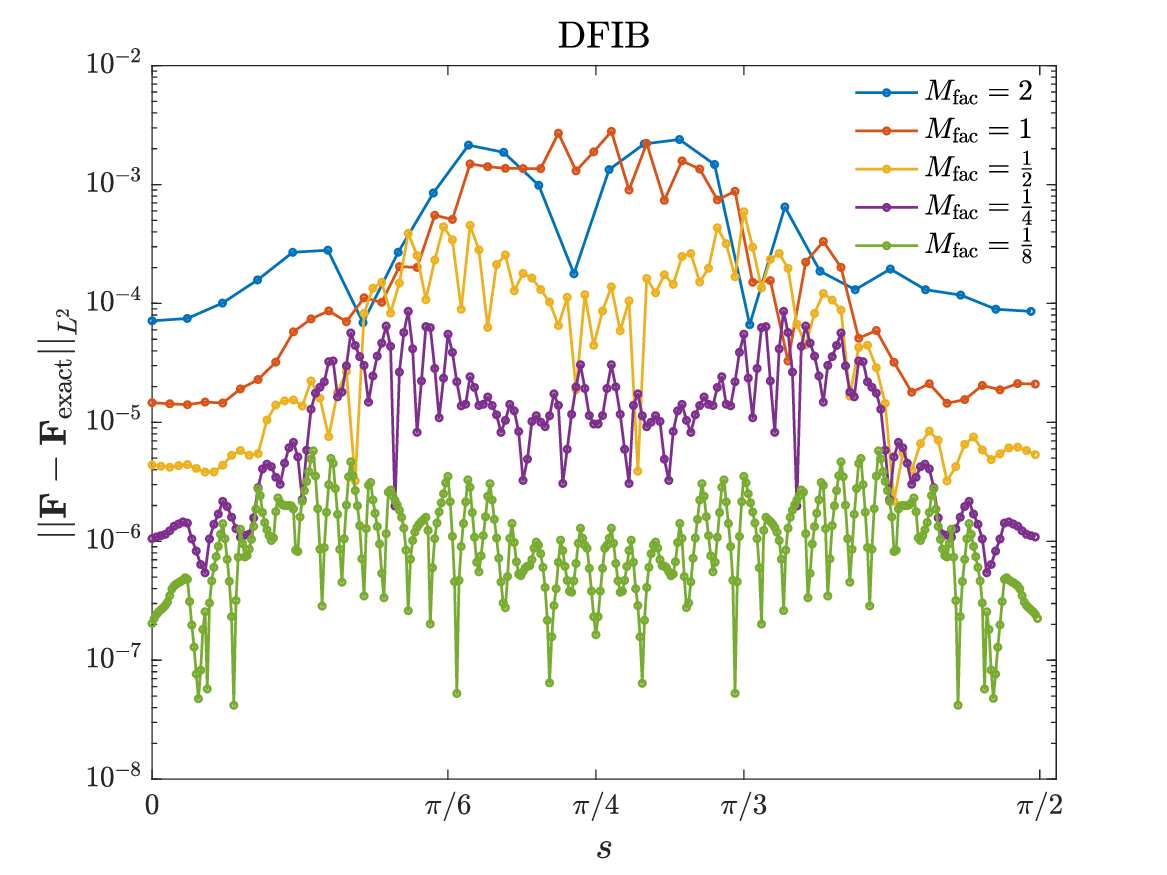}
    \includegraphics[width=0.475\textwidth]{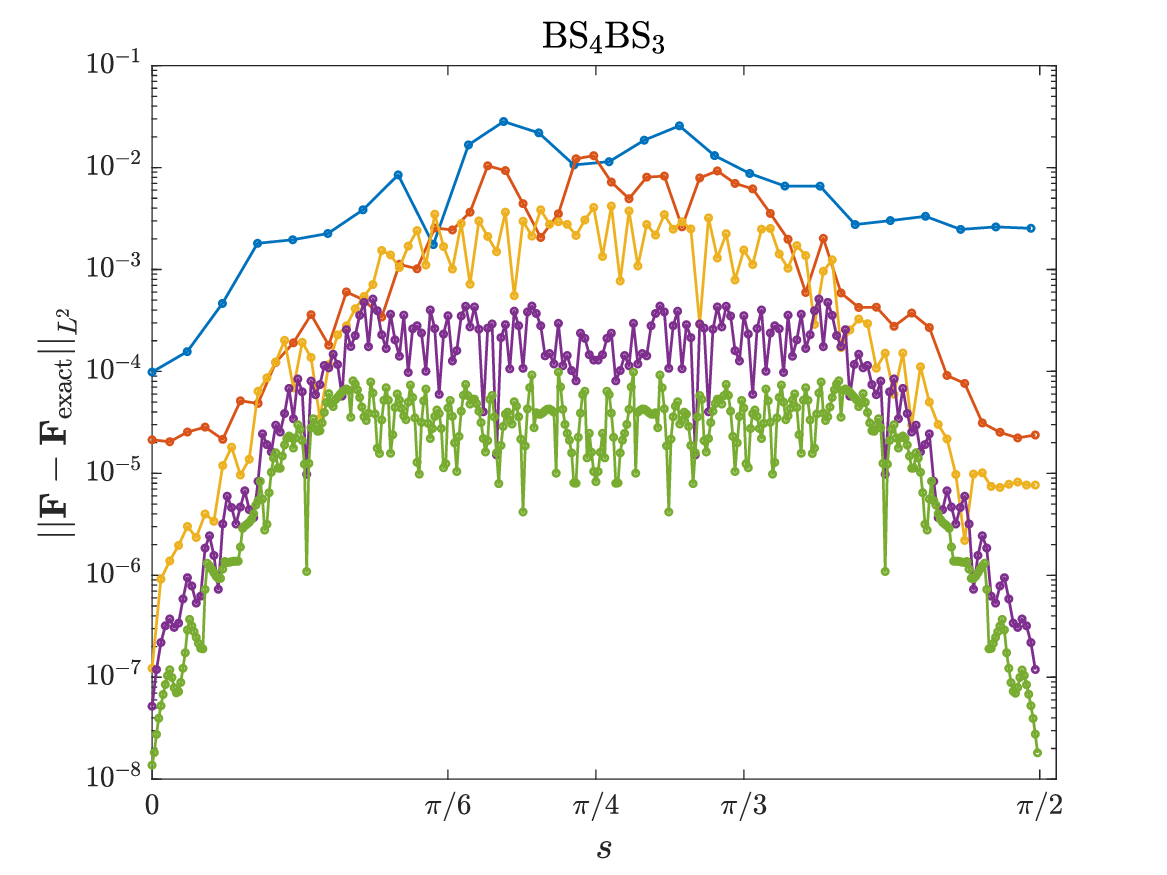}
    \includegraphics[width=0.475\textwidth]{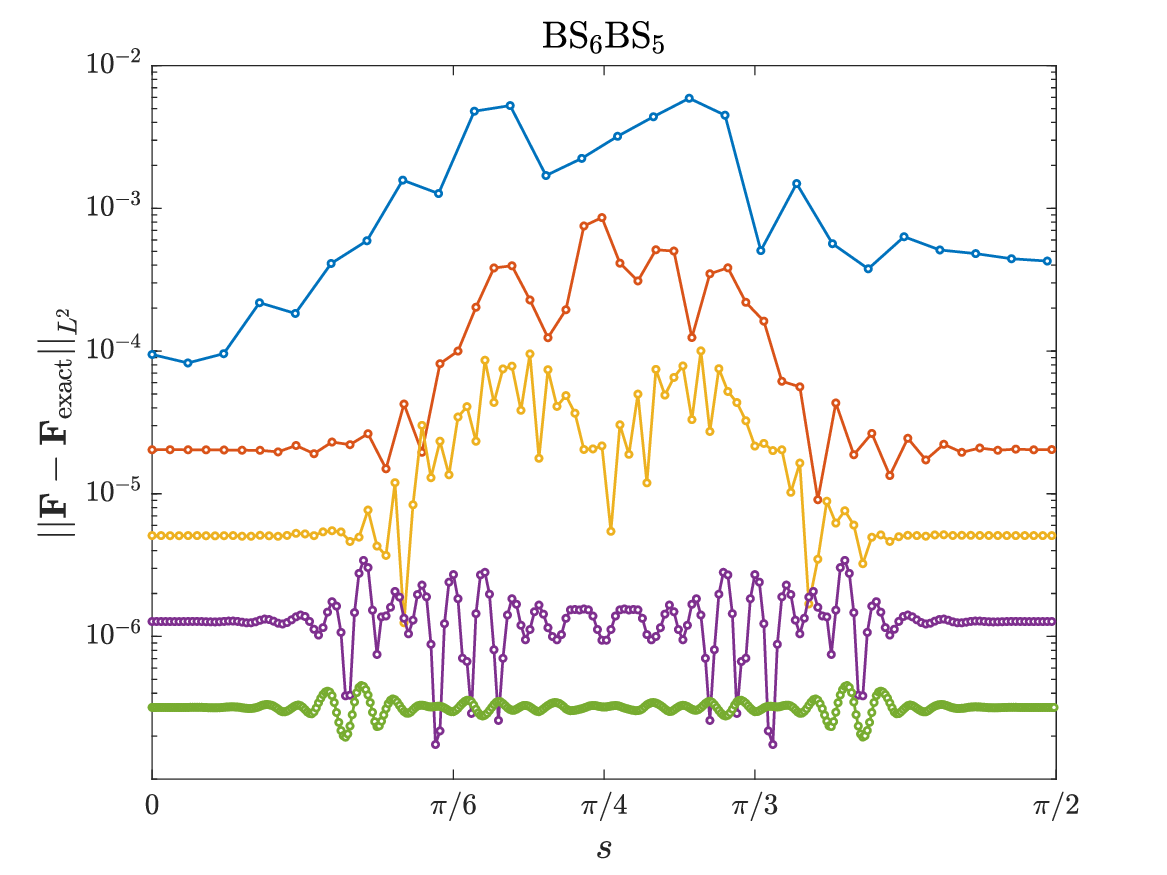}
    \caption{Semi-log plot of the error $||\lagForce-\lagForce_{\text{exact}}||_{L^2}$ for the \BSfour\BSthree, \BSsix\BSfive, \IBfour~regularized delta functions and the DFIB method for a uniform Cartesian grid discretization with $h = 1/128$ and Lagrangian mesh discretizations corresponding to $M_{\text{fac}} = 2, 1, \frac{1}{2}, \frac{1}{4}, \frac{1}{8}$.}
    \label{fig:force_accuraies}
\end{figure}

\begin{figure}
    \centering 
    \includegraphics[width=0.6\textwidth]{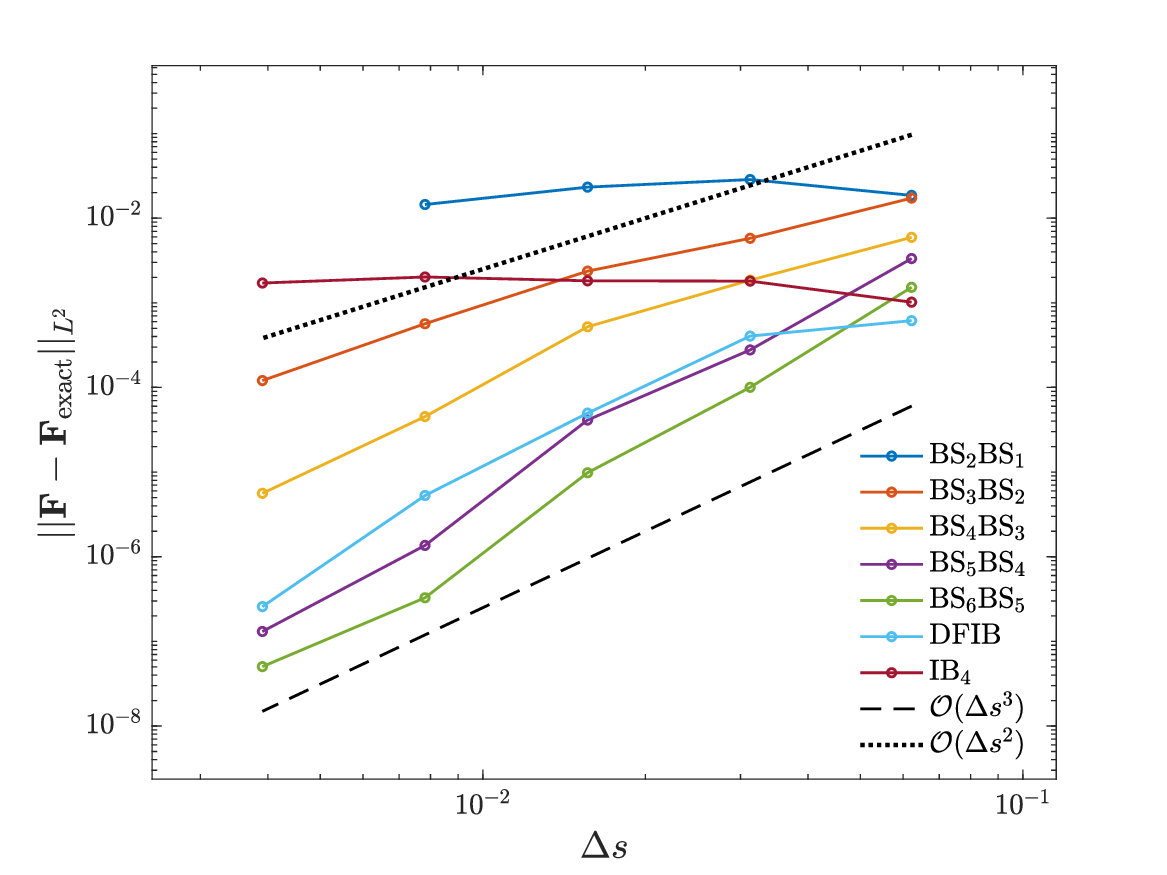}
    \caption{Log-log plot of the $L^2$ grid norm errors of the Lagrangian force densities corresponding for a uniform Cartesian grid discretization with $h = 1/128$ and Lagrangian mesh discretizations corresponding to $M_{\text{fac}} = 2, 1, \frac{1}{2}, \frac{1}{4}, \frac{1}{8}$. The results for the \BStwo\BSone~kernel at $M_{\text{fac}} = \frac{1}{8}$ were omitted because the simulation became unstable.}
    \label{fig:force_convergence}
\end{figure}

 
To analyze force computation accuracy, we compare numerical results with the exact Lagrangian force density $\lagForce_{\text{exact}}(s) = -\kappa r\hat{\mathbf{n}}(s)$ for the equilibrium circular configuration, with $\hat{\mathbf{n}}(s)$ denoting the outward unit normal. Figs \ref{fig:force_accuraies} and \ref{fig:force_convergence} show the pointwise and $L^2$ errors in the computed steady-state forces for various kernels, respectively. The $L^2$ grid norm errors associated with the Lagrangian force density is given by the formula
\begin{align}
    \left|\left|\lagForce-\lagForce_{\text{exact}}  \right|\right|_{L^2} = \Delta s\left(\sum_{k=0}^{M-1}\left(\lagForce_k-\lagForce_{k\,\text{exact}}\right)^2\right)^{1/2}.
\end{align}


\begin{figure}[htbp]
    \centering 
    \includegraphics[width=0.575\textwidth]{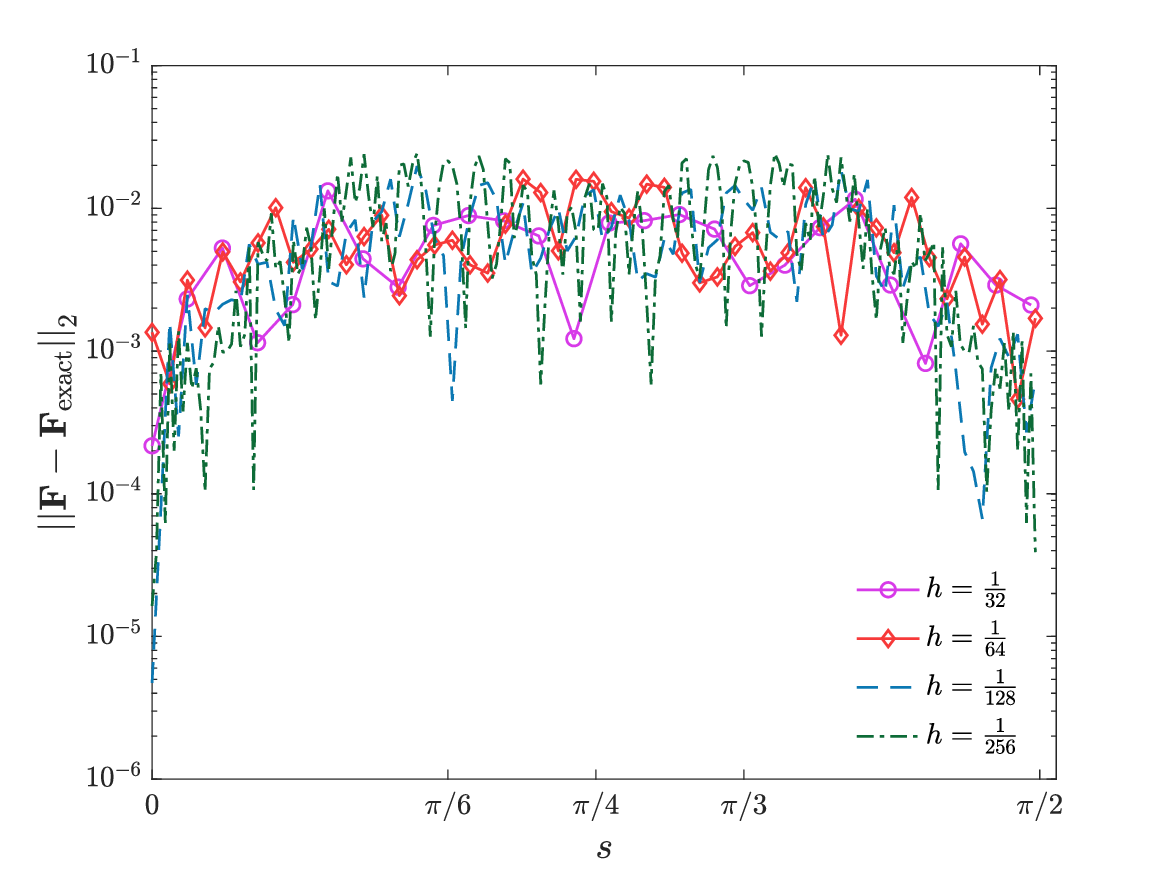}
    \caption{Euclidean error in the computed steady-state Lagrangian force densities using the \IBfour~regularized delta function. For each simulation, the time step size is set to $\Delta t = \frac{h}{8}$, and the Lagrangian grid is chosen such that $M_{\text{fac}} = \frac{1}{2}$. This ensures that both the time step and Lagrangian grid are simultaneously refined with the Cartesian grid increment $h$.}
    \label{fig:h_refine_IB4}
\end{figure}

We observe varying convergence behaviors in Lagrangian force densities across different kernels and methods. The composite B-splines, excluding only the discontinuous \BStwo\BSone~kernel, empirically produce Lagrangian force densities that converge pointwise under solely Lagrangian grid refinement. This pointwise convergence is demonstrated for the \BSsix\BSfive~and \BSfour\BSthree~kernels, as well as the DFIB method in Fig \ref{fig:force_accuraies}. In contrast, the \IBfour~kernel does not yield pointwise convergent Lagrangian force densities under these conditions. In fact, for the \IBfour~kernel, this lack of convergence persists even under simultaneous refinements of both grid sizes and time step size, as shown in Fig \ref{fig:h_refine_IB4}.\par
Examining the $L^2$ grid norm convergence of the Lagrangian force densities, as shown in Fig \ref{fig:force_convergence}, we find that the \BSsix\BSfive, \BSfive\BSfour, \BSfour\BSthree, and \BSthree\BStwo~regularized delta functions, along with the DFIB method, all provide $L^2$ convergence under Lagrangian grid refinement, with the more regular kernels generally producing smaller errors. The less regular \BSthree\BStwo~kernel appears to yield $L^2$ convergence rates proportional to $\Delta s^2$, while the more regular composite B-spline kernels and DFIB method produce asymptotic rates proportional to $\Delta s^3$. Notably, both the \IBfour~regularized delta function and the \BStwo\BSone~kernel fail to produce $L^2$ convergent Lagrangian force densities under solely Lagrangian grid refinement. The inaccuracies associated with the \BStwo\BSone~kernel can be attributed to the discontinuous nature of the interpolated velocity $\lagVel(s,t)$ in both time $t$ and the curvilinear coordinate $s$. These discontinuities in the interpolated velocity field propagate to the Lagrangian marker positions, resulting in jump discontinuities in their trajectories with respect to both $s$ and $t$. Consequently, large localized errors arise in the computed Lagrangian force densities near these discontinuity points. A more detailed discussion of the discontinuities in the interpolated velocity obtained using the \BStwo\BSone~kernel is presented in section \ref{sec:error_analysis_ts} of the appendix.\par

For the pressurized membrane problem, we note that as long as the time step is simultaneously refined with the Cartesian grid increment $h>0$, the errors in the Lagrangian force densities using the IB method with composite B-spline regularized delta functions (excluding \BStwo\BSone) and the DFIB method appear to depend only on the value of $\Delta s$.

\begin{figure}[htbp]
    \centering 
    \includegraphics[width=0.475\textwidth]{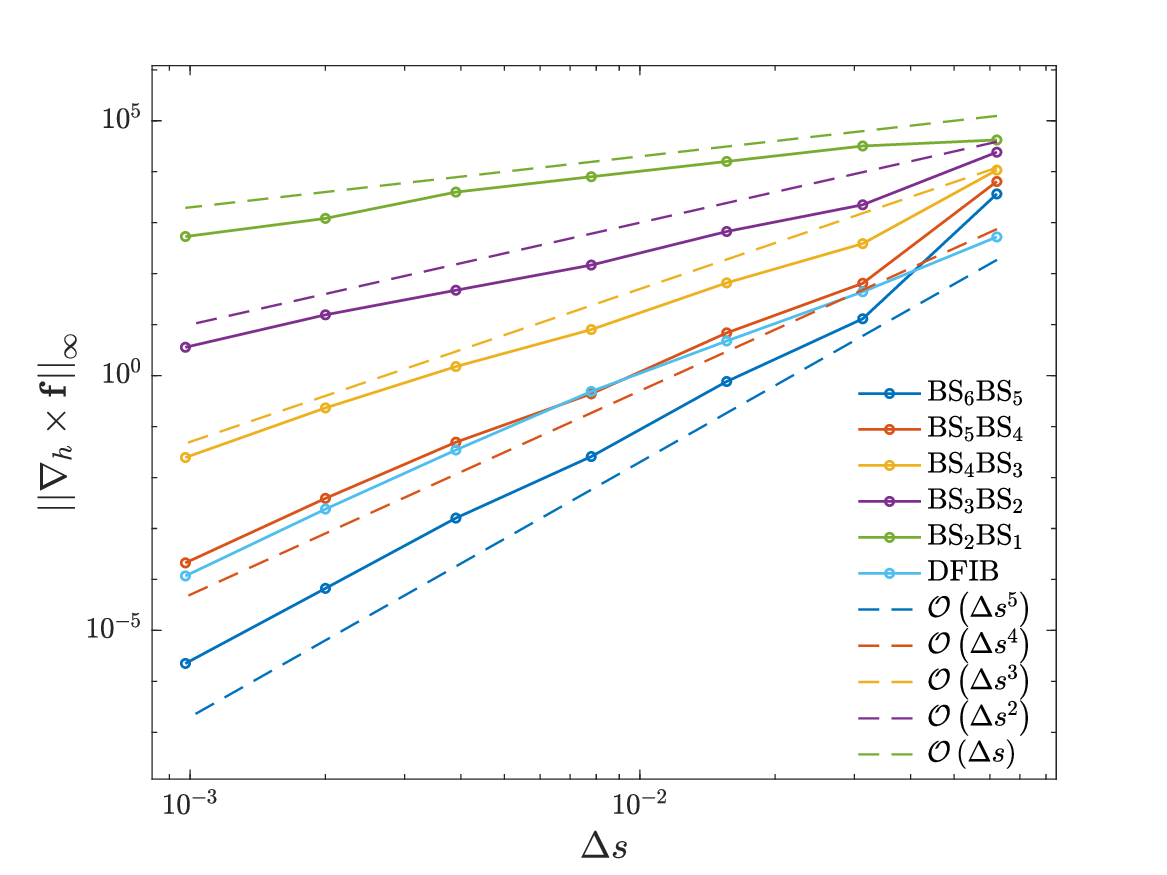}
    \includegraphics[width=0.45\textwidth]{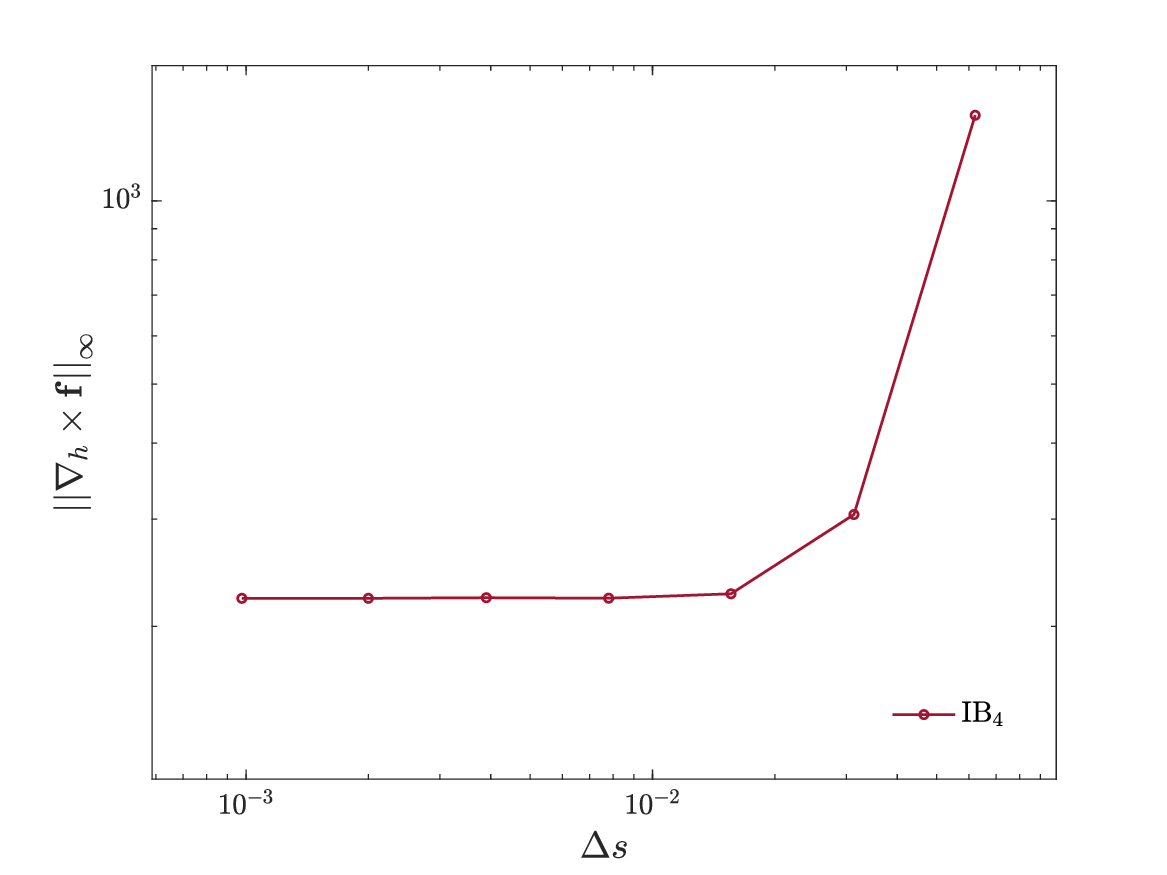}
    \caption{Plots of the maximum curl of the discrete Eulerian force density are presented. The figure on the right illustrates the convergence rates as the Lagrangian parameter grid size $\Delta s$ approaches zero. Kernels with higher regularity yield faster convergence rates, in accordance with the error estimates associated with the periodic trapezoidal rule. The loglog plot on the right demonstrates that the discrete curl of the IB4 force spreading operator does not vanish as $\Delta s$ is refined.}
    \label{fig:curl_errors}
\end{figure}

Having analyzed volume conservation and force accuracy, we next examine how these regularized delta functions spread forces to the grid. For the equilibrium configuration, the Lagrangian force density represents a zero-mean pressure jump across the interface which, in the continuous setting, corresponds to a distributional gradient and therefore induces no flow. To assess whether this gradient structure is maintained in the discrete setting, we examine the discrete curl of the spread force at the start of the simulation, before marker positions are significantly perturbed. Fig. \ref{fig:curl_errors} demonstrates that both composite B-spline kernels and the DFIB method exhibit convergent behavior under Lagrangian grid refinement, with rates improving with kernel regularity. As expected, this dependence on regularity reflects the properties of our use of the periodic trapezoidal rule used to discretize the force spreading operation --- smoother integrands yield higher-order accuracy in the quadrature approximation. For a $C^n$ kernel, the theoretical convergence rate scales as $\Delta s^{2+n}$\cite{boyd2001}, matching our observed results. The \IBfour~kernel, however, maintains persistent $\mathcal{O}(h^{-2})$ errors in the discrete curl. Further analysis regarding spurious flows and vorticity associated with spreading forces using isotropic regularized delta functions is provided in section \ref{sec:spur_flow_iso} of the appendix. \par

\begin{figure}[t]
    \centering
    \includegraphics[width=0.475\textwidth]{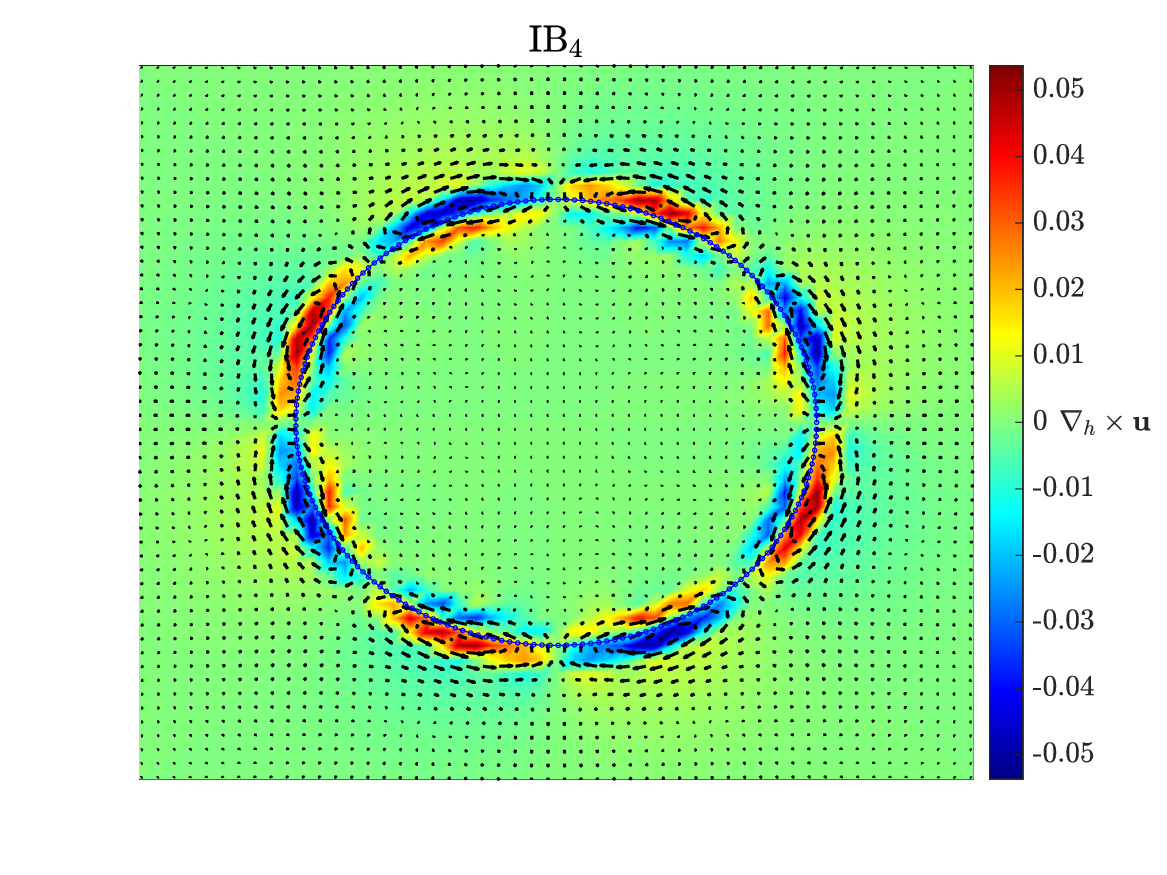}
    \includegraphics[width=0.475\textwidth]{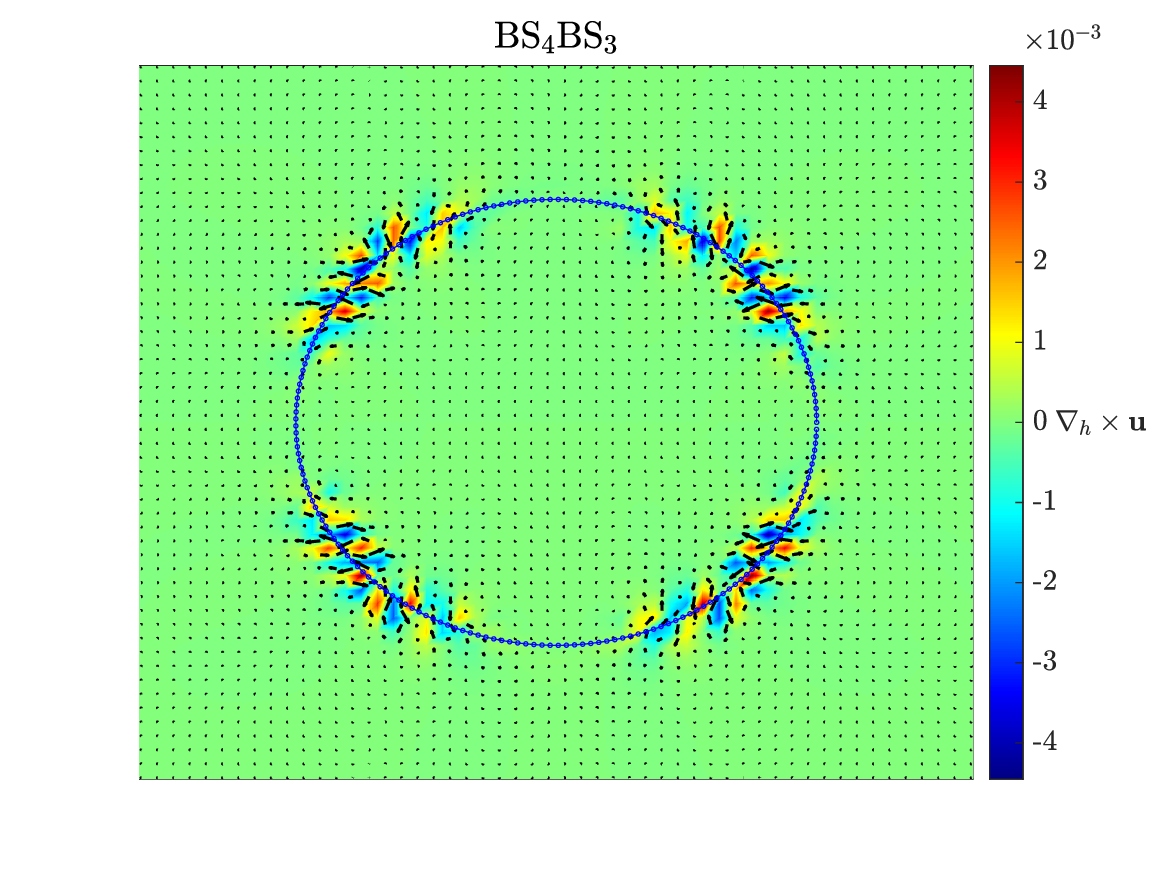}
    \includegraphics[width=0.475\textwidth]{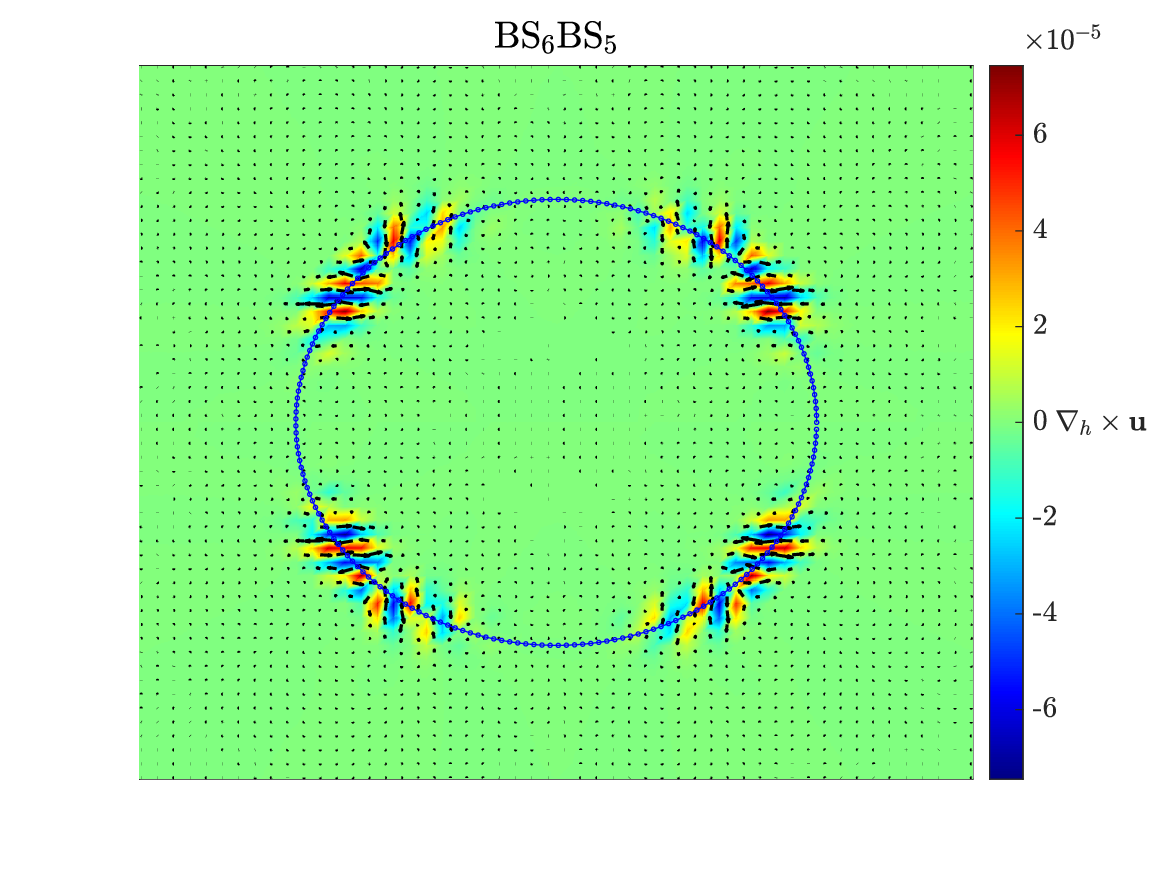}
    \includegraphics[width=0.475\textwidth]{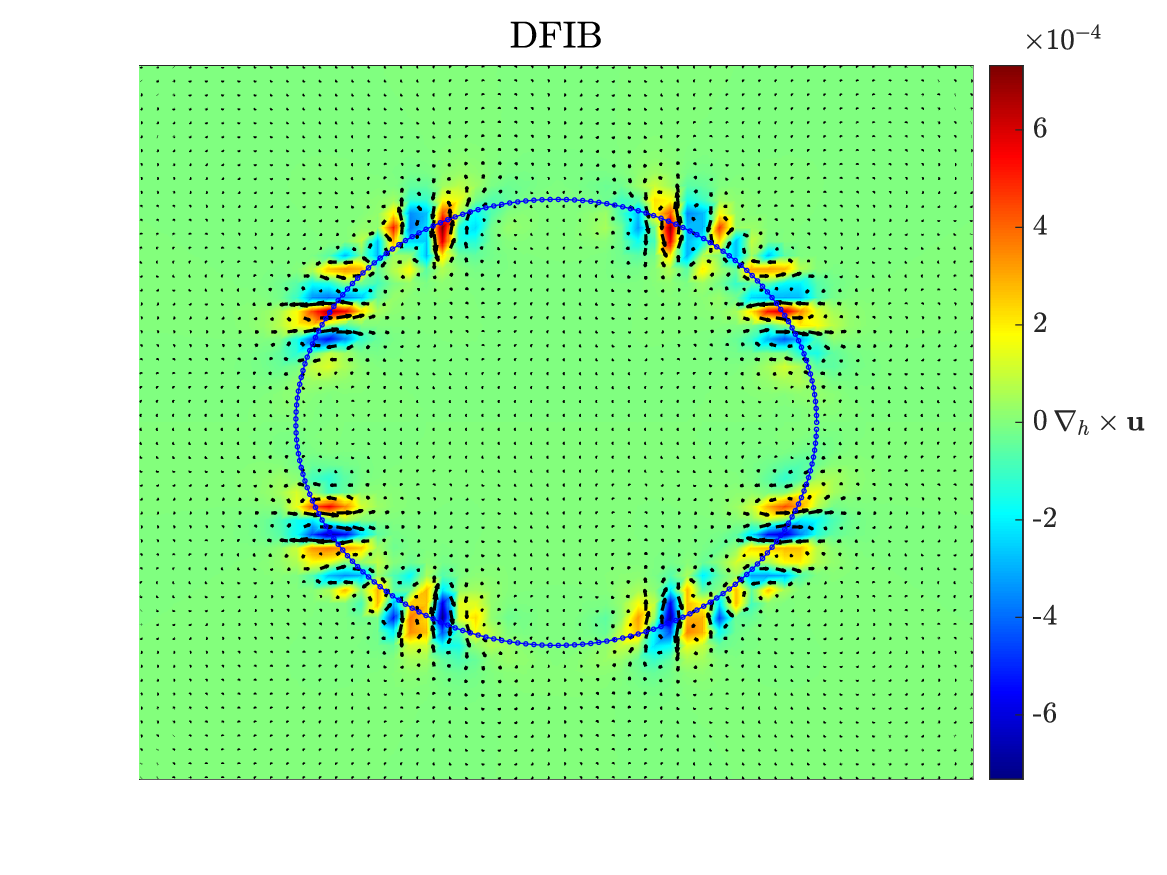}
    \caption{Flow fields at $t = 0.05$ for the IB method using the \IBfour, \BSfour\BSthree, and \BSsix\BSfive~regularized delta functions and the DFIB method. Each figure is a pseudo-color plot illustrating the vorticity $\omega = \nabla_h\times\eulVel$ and its associated vector field in black. The locations of the Lagrangian markers $\XX_k$ are plotted as blue markers. The color bar beside each plot indicates the magnitude of the vorticity.}
    \label{fig:vort_plots}
\end{figure}

These differences in force spreading manifest in the resulting flow fields. Fig. \ref{fig:vort_plots} illustrates the vorticity of the fluid velocity at $t = 0.05$ using different kernel choices. The $C^4$ \BSsix\BSfive kernel, which best maintains the gradient structure, generates minimal spurious vorticity. In contrast, the isotropic $C^1$ \IBfour kernel produces significant spurious flows distributed widely around the interface. The DFIB method and \BSfour\BSthree kernel show intermediate performance, with DFIB achieving approximately an order of magnitude reduction in spurious vorticity compared to \BSfour\BSthree.\par 

Given that our finite difference approximation of $\lagForce$ (equation \eqref{eq:discrete_lag_force}) is only second-order accurate, the close agreement between the expected periodic trapezoidal rule and observed error rates might seem surprising. However, the higher-order terms in the Taylor series expansion of the error are all proportional to even derivatives of the parameterization. Each error term is thus directly proportional to the circle's normal vector. Consequently, the divergence theorem implies that these higher-order error terms also vanish as the Lagrangian grid is refined.

\subsection{A Parametrically Excited Membrane}\label{sec:parametric_membrane}
In this section we use the IB and DFIB methods to simulate a parametrically excited membrane whose elastic stiffness varies periodic in time. The Lagrangian force density associated with the membrane is given by
\begin{align}
    \lagForce(s,t) &= \kappa(t)\frac{\partial^2\XX(s,t)}{\partial s^2},\\
    \kappa(t) &= \kappa_0\left(1 + 2\tau\sin\left(\omega_0 t\right)\right).
\end{align}
The membrane is initially configured as a perturbation of the a circle radius $r$ centered in a periodic square fluid domain $\Omega = \left[0,L\right]^2$. The initial configuration of the membrane is described by the parameterization
\begin{align}
    \XX(s,0) = L\left(\frac{1}{2} + r\left(1 + \epsilon\cos\left(ps\right)\right)\hat{\mathbf{n}}(s)\right),
\end{align}
where $\hat{\mathbf{n}}(s)$ is the unit outward pointing normal vector of the unit circle. 
This model problem, introduced by Cortez et al. \cite{cortex2004,ko2014correction}, serves as a simple representation of an active immersed material driven by a periodic forcing. Cortez et al. performed a Floquet analysis of this problem, in which the integer parameter $p > 1$ identifies the wavenumber of the Floquet mode under consideration. The parameters $\tau$ and $\omega_0$ control the amplitude and frequency of the stiffness oscillation, respectively. The parameter $\epsilon$ a small parameter used to linearize the equations of motion.\par 

Here, we consider two distinct time evolutions of the forced membrane using parameter values consistent with Bao et al.'s tests \cite{bao2017}. Guided by Cortez et al.'s Floquet analysis, Bao et al. examined two parameter combinations: one resulting in damped oscillations, and another leading to growing oscillations that are eventually stabilized by nonlinearities. These parameter values are reproduced in Table \ref{tab:membrane_params} for reference.
\begin{table}[t]
    \raggedright
    \caption{Parameter values for the forced membrane simulations}
    \centering
    \begin{tabular}{cccccccccc}
    \hline
    $\rho$ & $\mu$ & $L$ & $R$ & $\kappa_0$ & $\omega_0$ & $p$ & $\epsilon_0$ & \multicolumn{2}{c}{$\tau$} \\
    \hline
    1 & 0.15 & 5 & 1 & 10 & 10 & 2 & 0.05 & 0.4 & (damped oscillation) \\
    1 & 0.15 & 5 & 1 & 10 & 10 & 2 & 0.05 & 0.5 & (growing oscillation)\\
    \hline
    \end{tabular}
    \label{tab:membrane_params}
\end{table}

The background Cartesian grid is discretized using a uniform increment $h = \frac{L}{128}$ and the Lagrangian markers are initialized so that the physical distance separating them is roughly $h$ in the equilibrium configuration so that $M_{\text{fac}} = 1$. Similar to the static pressurized membrane problem we monitor errors in area conservation using the method described in section \ref{sec:area_measurements}. For each set of parameter values, we examine three different time step sizes: $\Delta t = \frac{h}{10}, \frac{h}{20},$ and $\frac{h}{40}$. Since we use a second-order time stepping scheme, we anticipate that the error in area computation will scale as $\mathcal{O}(\Delta t^2)$. However, this scaling may not be achieved if the force spreading approximation generates significant spurious flows, a problem exacerbated by kernels with insufficient regularity.

\begin{figure}[htpb]
    \hspace{0.5cm}
    \centering
    \includegraphics[trim={65pt 5pt 2pt 10pt}, clip, width=0.8\textwidth]{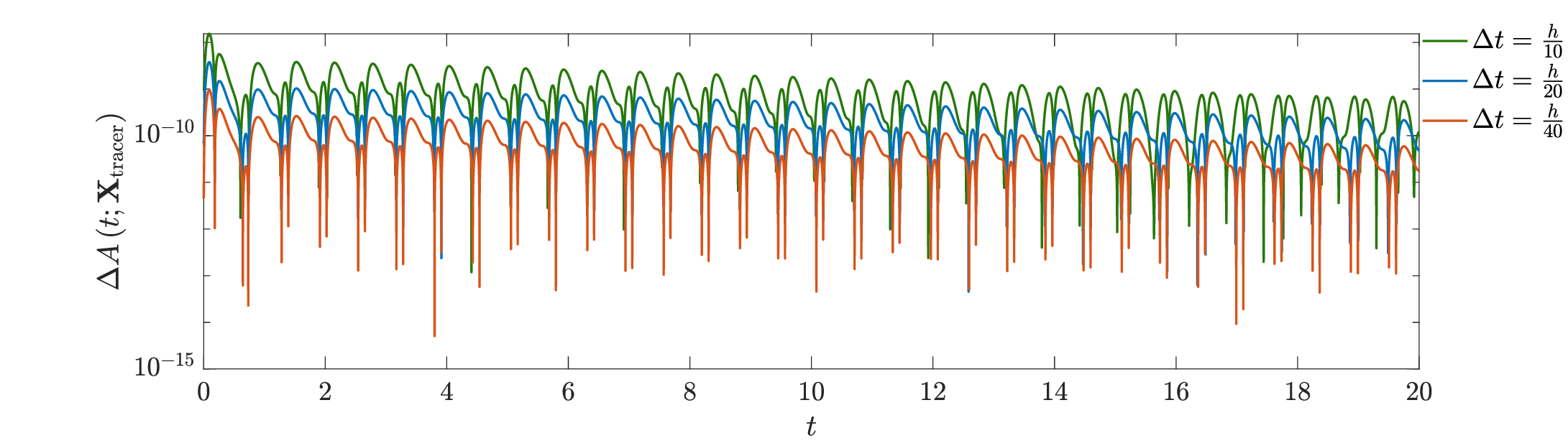}
    \caption*{(a) \BSfive\BSfour}
    \hspace{-0.9cm}
    \begin{minipage}[b]{0.495\textwidth}
     \centering
     \includegraphics[trim={65pt 5pt 10pt 10pt}, clip, width=1.2\textwidth]{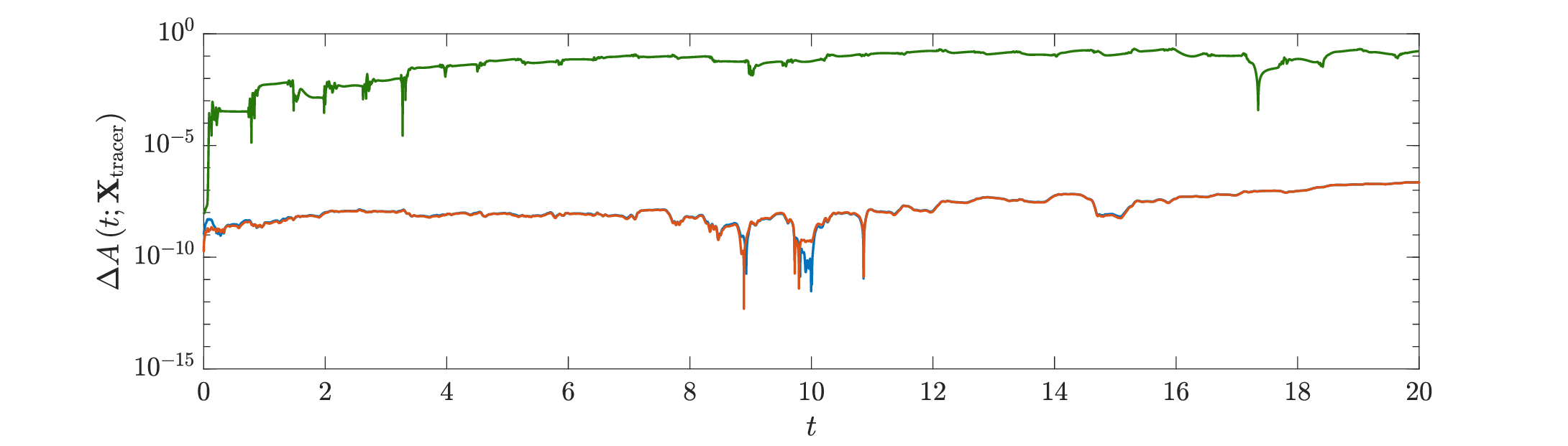}
     \caption*{(b) \BSthree\BStwo}
    \end{minipage}
    \hfill
    \begin{minipage}[b]{0.495\textwidth}
     \centering
     \includegraphics[trim={90pt 5pt 10pt 10pt}, clip, width=1.2\textwidth]{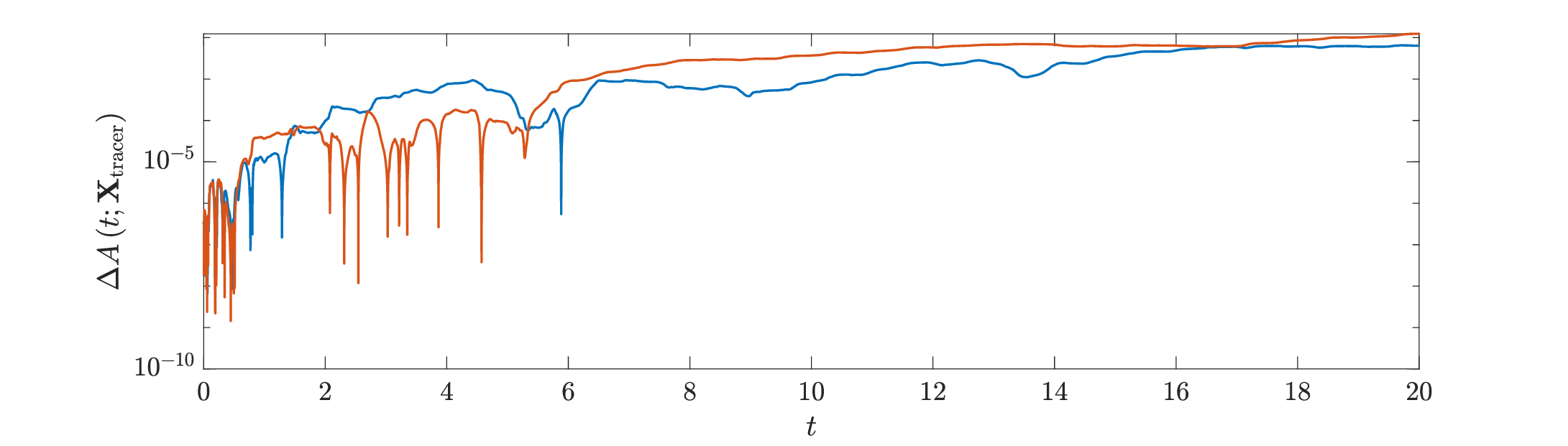}
     \caption*{(c) \BStwo\BSone}
    \end{minipage}
    
    \vspace{1em}
    \hspace{-0.9cm}
    \begin{minipage}[b]{0.495\textwidth}
     \centering
     \includegraphics[trim={65pt 5pt 10pt 10pt}, clip, width=1.2\textwidth]{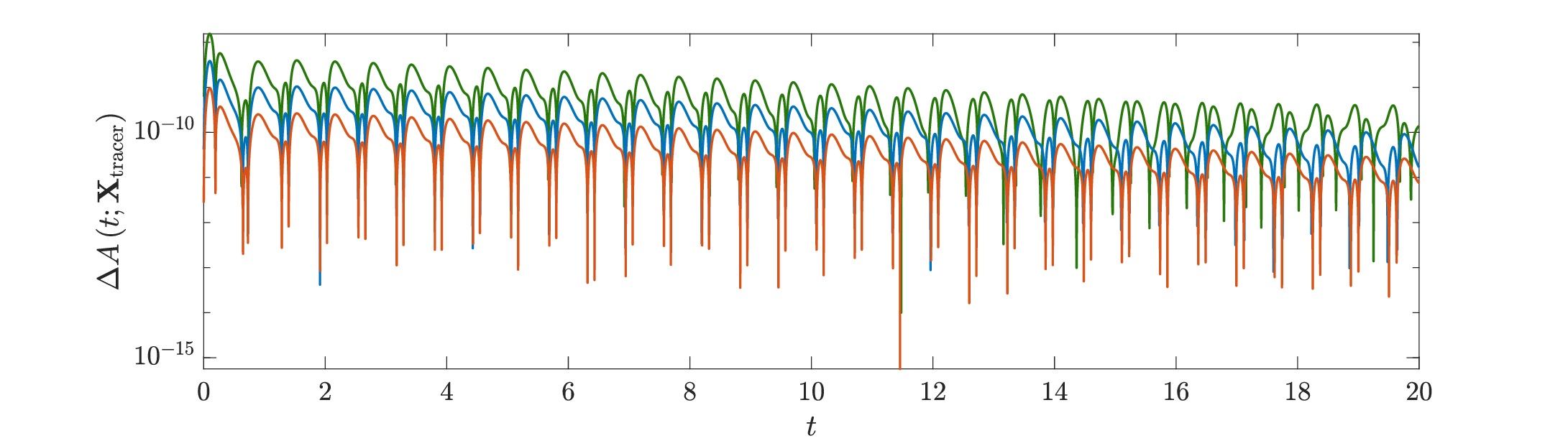}
     \caption*{(d) DFIB}
    \end{minipage}
    \hfill
    \begin{minipage}[b]{0.495\textwidth}
     \centering
     \includegraphics[trim={90pt 5pt 10pt 10pt}, clip, width=1.2\textwidth]{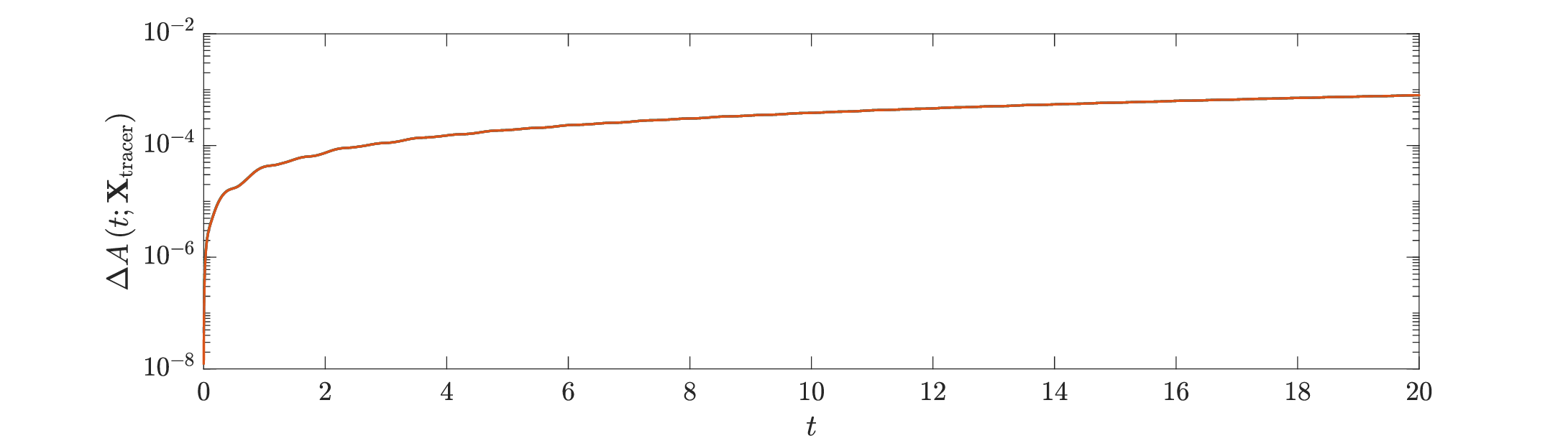}
     \caption*{(e) \IBfour}
    \end{minipage}
    
    \caption{Relative area errors over time for a membrane undergoing damped oscillations due to parametric excitation. time step sizes are $\Delta t = h/10$ (green), $h/20$ (blue), and $h/40$ (orange). Results for $\Delta t = h/10$ with the $\mathrm{BS}_2\mathrm{BS}_1$ kernel are omitted because the simulation became unstable. }
    \label{fig:damped_res_errors}
   \end{figure}

For each dynamic problem, we observed that area errors associated with composite B-splines of $C^1$ regularity or greater produced were roughly identical to the area errors produced by the DFIB method. Thus, to avoid redundancy, we present only the results of the \BStwo\BSone, \BSthree\BStwo,~and \BSfive\BSfour composite B-spline regularized delta functions. Fig \ref{fig:damped_res_errors} shows the relative area errors of these composite B-splines along side with the relative area errors produced by the \IBfour~kernel and DFIB method. 

\begin{figure}[htpb] 
    \hspace{0.5cm}
    \centering
    \includegraphics[trim={65pt 5pt 2pt 10pt}, clip, width=0.8\textwidth]{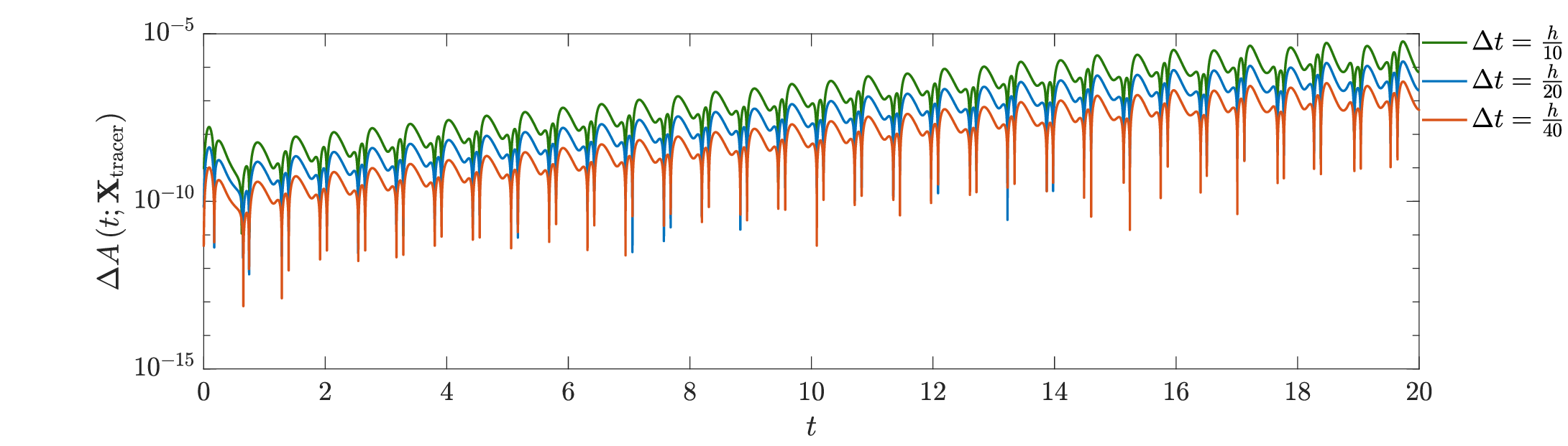}
    \caption*{(a) \BSfive\BSfour}
    \hspace{-0.9cm}
    \begin{minipage}[b]{0.495\textwidth}  
    \centering
    \includegraphics[trim={65pt 5pt 10pt 10pt}, clip, width=1.2\textwidth]{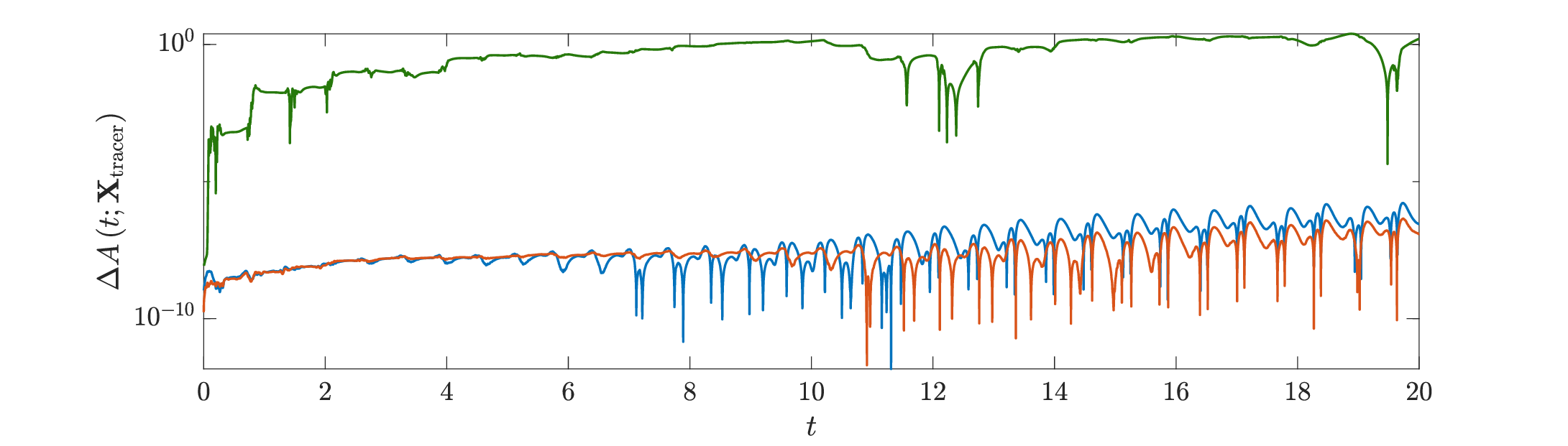}
    \caption*{(b) \BSthree\BStwo}
    \end{minipage}
    \hfill
    \begin{minipage}[b]{0.495\textwidth}
    \centering
    \includegraphics[trim={90pt 5pt 10pt 10pt}, clip, width=1.2\textwidth]{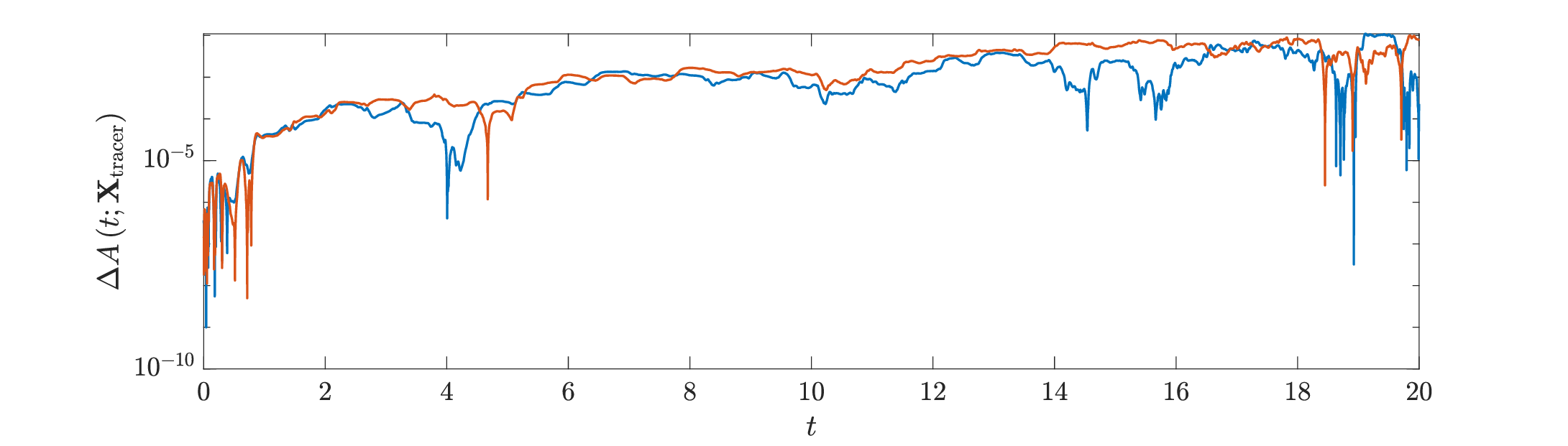}
    \caption*{(c) \BStwo\BSone}
    \end{minipage}
    
    \vspace{1.5em}  
    \hspace{-0.9cm}
    \begin{minipage}[b]{0.495\textwidth}
    \centering
    \includegraphics[trim={65pt 5pt 10pt 10pt}, clip, width=1.2\textwidth]{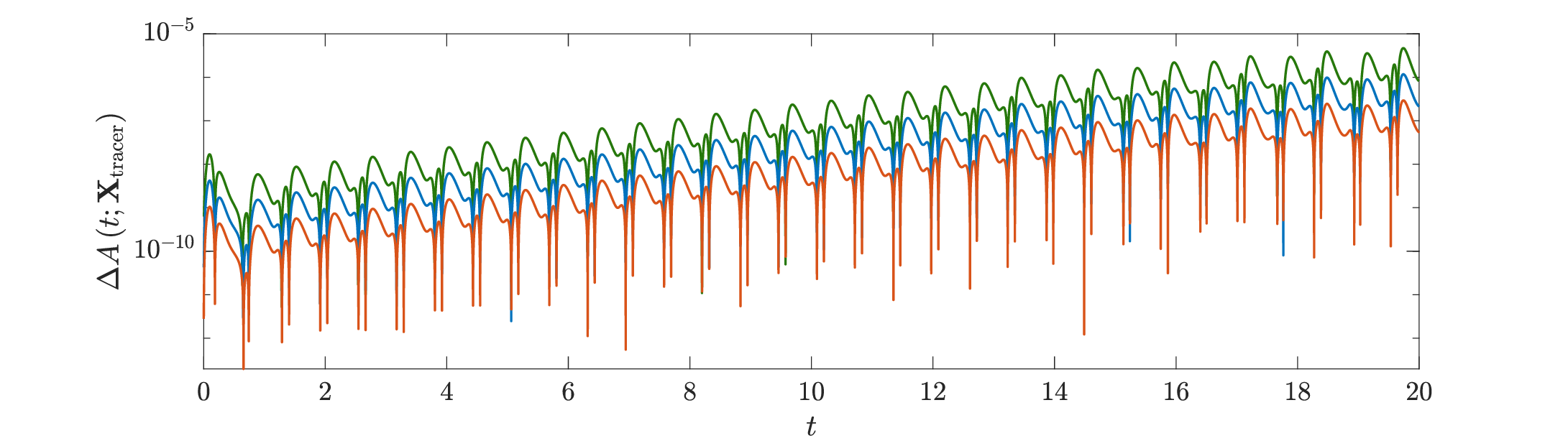}
    \caption*{(d) DFIB}
    \end{minipage}
    \hfill
    \begin{minipage}[b]{0.495\textwidth}
    \centering
    \includegraphics[trim={90pt 5pt 10pt 10pt}, clip, width=1.2\textwidth]{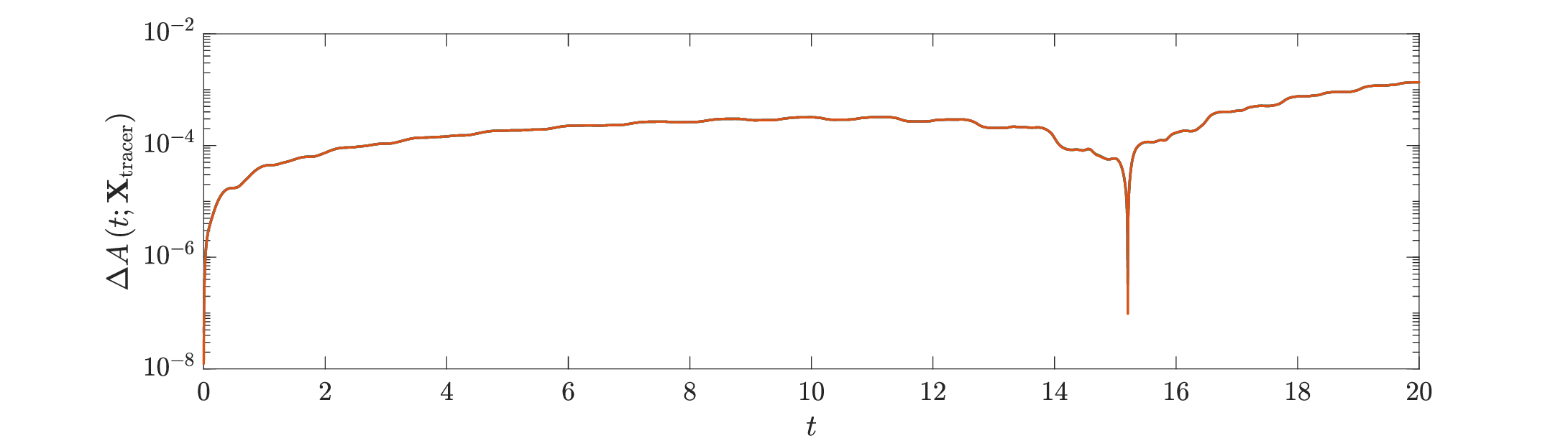}
    \caption*{(e) $\mathrm{IB}_4$}
    \end{minipage}
    
    \vspace{1em}
    \caption{Relative area errors as a function of time for a membrane undergoing growing oscillations due to parametric excitation. Each kernel function was utilized with three different choices of time step size $\Delta t = \frac{h}{10}$ (green), $\Delta t = \frac{h}{20}$ (blue), and $\Delta t = \frac{h}{40}$ (orange). The results for the time step size of $\Delta t = \frac{h}{10}$ were not included for the \BStwo\BSone~regularized delta functions because the simulation was not stable.}
    \label{fig:full_res}
\end{figure}

The area errors for the \BStwo\BSone~were not included for the time step size $\Delta t = \frac{h}{10}$ because the simulation became unstable. For composite B-spline regularized delta functions with $C^1$ regularity or higher, as well as for the DFIB method, we observe that the area errors are primarily influenced by the time-stepping scheme. This is evidenced by the consistent $\mathcal{O}\left(\Delta t^2\right)$ decrease in area error for both approaches. For the \IBfour~kernel, all of the relative area errors super impose atop each other indicating that the error is dominated by the velocity interpolation error associated with the \IBfour~kernel. The \BSthree\BStwo~composite kernel performs relatively poorly compared to the more regular composite B-splines at the largest time step size of $\Delta t = \frac{h}{10}$. When the time step size is reduced by a factor of 2, the \BSthree\BStwo~kernel produces more accurate error estimates. However, with further reduction in time step size, the errors remain constant, indicating that errors associated with force spreading approximation and lack of kernel regularity become dominant. The discontinuous \BStwo\BSone~kernel produces area errors of relatively consistent magnitude across all stable time step sizes. This consistency suggests that the lack of regularity of the \BStwo\BSone~kernel dominates the error. Recall that the \BStwo\BSone~kernel generates a discontinuous interpolated velocity field, resulting in first-order local truncation errors and thus $\mathcal{O}(1)$ global truncation errors over the simulation duration. These trends similarly reproduce for the resonant case of the parametrically excited membrane; however, the relative area errors are slightly larger due to larger time-stepping errors associated with the growing amplitude of the membrane's oscillation. These results are presented in Fig \ref{fig:full_res}.

 \FloatBarrier

\subsection{A Flexible Membrane in Lid-Driven Cavity Flow}\label{sec:lc_flow}

In this section, we use the IB method to simulate a flexible membrane immersed in lid-driven cavity flow. This example is inspired by the benchmark that Griffith and Luo explored with their then-nascent Immersed Finite Element Finite Difference (IFED) method, in which they simulated an immersed soft, neo-Hookean disc in lid-driven cavity flow~\cite{griffith2017}. We include this example to demonstrate that composite B-spline regularized delta functions effectively mitigate volume conservation errors even in the presence of physical boundary conditions. Since both the DFIB and FSIB methods, which excel at mitigating volume conservation errors, are limited to periodic boundary conditions, we focus our comparison on the IB method using composite B-spline regularized delta functions versus the standard four-point \IBfour~ and six-point \IBsix~regularized delta functions.

The computational domain is the unit square $[0,1]\times[0,1]$, discretized using a uniform Cartesian grid with meshwidth $h = \frac{1}{128}$. The Lagrangian force density is modeled using equation \eqref{eq:discrete_lag_force} with stiffness $\kappa = 1.0$. The membrane is initialized in its equilibrium configuration as a circle of radius $r = 0.2$ centered at $(0.6,0.5)$, with mesh factor $M_{\text{fac}} = \frac{1}{2}$.

The fluid is initialized to zero velocity. Following Griffith and Luo~\cite{griffith2017}, the horizontal component of velocity along the top wall is set to unity, $u(x,1) = 1.0$, while both velocity components satisfy homogeneous Dirichlet boundary conditions along the remaining walls.

For these parameter values, we utilize three different timestep sizes: $\Delta t = \frac{h}{4}$, $\frac{h}{8}$, and $\frac{h}{16}$ to investigate how volume conservation errors depend on timestep size. The simulation is run until a final time of $t = 10$. Since our temporal discretization is second-order accurate, we expect volume conservation errors to be dominated by timestepping errors and thus to decrease by approximately a factor of four each time the timestep size is halved. 

\begin{figure}[htbp]
\centering
\begin{subfigure}{0.46\textwidth}
    \centering
    \includegraphics[width=\textwidth]{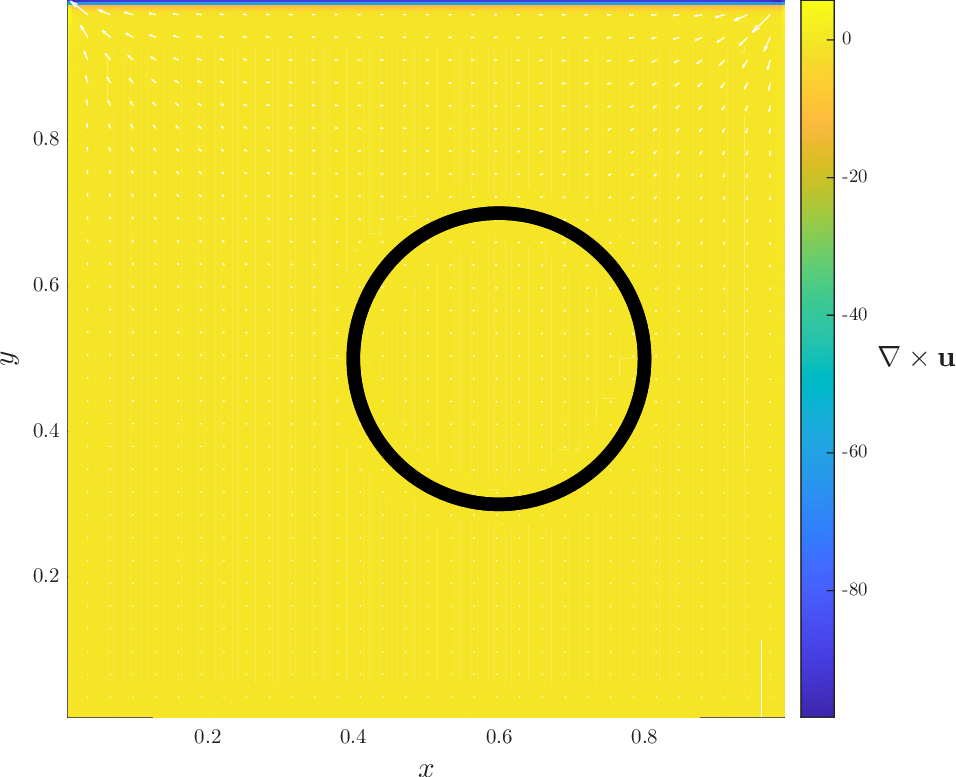}
    \caption{$t = 0$}
\end{subfigure}
\hfill
\begin{subfigure}{0.46\textwidth}
    \centering
    \includegraphics[width=\textwidth]{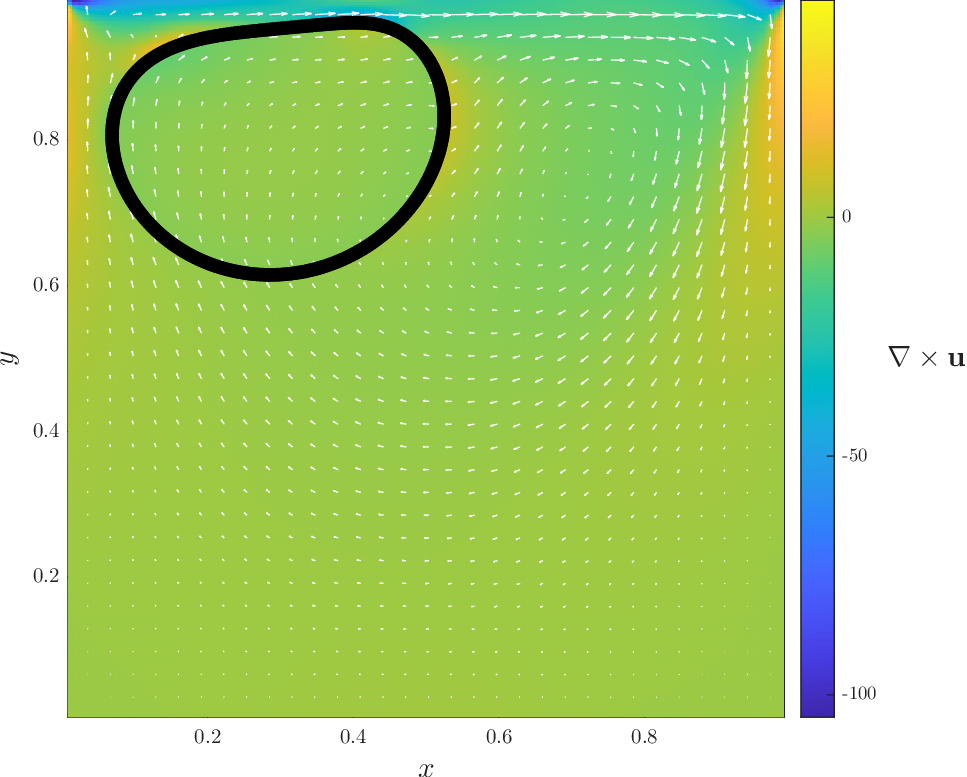}
    \caption{$t = 3$}
\end{subfigure}

\vspace{0.5cm}

\begin{subfigure}{0.46\textwidth}
    \centering
    \includegraphics[width=\textwidth]{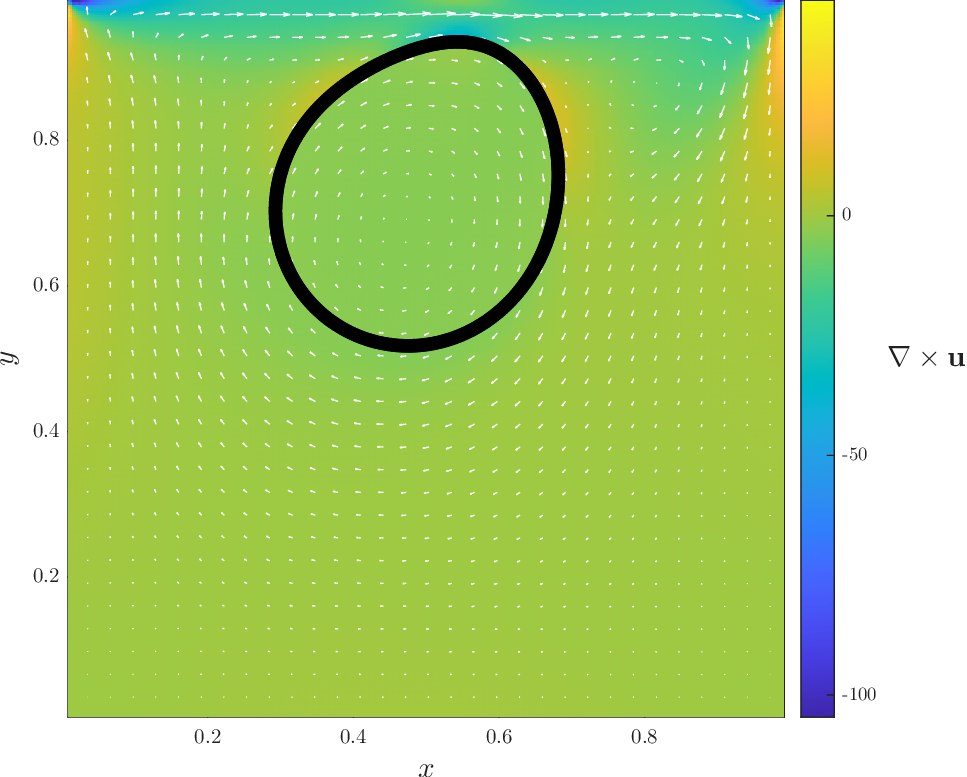}
    \caption{$t = 5$}
\end{subfigure}
\hfill
\begin{subfigure}{0.46\textwidth}
    \centering
    \includegraphics[width=\textwidth]{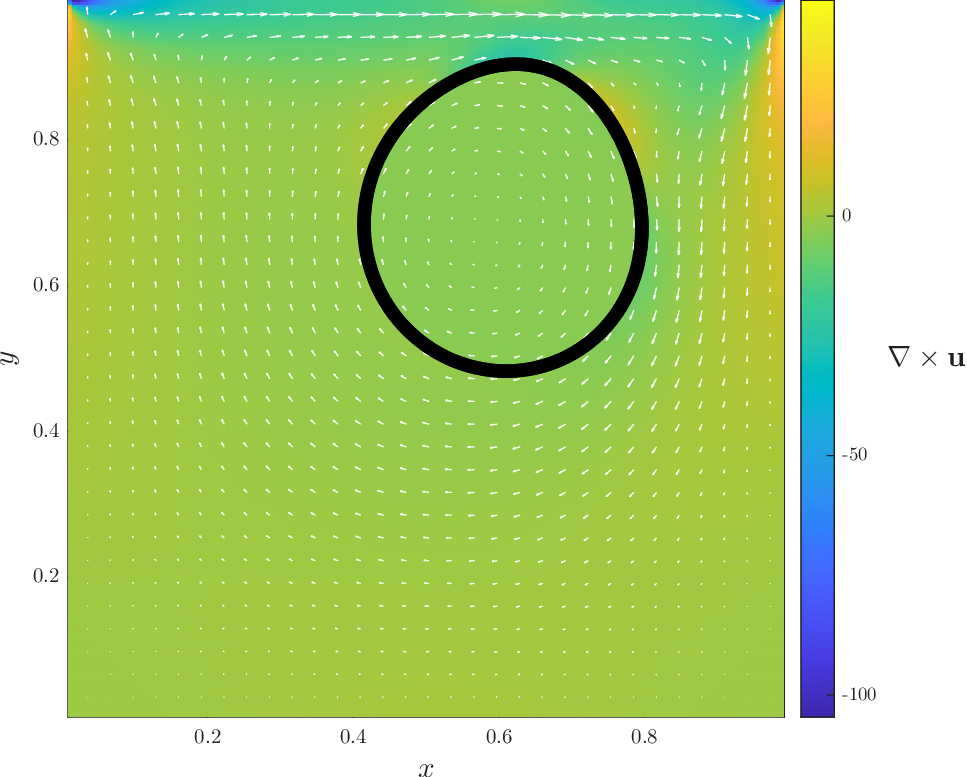}
    \caption{$t = 7$}
\end{subfigure}
\caption{Snapshots of membrane dynamics in lid-driven cavity flow over the time interval $0 \le t \le 7$.}
\label{fig:membrane_snapshots}
\end{figure}

Figure \ref{fig:membrane_snapshots} illustrates snapshots of the membrane dynamics at $t = 0$, $3$, $5$, and $7$ using the \BSfour\BSthree~regularized delta function. Over the course of the simulation the membrane becomes entrained in the flow field and migrates from its starting position towards the top of the cavity where it undergoes a relatively large deformation as it is compressed against the lid of the cavity. This near contact is automatically handled by the IB formulation using the approach outlined in the appendix of Kallemov et al.~\cite{kallemov2016}. 

\begin{figure}[htbp]
 \centering
 \includegraphics[width=0.95\textwidth]{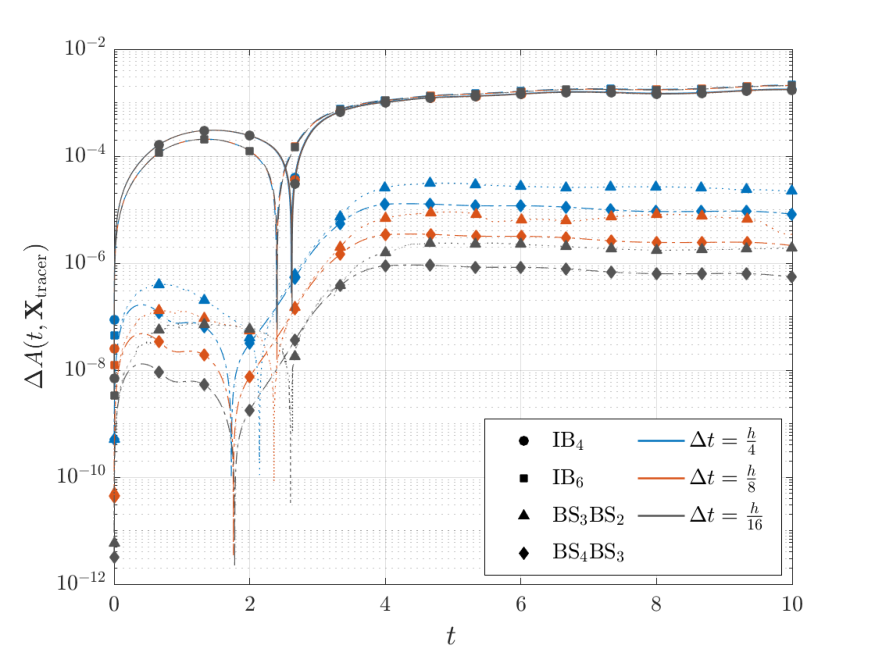}
 \caption{Relative area errors for the \IBfour~, \IBsix~, and \BSthree~\BStwo~ regularized delta functions across three timestep sizes: $\Delta t = \frac{h}{4}$, $\frac{h}{8}$, and $\frac{h}{16}$. The composite B-spline regularized delta functions exhibit volume conservation errors that decrease as $\mathcal{O}(\Delta t^2)$ with timestep refinement, while the \IBfour~ and \IBsix~ regularized delta functions show essentially identical behavior across all timestep choices.}
 \label{fig:area_errors_lc}
 \end{figure}

Figure \ref{fig:area_errors_lc} shows the relative area errors associated with the \IBsix~, \IBfour~, \BSthree\BStwo~, and \BSfour\BSthree~ regularized delta functions. The \IBfour~ and \IBsix~ regularized delta functions display relative area errors that essentially overlap with one another for each choice of timestep size. The \BSthree\BStwo~ and \BSfour\BSthree~ regularized delta functions exhibit relative area errors that decrease as $\mathcal{O}(\Delta t^2)$ as the timestep size is refined, indicating that for these kernels the volume conservation errors are dominated by timestepping errors. We did not display results for the higher-order composite B-splines because the results were essentially identical to those of the \BSfour\BSthree~ regularized delta function, indicating that additional regularity of the integrand in the force spreading operator no longer provides benefits since the timestepping errors are dominate. Similar to the parametrically excited membrane example (Section \ref{sec:parametric_membrane}), the discontinuous \BStwo~\BSone~ regularized delta function eventually yields $\mathcal{O}(1)$ relative area errors. This degradation occurs because the kernel's discontinuities reduce the accuracy of the temporal discretization when marker points begin to interact with new Eulerian degrees of freedom, effectively reducing what should be a second-order time integration scheme to one with first order local truncation errors as analyzed in Appendix \ref{sec:error_analysis_ts}. For this reason, we did not include results for the \BStwo\BSone~ regularized delta function in this example.

\FloatBarrier

\section{Conclusions}
In this study, we have introduced and analyzed the use of composite B-spline regularized delta functions to improve the volume conservation properties of the IB method. These B-spline kernels provide locally divergence-free velocity interpolants and effective force spreading operators. Composite B-spline regularized delta functions address a long-standing challenge associated with the IB method: the poor volume conservation of closed immersed structures, particularly evident in simulations of pressurized, closed membranes which exhibit volume loss at a rate proportional to the pressure jump across the immersed boundary.\cite{peskin1993,griffith2012} \par

Our numerical experiments demonstrate that composite B-spline regularized delta functions significantly enhance volume conservation. One mathematical insight that emerges from our analysis is that composite B-spline regularized delta functions, via their derivative property \eqref{eq:deriv_prop}, spread Lagrangian force densities which correspond continuous gradients to discrete gradients on the staggered grid. In particular, when the Lagrangian force density contains a zero-mean component oriented normal to the boundary, which in the continuous setting represents the distributional gradient of a function with jump discontinuity at the boundary, this component is spread (up to quadrature errors) to a discrete gradient field that is properly balanced by the pressure Lagrange multiplier. In contrast, isotropic kernels project a significant portion of such force onto the space of discretely divergence-free vector fields, which ultimately drives spurious currents and persistent volume loss even when the fluid velocity should be at steady-state and quiescent.

We observed for tests concerning the pure advection of Lagrangian markers (Section \ref{sec:advect_test_sec}), composite B-spline regularized delta functions of $C^0$ or higher regularity achieved area errors determined solely by second-order time discretization, matching the DFIB method's performance. While the \IBfour~kernel demonstrated second-order convergence for relatively large timestep sizes, further timestep refinement revealed that errors from the non-divergence-free interpolated velocity field ultimately dominate area conservation. The discontinuous composite B-spline, \BStwo\BSone, similarly shows initial second-order convergence in timestep size, but its convergence rate eventually degrades to first order. This degradation, though initially surprising given that \BStwo\BSone provides a continuously divergence-free interpolant, occurs because the discontinuity in the kernel reduces the traditionally second-order midpoint rule to a first-order method. A detailed analysis of this phenomenon is provided in Appendix \ref{sec:error_analysis_ts}.

Concerning the quasi-static pressurized membrane problem (Section \ref{sec:pres_membrane}), the more regular composite B-splines kernels, \BSfour\BSthree, \BSfive\BSfour, and \BSsix\BSfive, maintained the initial area to within machine precision, similar to the DFIB method. Less regular splines achieved better area conservation compared to the \IBfour~regularized delta function, which produces a linear rate of area loss in time\cite{peskin1993,griffith2012}. Our analysis revealed that these differences in performance stem from the approximation of the force spreading operator. By examining the discrete curl of the spreading operator, we analyzed how much of the spread force projects onto the space of discretely divergence-free vector fields. These divergence-free components of the spread force directly drive spurious flows at steady state. Both the IB method with composite B-spline regularized delta functions and the DFIB method show that this problematic projection converges to zero at rates consistent with periodic trapezoidal rule estimates, with more regular kernels achieving faster convergence and consequently generating smaller spurious velocities. For isotropic kernels like \IBfour, the discrete curl of the force spreading operator at steady state contains persistent $\mathcal{O}(h^{-2})$ errors, generating spurious flows and sustained volume conservation errors. These force spreading properties also impact the accuracy of computed Lagrangian forces. Both the DFIB method and IB method with composite B-spline regularized delta functions produced pointwise convergent Lagrangian forces, with errors decreasing for more regular kernels. In contrast, the \IBfour~kernel failed to achieve convergent Lagrangian force densities even under simultaneous temporal and spatial grid refinement.

In the dynamic simulations of parametrically excited membranes (Section \ref{sec:parametric_membrane}), we observed that composite B-splines of regularity class $C^1$ or higher produced results essentially identical to the DFIB method. For these higher-regularity kernels, errors were primarily dominated by time-stepping errors rather than spatial discretization errors. Lower-regularity composite B-splines were less competitive, as they provided poorer approximations of the force spreading operator. Not surprisingly, the discontinuous \BStwo\BSone~kernel performed the worst, yielding the poorest volume conservation estimates. Moreover, for the largest timestep size of $\Delta t = \frac{h}{10}$, simulations using this kernel became unstable. The interpolation and force spreading operator errors were dominant for the \IBfour~kernel at each choice of time-step size.

The lid-driven cavity flow example (Section \ref{sec:lc_flow}) further demonstrates the robustness of our approach by showing that composite B-spline regularized delta functions effectively mitigate volume conservation errors even in the presence of physical boundary conditions. This capability extends the applicability of our method beyond the periodic boundary conditions required by both the DFIB and FSIB methods. The results confirmed the expected $\mathcal{O}(\Delta t^2)$ convergence in volume conservation errors with timestep refinement, indicating that timestepping errors dominated over spatial discretization errors when using composite B-spline kernels. This example underscores the practical advantages of our approach, as it maintains the computational efficiency of the standard IB method while providing enhanced volume conservation across a broader range of boundary condition types than existing alternatives.

Our findings indicate that composite B-spline regularized delta functions require a minimum level of regularity to be competitive with the DFIB method developed by Bao et al \cite{bao2017}. For practical implementations of the IB method, we recommend the $C^1$ \BSfour\BSthree~kernel. This kernel represents the least regular B-spline with the smallest set of support that consistently performed comparably to the DFIB method across all our tests. However, we note that the dependence on B-spline regularity is perhaps largely a consequence of our use of the periodic trapezoidal rule to approximate the line integral associated with the force spreading operator. An important direction for future work is to analyze the performance of the method using quadrature rules that are less sensitive to the regularity of the integrand. This would allow force spreading operations to be made more accurate for the less regular composite B-spline kernels, potentially offering a better balance between accuracy and computational efficiency.\par 

Implementing the IB method using composite B-spline regularized delta functions offers a compelling alternative to existing approaches for improving volume conservation in the IB method. It achieves greatly improved volume conservation properties without the computational overhead of the DFIB method, making it particularly attractive for large-scale, three-dimensional simulations where computational efficiency is crucial. Importantly, adopting this approach requires only a single modification to existing IB code: changing the delta function implementation. This minimal change allows users to significantly enhance volume conservation of IB simulations with minimal effort, making it an accessible improvement for a wide range of IB method applications.

\section{Acknowledgments}
Cole Gruninger is grateful for support from the Department of Defense (DoD) through the National Defense Science and Engineering Graduate (NDSEG) Fellowship Program. Boyce E. Griffith gratefully acknowledges support from the National Institutes of Health (NIH) under grant numbers NIH U01HL143336, NIH R01HL157631, and NIH R41GM136084, as well as the National Science Foundation (NSF) under grant numbers OAC 1652541 and OAC 1931516. This manuscript is the result of funding in whole or in part by the National Institutes of Health (NIH). It is subject to the NIH Public Access Policy. Through acceptance of this federal funding, NIH has been given a right to make this manuscript publicly available in PubMed Central upon the Official Date of Publication, as defined by NIH.

% main text bibliography:
\renewcommand\refname{REFERENCES}
\bibliographystyle{unsrt}
\bibliography{refs}
\newpage
\appendix
\section{Appendix}
\subsection{Analysis of spurious flows associated with isotropic kernels}\label{sec:spur_flow_iso}
In section \ref{sec:pres_membrane}, we observed that the force density spread by the isotropic regularized delta function constructed using the \IBfour kernel generated significant spurious flows even when the Lagrangian force density was conservative. More specifically, our analysis focused on a pressurized membrane at equilibrium, where the Lagrangian force density is given by $\lagForce = \kappa r\hat{\mathbf{n}}(s)$, with $\kappa$ representing the membrane's stiffness and $r$ the radius of the membrane's circular equilibrium configuration. We observed that using isotropic regularized delta functions to spread this force resulted in substantial spurious flows near the circle's boundary.\par 
In the subsequent discussion, we posited that the induced vorticity resulting from the spread force does not converge under grid refinement and scales proportionally to $\mathcal{O}\parens{\frac{\kappa}{\mu}}$, where $\mu$ is the fluid's dynamic viscosity. This section of the appendix aims to provide a heuristic argument supporting this claim and present additional empirical evidence, further elaborating on our findings for isotropic regularized delta functions.\par 
To analyze these spurious flows, we work in the context of the equilibrium pressurized membrane problem. We once again examine the discrete curl of the force spreading operator, assuming the line integral associated with it is exact. This approach is chosen because, theoretically, for no flow to be induced, the discrete curl of the force spreading operator should vanish assuming the resulting force being spread is indeed a conservative vector field. Taking the discrete curl of the force spreading operator when using an isotropic regularized delta function yields: 
\begin{align}
    \nabla_h\times\eulForce_{i,j} = \frac{\kappa}{h^3}\int_0^{2\pi}\hat{\mathbf{n}}(s)\cdot\begin{bmatrix}
        -\phi\left(\frac{x_{i-\frac{1}{2}} - X(s)}{h}\right)\left(\phi\parens{\frac{y_{j}-Y(s)}{h}} - \phi\parens{\frac{y_{j-1}-Y(s)}{h}}\right) \\
        \phi\left(\frac{y_{j-\frac{1}{2}} - Y(s)}{h}\right)\left(\phi\parens{\frac{x_{i}-X(s)}{h}} - \phi\parens{\frac{x_{i-1}-X(s)}{h}}\right)
    \end{bmatrix}
    \,r\text{d}s.
\end{align}
Applying the divergence theorem to this equation, we can see that the the curl of the Eulerian force density is given by 
\begin{align}
    \nabla_h\times\eulForce_{i,j} = -\frac{\kappa}{h^4}\iint &\left[\left(\phi\parens{\frac{x_{i}-X}{h}} - \phi\parens{\frac{x_{i-1}-X}{h}}\right)\phi^{\prime}\parens{\frac{y_{j-\frac{1}{2}} - Y}{h}} \right. \notag \\
    &\left. - \left(\phi\parens{\frac{y_{j}-Y}{h}} - \phi\parens{\frac{y_{j-1}-Y}{h}}\right)\phi^{\prime}\left(\frac{x_{i-\frac{1}{2}} - X}{h}\right)\right]\,\text{d}X\text{d}Y,\\
    = -\frac{\kappa}{h^4}\iint &\left[\left(\phi\parens{\frac{x_{i-\frac{1}{2}}-X}{h}+ \frac{1}{2}} - \phi\parens{\frac{x_{i-\frac{1}{2}}-X}{h}- \frac{1}{2}}\right)\phi^{\prime}\parens{\frac{y_{j-\frac{1}{2}} - Y}{h}} \right. \notag \\
    &\left. - \left(\phi\parens{\frac{y_{j-\frac{1}{2}}-Y}{h}+ \frac{1}{2}} - \phi\parens{\frac{y_{j-\frac{1}{2}}-Y}{h}-\frac{1}{2}}\right)\phi^{\prime}\left(\frac{x_{i-\frac{1}{2}} - X}{h}\right)\right]\,\text{d}X\text{d}Y,
\end{align}
in which the double integral above is taken of the interior of the circle. Clearly, if the central differences of the kernels, e.g., $\phi\parens{\frac{y_{j-\frac{1}{2}}-Y}{h}+ \frac{1}{2}} - \phi\parens{\frac{y_{j-\frac{1}{2}}-Y}{h}-\frac{1}{2}}$, were equivalent to their true derivatives or even good approximations of those derivatives, the curl of the Eulerian force density would be zero or relatively small, respectively. However, this is not the case. These central difference approximations generally contain errors which are $\mathcal{O}(1)$ even when the kernel is smooth. As a result, errors in the curl of the Eulerian force density will be supported along the circle and will be of size $\mathcal{O}\left(\frac{\kappa}{h^2}\right)$. This scaling can be inferred from the fact that the integrand is supported on a set with area proportinal to $h^2$.\par 
Because the discrete curl of the Eulerian force density scales like $\mathcal{O}\left(\frac{\kappa}{h^2}\right)$, an order of magnitude analysis indicates that the resulting vorticity ought to scale like $\mathcal{O}\left(\frac{\kappa}{\mu}\right)$. This is because, at steady-state, assuming the nonlinear convective term is small compared to the other terms, taking the discrete curl of the equations tells us 
\begin{align}
\mu\Delta_h\omega \sim -\nabla_h\times \eulForce_{i,j}.
\end{align}
Since the $\Delta_h$ operator is of magnitude $\mathcal{O}\left(\frac{1}{h^2}\right)$, we expect $\omega$ to scale like $\mathcal{O}\left(\frac{\kappa}{\mu}\right)$.\par 
Similarily, we can determine the scaling of the velocity field. Since $\eulVel_{i,j}$ may be obtained by solving $\Delta_h\eulVel_{i,j} = \nabla^{\perp}_h\omega$, we expect the spurious velocity to scale like $\mathcal{O}\left(h\frac{\kappa}{\mu}\right)$. Consequently, the velocity is still expected to converge under grid refinement. A result which was proven by Mori for the IB method applied to Stokes flow\cite{mori2008}.\par 
In the remainder of this section, we provide empirical evidence confirming our claims regarding the magnitude of the vorticity and spurious velocity induced by using isotropic regularized delta functions. To test the heuristic analysis above, we use two isotropic regularized delta functions to spread the equilibrium Lagrangian force density $\lagForce = -\kappa r\hat{\mathbf{n}}(s)$ onto the background grid and measure the maximum vorticity and spurious flow induced by the spread force. The isotropic delta functions we employ are constructed using the \IBfour~kernel function and the ``Gaussian-like" \IBsix~kernel function \cite{bao2017}. Each test uses the same numerical implementation as discussed in section \ref{sec:numerical implementations} and each simulation was run until a final time of $0.05$ using a Lagrangian grid spacing corresponding to $M_{\text{fac}} = \frac{1}{8}$. We've observed that the velocity and vorticity do not change substantially on longer timescales. \par
\begin{figure}[h]
    \centering
    \begin{minipage}[b]{0.495\textwidth}
     \centering
     \includegraphics[width=\textwidth]{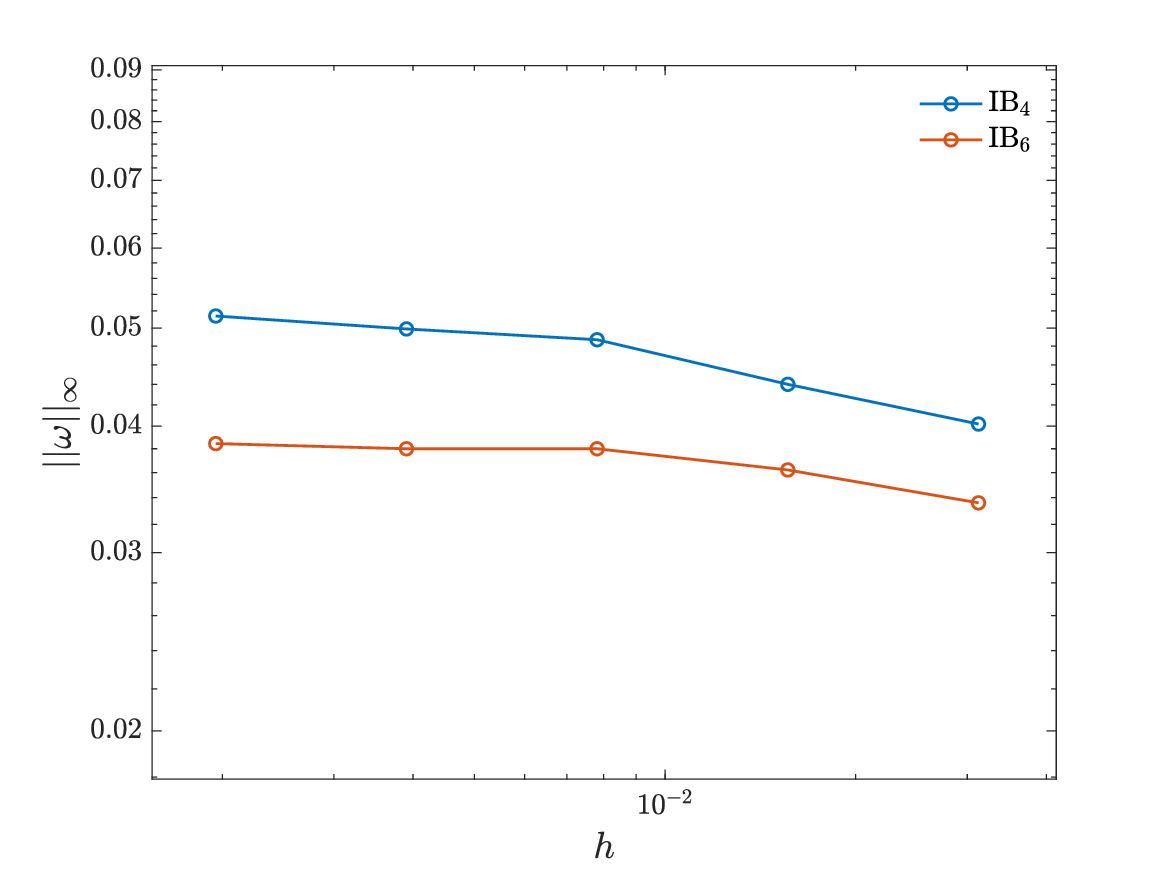}
    \end{minipage}
    \hfill
    \begin{minipage}[b]{0.495\textwidth}
     \centering
     \includegraphics[width=\textwidth]{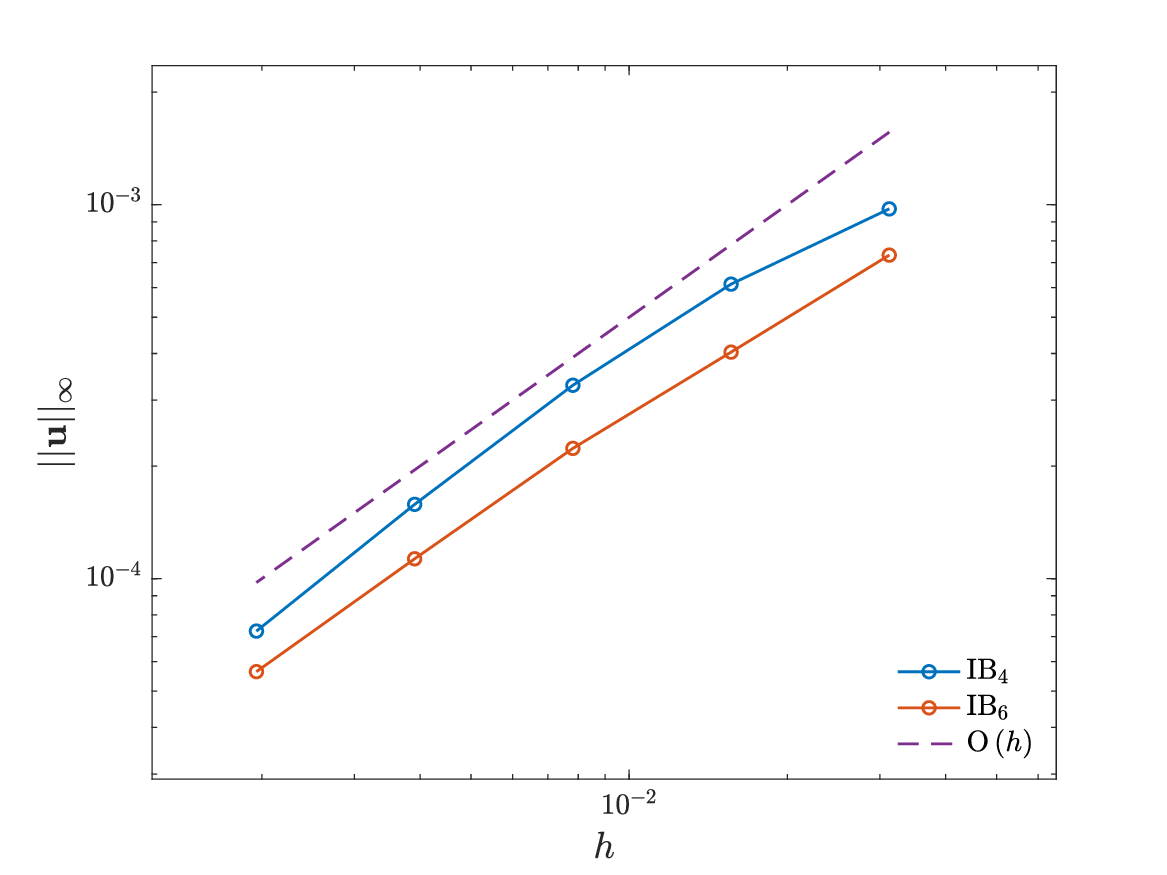}
    \end{minipage}
    \caption{Dependence of the magnitude of spurious velocity and vorticity on the background mesh width $h$. The membrane stiffness $\kappa$ and fluid viscosity $\mu$ are kept constant at values of $\kappa = 1$ and $\mu = 0.1$, respectively.}
    \label{fig:h_dependence}
\end{figure}
\begin{figure}[h]
    \centering
    \begin{minipage}[b]{0.495\textwidth}
     \centering
     \includegraphics[width=\textwidth]{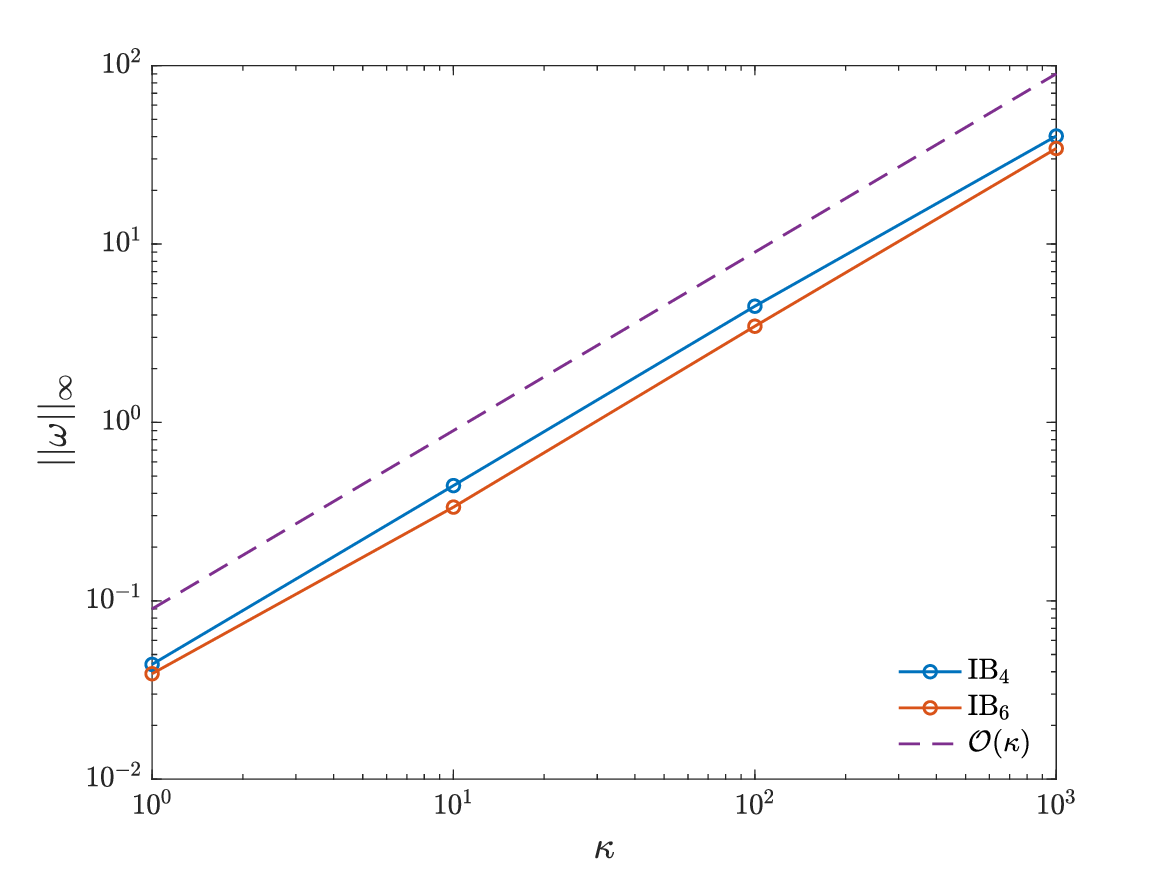}
    \end{minipage}
    \hfill
    \begin{minipage}[b]{0.495\textwidth}
     \centering
     \includegraphics[width=\textwidth]{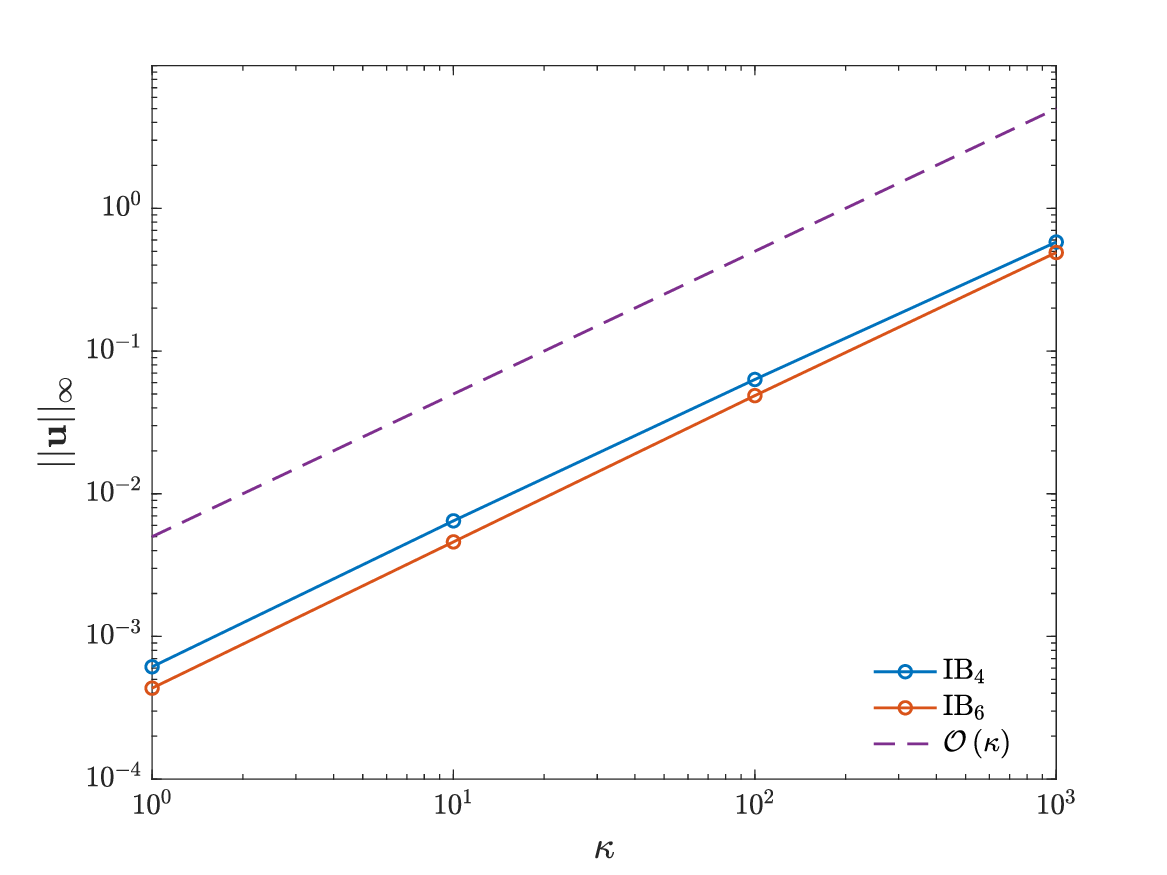}
    \end{minipage}
    \caption{Dependence of the magnitude of spurious velocity and vorticity on the value of $\kappa$. The background meshwidth $h$ and fluid viscosity $\mu$ are kept constant at values of $h = \frac{1}{64}$ and $\mu = 0.1$, respectively.}
    \label{fig:kappa_dependence}
\end{figure}
\begin{figure}[h]
    \centering
    \begin{minipage}[b]{0.495\textwidth}
     \centering
     \includegraphics[width=\textwidth]{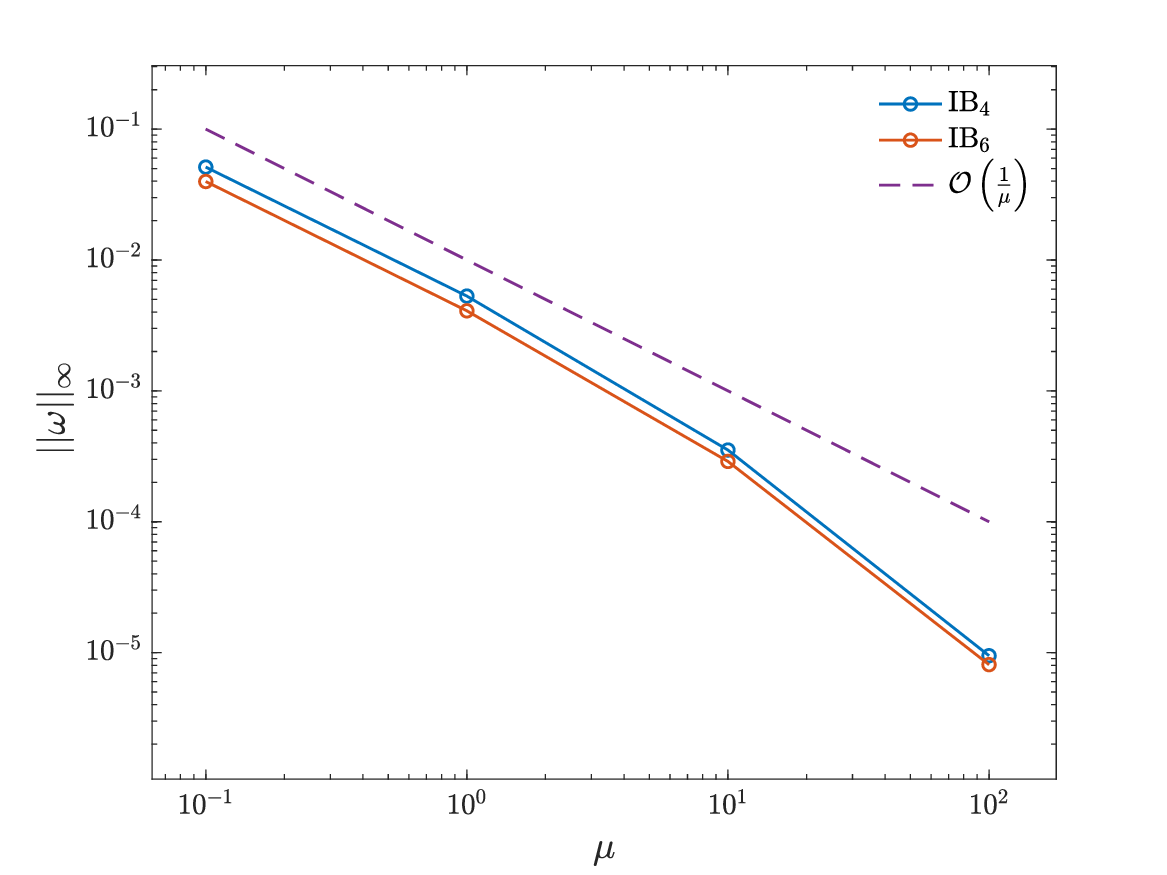}
    \end{minipage}
    \hfill
    \begin{minipage}[b]{0.495\textwidth}
     \centering
     \includegraphics[width=\textwidth]{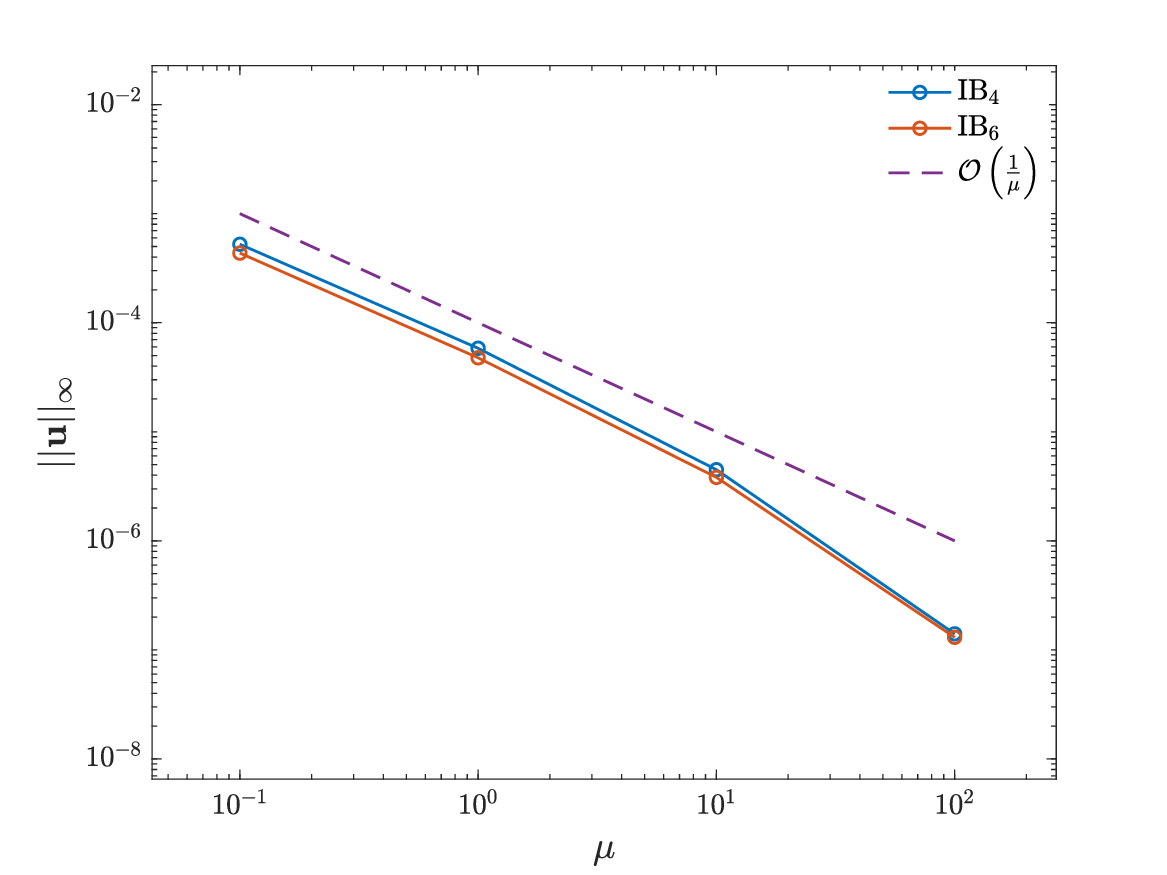}
    \end{minipage}
    \caption{Dependence of the magnitude of spurious velocity and vorticity on the value of $\mu$. The background meshwidth $h$ and membrane stiffness $\kappa$ are kept constant at values of $h = \frac{1}{64}$ and $\kappa = 1$, respectively.}
    \label{fig:mu_dependence}
\end{figure}

Figures \ref{fig:h_dependence}, \ref{fig:kappa_dependence}, and \ref{fig:mu_dependence} illustrate the magnitude of vorticity $\omega$ and velocity $\mathbf{u}$ as a function of the parameters $h$, $\kappa$, and $\mu$, respectively. In all cases, we find that the scalings predicted by our heuristic analysis above are consistent with the simulation results.

\FloatBarrier

\subsection{Error analysis of explicit time integration methods for the discontinuous velocity fields generated by the $\text{BS}_2\text{BS}_1$ kernel} \label{sec:error_analysis_ts}
In this section of the appendix, we investigate why the discontinuous \BStwo\BSone~kernel produces area errors for the pure advection problem that exhibit $\mathcal{O}(1)$ errors after initially demonstrating second-order convergence in time. To explain this behavior, we analyze the local truncation error associated with applying the forward Euler method to pure advection problem where the velocity field is interpolated using the \BStwo\BSone~regularized delta function. We demonstrate that the discontinuous nature of the interpolated velocity field results in a first order location truncation error of the time stepping scheme, resulting in a global truncation error which is $\mathcal{O}(1)$.\par 
The ordinary differential equations for the coordinates of a given Lagrangian marker $\XX(s,t) = \left(X(s,t),Y(s,t)\right)$ associated with the pure advection problem are:
\begin{align}
\label{eq:X_pos}
\frac{\partial X}{\partial t}(s,t) &= U(X(s,t),Y(s,t),t) = \sum_{i,j}u_{i,j}(t)\text{BS}_2\left(\frac{x_{i-\frac{1}{2}} - X(s,t)}{h}\right)\text{BS}_1\left(\frac{y_{j} - Y(s,t)}{h}\right), \\
\label{eq:Y_pos}
\frac{\partial Y}{\partial t}(s,t) &= V(X(s,t),Y(s,t),t) = \sum_{i,j}v_{i,j}(t)\text{BS}_1\left(\frac{x_{i} - X(s,t)}{h}\right)\text{BS}_2\left(\frac{y_{j-\frac{1}{2}} - Y(s,t)}{h}\right).
\end{align}
Since the function $\text{BS}_1$ is discontinuous, the interpolated velocities $U\parens{X(s,t),Y(s,t),t}$ and $V\parens{X(s,t),Y(s,t),t}$ are discontinuous in their $Y(s,t)$ and $X(s,t)$ arguments, respectively. Consequently, the time derivatives $\frac{\partial X}{\partial t}$ and $\frac{\partial Y}{\partial t}$ are discontinuous functions. When using an explicit time stepping scheme that is neither adaptive nor explicitly designed to handle these discontinuities, the local truncation error is reduced to first order.
To illustrate this, consider the forward Euler method applied to the ordinary differential equation:
\begin{align}
\frac{\partial w}{\partial t} = f\left(w(t),t\right),
\end{align}
in which $f(w(t),t)$ is discontinuous with respect to $w(t)$. The forward Euler method at starting time $t$ is
\begin{align}
w(t+\Delta t) = w(t) + \Delta t f(w(t),t).
\end{align}
If $f$ is smooth in the interval $[t,t+\Delta t]$, the local truncation error is $\mathcal{O}(\Delta t^2)$. However, if $f(w(t),t)$ is discontinuous at $t^* \in \left(t,t+\Delta t\right)$, the local truncation error reduces to first order and is proportional to the jump in the time derivative. Assuming $f(w(t),t)$ is smooth before and after $t^*$, we can demonstrate this using a Taylor series expansion. Let $t^*_+ = t^* + \varepsilon$ and $t^*_- = t^* - \varepsilon$, where $0 < \varepsilon \ll 1$ is on the order of $\Delta t^2$ or smaller. Expanding $w(t+\Delta t)$ at $t^*_+$ and $w(t)$ at $t^*_-$, we have
\begin{align}
w(t+\Delta t) &= w(t^*_+) + \frac{\partial w}{\partial t}(t^*_+)\left(t+\Delta t - t^*\right) + \mathcal{O}(\Delta t^2), \\
w(t) &= w(t^*_-) + \frac{\partial w}{\partial t}(t^*_-)\left(t - t^*\right) + \mathcal{O}(\Delta t^2). \label{eq:prev_time}
\end{align}
Defining $[\![w]\!](t^*) = \lim_{\varepsilon\to 0^+} w(t^*_+) - w(t^*_-)$ as the jump in $w(t)$ at $t^*$, we can rewrite the expansion of $w(t+\Delta t)$ as
\begin{align}
w(t+\Delta t) = [\![w]\!](t^*) + w(t^*_-) + \left[\!\left[\frac{\partial w}{\partial t}\right]\!\right](t^*)\left(t + \Delta t - t^*\right) + \frac{\partial w}{\partial t}(t^*_-)\left(t + \Delta t - t^*\right) + \mathcal{O}(\Delta t^2). \label{eq:tay_jumps}
\end{align}
Subtracting equation \eqref{eq:prev_time} from \eqref{eq:tay_jumps}, we get:
\begin{align}
w(t+\Delta t) - w(t) = [\![w]\!](t^*) + \Delta t \frac{\partial w}{\partial t}(t^*_-) + \left[\!\left[\frac{\partial w}{\partial t}\right]\!\right](t^*)\left(t + \Delta t - t^*\right) + \mathcal{O}(\Delta t^2).
\end{align}
Taylor expanding $f(w(t),t)$ about $t^*_-$ and subtracting from the above equation, we find the local truncation error for the forward Euler method:
\begin{align}
w(t+\Delta t) - w(t) - \Delta t f(w(t),t) = [\![w]\!](t^*) + \left[\!\left[\frac{\partial w}{\partial t}\right]\!\right](t^*)\left(t + \Delta t - t^*\right) + \mathcal{O}(\Delta t^2).
\end{align}
Assuming $w$ is continuous at $t^*$, $[\![w]\!](t^*)$ vanishes, and the local truncation error becomes first order accurate. The principal error term is proportional to the jump in the derivative at $t^*$:
\begin{align}
\left[\!\left[\frac{\partial w}{\partial t}\right]\!\right] = \lim_{\varepsilon\to 0^+} f(w(t^*_+),t^*_+) - f(w(t^*_-),t^*_-).
\end{align}
For the pure advection problem, these jump terms are proportional to the differences in velocity sampled at different grid values as a Lagrangian marker crosses grid cell boundaries during a time step. Assuming the background velocity field is smooth, Taylor series analysis suggests that the jump terms should be proportional to the magnitude of the velocity's partial derivatives multiplied by the background grid's meshwidth. Thus, we expect the local truncation error for the forward Euler method applied to the pure advection problem is $\mathcal{O}\left(h\Delta t\right)$.\par
\begin{figure}[h]
    \centering
    \begin{minipage}[b]{0.49\textwidth}
     \centering
     \includegraphics[width=\textwidth]{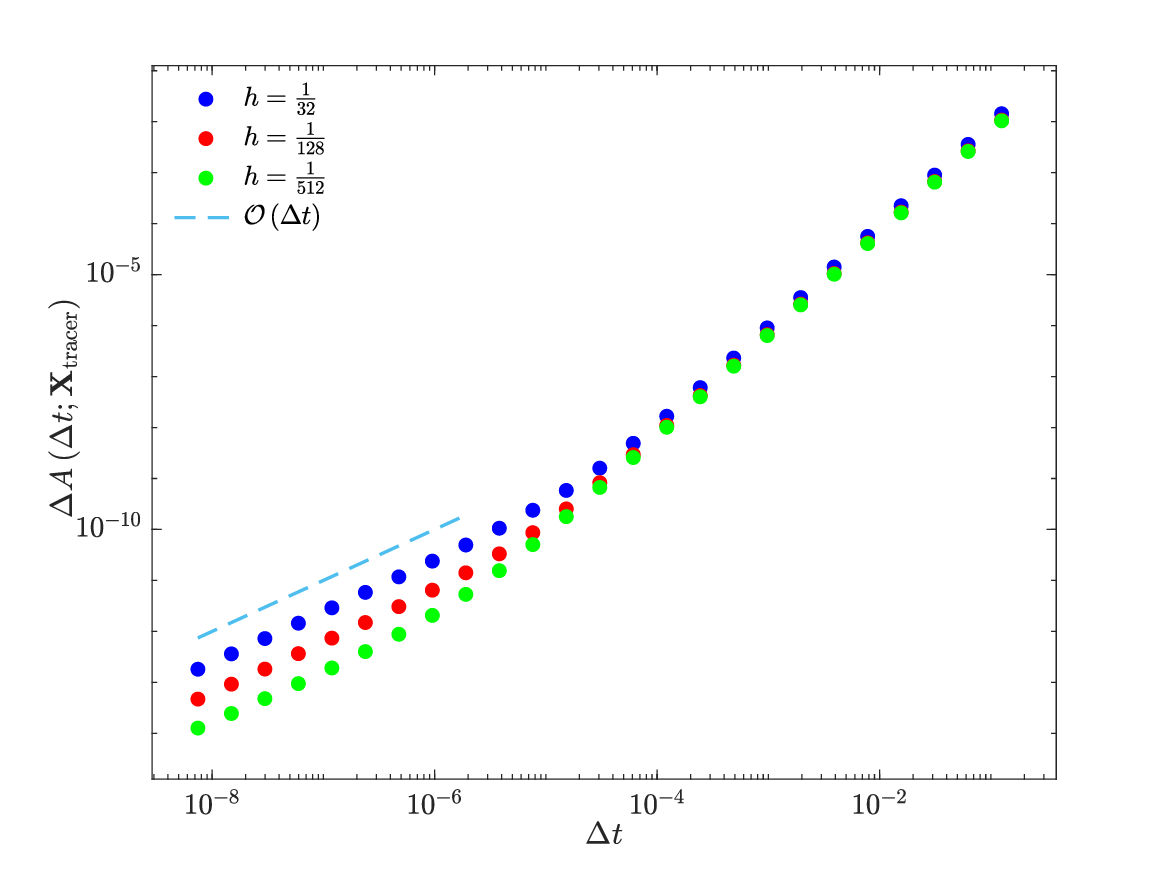}
     \caption*{(a) forward Euler}
    \end{minipage}
    \hfill
    \begin{minipage}[b]{0.49\textwidth}
     \centering
     \includegraphics[width=\textwidth]{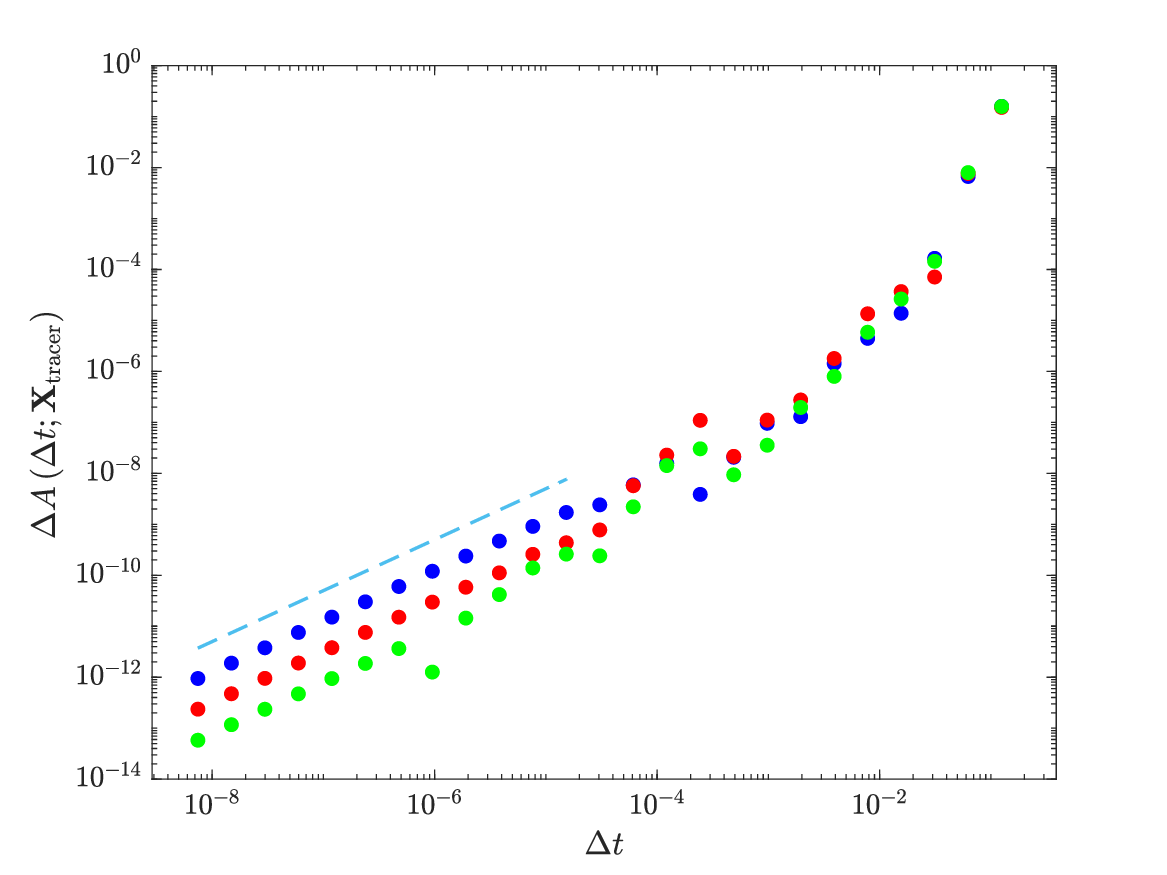}
     \caption*{(b) Midpoint}
    \end{minipage}
    \caption{Plots demonstrating the local truncation errors in the computed area of a circle advected using (a) the forward Euler method and (b) the explict midpoint rule.}
    \label{fig:LTE_comparisons}
\end{figure}

 
Our analysis focuses on the local truncation error of the forward Euler method, but the findings are applicable to any explicit time stepping scheme not adapted for discontinuities in the time derivative. Figure \ref{fig:LTE_comparisons} shows the error in the computed area after a single time step. In this test, Lagrangian markers were initially arranged to form a circle with radius $\frac{1}{4}$ centered at $\left(\frac{1}{2},\frac{1}{2}\right)$. These markers were then advected using the velocity field:
\begin{align}
    \mathbf{u} = \left(3\cos\left(4\pi\left(y-\frac{\pi}{4}\right)\right), 2\sin\left(2\pi\left(x-\frac{\pi}{4}\right)\right)\right).
\end{align}
For both the forward Euler and explicit midpoint rule methods, the local truncation error exhibits different convergence regimes depending on the choice of time step size $\Delta t$. At relatively larger time step sizes, the error is higher order and behaves roughly as one would expect when applying time integration schemes to problems without discontinuities in the time derivative. Specifically, the forward Euler method initially demonstrates apparent quadratic convergence, while the explicit midpoint rule shows roughly cubic convergence for larger time steps before quickly transitioning to a quadratic rate. Additionally, in this regime, the error in the computed area appears to be independent of the background grid spacing.
However, as anticipated by our analysis above, the asymptotic local truncation error rate is indeed first order for each method. For the smallest choices of time step size, the error appears to be linearly proportional to $h$ which aligns with our theoretical predictions based on the presence of discontinuities in the interpolated velocity as a result of the \BStwo\BSone~kernel.\par 
This contrast between the two regimes highlights that for relatively large choices of time step size, the local truncation error is roughly the expected rate for each method when applied to problems without discontinuities in the time derivative. However, as the time step becomes sufficiently small, the linear error rate associated with the discontinuity in the time derivative becomes dominant and the local truncation error associated with the area reduces to first order yielding a global error rate which is $\mathcal{O}(1)$. While the analysis provided here gives a reasonable explanation as to why the area error for the \BStwo\BSone~kernel in Fig \ref{fig:Advect_test} reduces to $\mathcal{O}(1)$ for the smaller choices of time step size, we note that there are also conceivably errors which are of similar size which are committed when fitting the cubic splines to the updates in the positions of the Lagrangian markers. However, empirical tests not reported herein have indicated that using splines to interpolate the positions of the Lagrangian markers advected by a discontinuous velocity field incorporates errors which are much smaller than the local truncation errors admitted by our choice of timestepping scheme.\par 
The analysis provided here makes explicit the critical role of the \BStwo\BSone~kernel in generating a discontinuous interpolated velocity field, which in turn yields discontinuities in the time derivative of the Lagrangian marker positions. The errors made by applying time-integration schemes naive to the presence of discontinuities in the time derivative are the root cause of first order local truncation errors which lead to $\mathcal{O}(1)$ global truncation errors in the area over the course of a simulation. Our findings emphasize the limitations of applying standard explicit timestepping methods to problems involving discontinuous interpolated velocity fields, such as those produced by the \BStwo\BSone~kernel. Future work might explore alternative interpolation methods or adaptive timestepping strategies that can better maintain higher-order accuracy in the presence of such discontinuities, potentially improving the overall accuracy of IB method simulations using the \BStwo\BSone~kernel.

\end{document}